\documentclass[11pt]{article}
\usepackage{amsmath,amssymb,amsthm}
\usepackage[T1]{fontenc}
\usepackage[latin9]{inputenc}
\usepackage[a4paper]{geometry}
\usepackage{setspace}
\usepackage{comment}
\usepackage[affil-it]{authblk}
\usepackage[english]{babel}
\usepackage{enumerate}
\usepackage{graphicx}
\usepackage{mathrsfs}
\usepackage{mathabx}
\usepackage{xcolor}
\usepackage{natbib}
\usepackage{url}
\usepackage{upgreek}
\usepackage{subcaption} 
\usepackage{caption}
\usepackage{bm}
\usepackage{natbib}
\usepackage[colorlinks = true,
            linkcolor = blue,
            urlcolor  = blue,
            citecolor = blue,
            anchorcolor = blue]{hyperref}
\usepackage{lineno}
 
\usepackage[ruled,vlined,linesnumbered]{algorithm2e}

\newtheorem{theorem}{Theorem}
\newtheorem{lemma}{Lemma}
\newtheorem{po}{Proposition}
\newtheorem{co}{Corollary}
\theoremstyle{definition}
\newtheorem{definition}{Definition}
\newtheorem{assumption}{Assumption}
\newtheorem{remark}{Remark}
\usepackage{amssymb}

\newcommand*{\p}{\mathbb{P}}
\newcommand{\T}{\top}
\newcommand{\RR}{\mathbb{R}}
\newcommand{\ee}{\end{aligned} \end{equation}}
\newcommand{\eq}{\end{quote}}
\newcommand{\diag}{\mathrm{diag}}
\newcommand{\sk}{\tilde{k}}
\newcommand{\ep}{\end{parts}}

\newcommand{\bqp}{\begin{quote}\begin{parts}}
\newcommand{\bds}{\boldsymbol}
\newcommand{\epq}{\end{parts}\end{quote}}
\newcommand{\sn}{S_n}
\newcommand{\en}{E_n}
\newcommand{\twoinf}{2,\infty}
\newcommand{\E}{\mathbb{E}}
\newcommand{\Eq}{Eq.~}

\newcommand{\Uk}{\M U}

\newcommand{\be}{\beta_{E}}

\newcommand{\mGs}{\bds{\M Q}^*}
\newcommand{\bd}{\beta_{\Delta}}
\DeclareMathOperator*{\argmin}{argmin}
\newcommand{\Rom}[1]{\text{\uppercase\expandafter{\romannumeral #1\relax}}}
\newcommand{\bee}{\begin{equation}\begin{aligned}}
\newcommand{\M}[1]{{{\mathbf{\MakeUppercase{#1}}}}}

\newcommand{\vertiii}[1]{{\vert\kern-0.25ex\vert\kern-0.25ex\vert #1 
    \vert\kern-0.25ex\vert\kern-0.25ex\vert}}
\newcommand{\emm}{\end{bmatrix}}

\newcommand{\F}{\mathrm{F}}

\newcommand{\argmax}{\operatornamewithlimits{argmax}}
\numberwithin{equation}{section}
\usepackage{babel} 

\let\hat\widehat
\let\tilde\widetilde
\let\check\widecheck

\begin{document}

 \title{\huge\textbf{Perturbation Analysis of Randomized SVD and its Applications  to Statistics}}
\author{Yichi Zhang\thanks{Department of Statistics, Indiana University Bloomington. Email:~\url{yiczhan@iu.edu}}$\ \ $and Minh Tang\thanks{Department of Statistics, North Carolina State University.}}\date{}
  \maketitle
 \begin{abstract} Randomized singular value decomposition (RSVD) is a
class of computationally efficient algorithms for computing the
truncated SVD of large data matrices.  Given an $m \times n$ matrix $\hat{\M
M}$, the prototypical RSVD algorithm outputs an approximation of the
$k$ leading left singular vectors of $\hat{\M M}$ by computing the SVD of
$\hat{\M M} (\hat{\M M}^{\top} \hat{\M M})^{g} \M G$; here $g \geq 1$ is an integer and
$\M G \in \mathbb{R}^{n \times \sk}$ is a random Gaussian sketching
matrix with $\sk \geq k$.  In this paper we derive upper bounds for
the $\ell_2$ and $\ell_{\twoinf}$ distances between the
exact left singular vectors $\hat{\M U}$ of $\hat{\M M}$ and its
{\em approximation} $\hat{\M U}_g$ (obtained via RSVD), as well as
entrywise error bounds when $\hat{\M M}$ is projected onto $\hat{\M
U}_g \hat{\M U}_g^{\top}$. These bounds depend on the singular values
gap and number of power iterations $g$, and smaller gap requires
larger values of $g$ to guarantee the convergences of the $\ell_2$ and $\ell_{\twoinf}$ distances.
We apply our theoretical results to
settings where $\hat{\mathbf{M}}$ is an additive perturbation of some
unobserved signal matrix $\mathbf{M}$. In particular, we obtain the nearly-optimal
convergence rate and asymptotic normality for RSVD on three inference problems, namely,
subspace estimation and community detection in random graphs, noisy
matrix completion, and PCA with missing data. 

 \bigskip{}
\end{abstract}

 \tableofcontents
 

 \doublespacing

\section{Introduction}
\label{sec:intro}

Spectral methods are popular in statistics and machine
learning as they provide simple algorithms with  
strong theoretical guarantees for a diverse number of inference problems including network analysis
\citep{rohe2011spectral}, matrix completion and denoising
\citep{achlioptas2007fast,chatterjee2015matrix}, covariance
estimation/principal component analysis (PCA), non-linear dimension reduction
and manifold learning \citep{belkin03:_laplac}, ranking
\citep{chen2019spectral}, etc. A common unifying theme for spectral algorithms is, given a $\hat{\M M}$ of dimensions $m \times n$, first compute a factorization of $\hat{\M M}$
via singular value decomposition (SVD), keep only the $k $ leading singular values and
singular vectors, and finally perform inference using the truncated SVD
representation. The value $k$ is usually chosen to
be as small as possible while still preserving most of the information.

For many inference problems in statistics, the observed matrix
$\hat{\M M}$ is generally noisy due to sampling and/or perturbation
errors, i.e., $\hat{\M M}$ is 
generated from a ``signal-plus-noise'' model
$\hat{\M M}  = \M M + \M E$
where $\M M$ is assumed to be the underlying true signal matrix 
with certain structure such as being (approximately) low rank and/or sparse, and $\M
E$ is the unobserved perturbation noise. 
Let $\hat{\M U}$ and $\M U$
be the leading singular vectors of $\hat{\M M}$ and $\M M$, respectively. 
As $\hat{\M M}$ is a noisy
realization of $\M M$,  $\hat{\M U}$ will also be 
a noisy estimate of $\M U$ and thus the main aim now is to bound the distance
between $\hat{\M U}$ and $\M U$ or between $\hat{\M U} \hat{\M U}^{\top} \hat{\M M}$ and
$\M M$. 

Error bounds for $\hat{\M U}$ is a fundamental
topic in matrix perturbation theory. Classical
results include the Davis--Kahan 
and Wedin's Theorems
\citep{davis1970rotation,wedin1972perturbation} for eigenvectors and singular vectors subspaces; these results make minimal
assumptions on $\M E$. 
The last decade has witnessed further study of matrix perturbations
from more statistical perspectives by introducing additional
assumptions on $\M E$ and $\M M$ such as (1) the entries of $\M E$
are independent random variables and/or (2) the leading singular vectors
of $\M M$ has bounded coherence. Examples include more refined
matrix concentration
inequalities \citep{tropp,oliveira},
rate-optimal subspace perturbation bound \citep{cai2018rate}, and $\ell_{\twoinf}$
perturbation bounds
\citep{abbe2020entrywise,damle,fan2018eigenvector,
eldridge2018unperturbed,mao2020estimating,cape2019signal, cape2019two,
lei2019unified}. In particular, error bounds for $\hat{\M U}$ in $\ell_{\twoinf}$ norm (see Section~\ref{sec:notation} for a definition) yields finer and more uniform control between $\M U$ and $\hat{\M U}$, and thus can be used to derive limiting distributions for the rows of $\hat{\M U}$ and entrywise confidence intervals for $\hat{\M U} \hat{\M U}^{\top} \hat{\M M} - \M M$;
see \cite{chen2021spectral} for a recent survey with several illustrative examples. 

While $\hat{\M U}$ has many desirable statistical properties, its computation can be quite challenging when the dimensions of $\hat{\M M}$ are large. Indeed, many classical algorithms for SVD, such as those based on pivotal QR
decompositions and/or Householder transformations, require
$O(mn\min\{m,n\})$ floating-point operations (flops) and return
the full set of $\min\{m,n\}$ singular values and vectors, even when
only the leading $k < m$ of them are desired; see Sections~5.4 and 8.3
of \cite{golub2013matrix}. 
These algorithms also require random access to the entries of $\hat{\M M}$
and are thus very inefficient when $\hat{\M M}$ is too large to
store in RAM due to the need for frequent data transfer between slow and fast memory.
Recently in the numerical linear algebra community, 
randomized SVD (RSVD) \citep{rokhlin2010randomized,
halko2011finding, musco2015randomized} had been widely
studied with the aim of providing fast, memory efficient, and accurate
approximations for the truncated SVD of large data matrices. 
The prototypical RSVD algorithm \citep{halko2011finding} first sketches $\hat{\M M}$
into a smaller matrix $\mathbf{Y} = \hat{\M M}(\hat{\M M}^{\top} \hat{\M M})^{g}
\M G\in\RR^{m\times
  \tilde{k}}$ where $\M G \in\RR^{n\times
  \tilde{k}}$ is a random matrix, and then uses the $k$ leading left singular
vectors of $\mathbf{Y}$, namely $\hat{\M U}_g$, as an approximation to $\hat{\M U}$. The parameters $g$ and $\sk$ are user-specified, with $g$ usually a
small integer and $\sk$ being slightly larger
than $k$. There are numerous choices for $\M G$ including
Gaussian, Rademacher, comlumn-subsampling, and random orthogonal matrices; see 
\cite{mahoney2011randomized,woodruff2014sketching,
  kannan2017randomized} and the references therein. 

The {\it sketch-and-solve} strategy of RSVD yields an algorithm with computational complexity of
$O(mn\sk)$ flops
and furthermore
these mainly involve the matrix-matrix products $\hat{\M M} \M X$ and $\hat{\M M}^{\top} \M X$, where $\M X$
is of dimensions $m \times \sk$ or $n \times \sk$, which are
highly-optimized operations on almost all computing platforms. 
RSVD is also ``pass efficient"
\citep{drineas2006fast, halko2011finding} and requires at most $(2g + 1)$
passes through the data; this dramatically reduces memory storage
\citep{golub2013matrix, lopes2020error}. 
Finally, RSVD allows for data compression \citep{cormode2011sketch}
and can be adapted to a streaming setting
\citep{tropp2019streaming}. Many recent works have replaced 
classical SVD with RSVD; see e.g.,
\cite{tsiligkaridis2013covariance,davenport2016overview,tsuyuzaki2020benchmarking,zhang2018billion,kumar2019target,hie2019efficient}
 for 
examples in covariance matrix estimation, matrix completion, and network
embeddings.

Existing theoretical results for RSVD, such as those in  
\cite{rokhlin2010randomized,halko2011finding,saibaba,lopes2020error}, 
focused exclusively on the setting where $\hat{\M M}$ is assumed to be
noise-free, i.e., these results either bound $\vertiii{\sin \Theta(\hat{\M U}_g, \hat{\M U})}$ or
$\vertiii{(\M I - \hat{\M U}_g \hat{\M U}_g^{\top}) \hat{\M M}}$ where $\vertiii{\cdot}$ denote some unitarily invariant (UI) norm. 
In particular these error bounds for $\hat{\M U}_g$ 
decreases as $g$ (the number of power iterations) and/or $\sk$
(the sketching dimensions) increases. 
However, if $\hat{\M M}$ is noisy then its leading singular vectors $\hat{\M U}$ will also
be a noisy estimate of $\M U$.  
Under this perspective, the main aim now should be to bound the difference between
$\hat{\M U}_g$ and $\M U$, and we thus need to balance
between the approximation error of $\hat{\M U}_g$ to $\hat{\M U}$, and the estimation
error of $\hat{\M U}$ to ${\M U}$.
This is straightforward and generally leads to sharp bounds when both approximation and estimation errors are 
in terms of UI norms. However, if results for $\hat{\M U}$ use
$\ell_{\twoinf}$ or entrywise norms while those for $\hat{\M U}_g$ use UI norms then
their combinations will surely be sub-optimal. In other words, 
more refined results for $\hat{\M U}$, such as row-wise limiting distributions of $\hat{\M U}$
and entrywise concentration for $\hat{\M U} \hat{\M U}^{\top} \hat{\M M} - \M M$, do not extend directly
to that for $\hat{\M U}_g$ and $\hat{\M U}_g \hat{\M U}_g^{\top} \hat{\M M} - \M M$. 

This paper addresses the above-mentioned gap, i.e.,
we derive, under minimal assumptions, upper bounds
for the $\ell_{\twoinf}$ difference
between $\hat{\M U}_g$
 and $\hat{\M U}$, as well as 
entrywise concentration for 
$(\hat{\M U}_g \hat{\M U}_g^{\top} - \hat{\M U} \hat{\M U}^{\top})
\hat{\M M}$. 
As a by-product of our analysis we also obtain bounds for $\vertiii{\sin \Theta(\hat{\M U}_g, \hat{\M U})}$
  comparable to those in the RSVD literature but with a markedly simple proof. 
  In the setting where $\hat{\M M}$ is generated under a ``signal-plus-noise''
  model, we show that the $\ell_2$ and 
  $\ell_{\twoinf}$ bounds between $\hat{\M U}_g$ and $\M U$
  exhibit a phase-transition phenomenon in that if the signal-to-noise ratio (SNR) decreases then $g$ need to
  increase to guarantee sharp convergence rates and {\color{black}asymptotic normality}. Precise values of $g$ where this transition occurs
  can also be determined    provided that $\hat{\M M}$ satisfy a
  certain trace growth conditions. 
  Finally we apply our
  theoretical results to three inference
  problems: subspace estimation and community detection in random graphs, matrix completion, and PCA with
  missing data. By combining our results for $\hat{\M U}_g$ with existing
  results for $\hat{\M U}$, we show that
  $\hat{\M U}_g$ has the same theoretical guarantees as $\hat{\M U}$, thus our results  provide a bridge between the numerical linear algebra and statistics communities. 

{\color{black}Empirically, the effectiveness of RSVD has been
demonstrated across a wide range of application domains, including
single-cell RNA sequencing (scRNA-seq)
\citep{hie2019efficient,tsuyuzaki2020benchmarking}, geophysical
imaging \citep{kumar2019target}, and large-scale network analysis
\citep{zhang2018billion}. Thus, the primary objective of our numerical
studies is not to extend the scope of RSVD's empirical use, but rather
to illustrate the theoretical insights developed in this
paper and showcase the statistical validation of RSVD in applications. We conduct simulations and real-data analyses using
RSVD for random graph inference, PCA, and matrix completion. These
empirical studies consistently support our theoretical predictions, namely that
to achieve sharp convergence rates and asymptotic normality under
certain SNR regimes requires sufficiently large $g$ and $\sk$. We
present the random graph simulations and scRNA-seq data analysis
in the main text; the remaining numerical results are included in
the Supplementary File. }
\subsection{Notation}
\label{sec:notation}
Let $a$ be a positive integer. We  write $[a]$ to denote the set
$\{1,\dots,a\}$. For two {\color{black}non-negative} sequences
$\{a_n\}_{n \geq 1}$ and $\{b_n\}_{n \geq 1}$, we write $a_n \precsim
b_n$, {\color{black}$b_n \succsim a_n$, $a_n = O(b_n)$ or $b_n = \Omega(a_n)$} if there exists a constant $c > 0$ not depending on $n$ such that
$a_n \leq c b_n$ for all but finitely many $n \geq 1$.  We write $a_n \asymp b_n$ if $a_n\precsim
b_n$ and $a_n\succsim b_n$. {\color{black}We   write $a_n = o(b_n)$ or $b_n = \omega(a_n)$ if $\lim_{n\rightarrow \infty}a_n/b_n = 0$.} The set of $d\times d'$ matrices with orthonormal columns is denoted as
$\mathbb{O}_{d\times d'}$ when $d \neq d'$ and is denoted as
$\mathbb{O}_d$ otherwise. 
Let $\M N$ be an arbitrary matrix. We 
denote the $i$th row of $\M N$ by $[\M N]_i$, and the $ij$th entry
of $\M N$ by $[\M N]_{ij}$. We write $\mathrm{tr} \, \M N$ and $\mathrm{rk}(\M N)$ to
denote the trace and rank of a matrix $\M N$, respectively,
and write $\sigma_k(\M N)$ as the $k$th largest singular value of $\M N$. 
The spectral and Frobenius norm of $\M N$
are denoted as $\|\mathbf{N}\|$ and $\|\mathbf{N}\|_{\F}$,
respectively. The maximum (in modulus) of the entries of $\M N$ is
denoted as $\|\M N\|_{\max}$. In addition we denote the $2 \to \infty$ norm of
$\mathbf{N}$ by
$\|\mathbf{N}\|_{\twoinf} = \max_{\|\bm{x}\| = 1}
\|\mathbf{N} \bm{x}\|_{\infty} = \max_{i} \|[\M N]_{i}\|,$
i.e., $\|\mathbf{N}\|_{\twoinf}$ is the maximum of the $\ell_2$
norms of the rows of $\M N$. We have the relationships
$n^{-1/2}\|\mathbf{N}\| \leq \|\mathbf{N}\|_{\twoinf} \leq \|\mathbf{N}\|$ and $
\|\M N\|_{\max} \leq \|\M N\|_{\twoinf} \leq d^{1/2} \|\M
 N\|_{\max},
$
where $n$ and $d$ are the number of rows and columns of
$\mathbf{N}$, respectively.
For two matrices $\M U_1 \in
\mathbb{O}_{n \times d}$ and $\M U_2\in \mathbb{O}_{n \times d}$, we define their
$\ell_2$ and $\ell_{\twoinf}$ distances as
\[
d_{2}(\M U_1, \M U_2) := \inf_{\M W\in \mathbb{O}_d}\|\M U_1- \M
U_2\M W\|, \quad 
d_{\twoinf}(\M U_1, \M U_2) := \inf_{\M W\in \mathbb{O}_d}\|\M U_1- \M U_2\M W\|_{\twoinf}.
\]
Note that $d_2(\M U_1, \M U_2) \leq \sqrt{2} \|\sin \Theta(\M U_1, \M
U_2)\|$ where $\sin \Theta(\M U_1, \M U_2)$ is the diagonal matrix whose elements are the singular values of $(\M I - \M U_1 \M U_1^{\top}) \M U_2$. 

\subsection{Randomized SVD}
Let $\hat{\mathbf{M}}$ be a $m \times n$ matrix and suppose that we want to compute the singular vectors
$\hat{\mathbf{U}}^{(k)}$ associated with the $k$ largest singular values of
$\hat{\mathbf{M}}$. One popular and widely
used approach for computing $\hat{\mathbf{U}}^{(k)}$ is via randomized
subspace iteration. More specifically we first sample an $n \times\sk$ matrix $\M G$ whose entries are iid standard normals and compute $\mathbf{Y}_g =
\hat{\mathbf{M}}^{g} \M G$ if $\hat{\M M}$ is symmetric and $\M Y_g = \hat{\M M} (\hat{\M M}^{\top} \hat{\M M})^{g} \M G$ otherwise, where $g$ is a positive integer. 
Let $\hat{\M U}_g^{(k)}$ be the $n \times k$ matrix whose columns form an
orthonormal basis for the $k$ leading left singular vectors of $\mathbf{Y}_g$.
Then $\hat{\M U}_g^{(k)}$ is an approximation to $\hat{\mathbf{U}}^{(k)}$
and we can take $\hat{\M U}_g^{(k)} \hat{\M U}_g^{(k)\top} \hat{\mathbf{M}}$ as a low rank approximation to $\hat{\mathbf{M}}$; see Algorithm~\ref{RSRS} for a formal description. The value of $\sk$, the number of columns of $\M G$, is often chosen to be slightly larger than $k$ in order to increase the probability that the column space of
$\mathbf{Y}_g$ is closely aligned with
$\hat{\mathbf{U}}^{(k)}$, and empirical observations suggest that
$5 \leq \sk - k \leq 10$ is sufficient for most practical applications \citep{halko2011finding}.
Algorithm~\ref{RSRS} is algebraically equivalent
to a version wherein one periodically orthonormalizes 
$\hat{\M M} (\hat{\M M}^{\top} \hat{\M M})^{g'} \M G$ (via QR decomposition) for $g' < g$ before computing 
$\hat{\M M} (\hat{\M M}^{\top} \hat{\M M})^{g'+1} \M G$; see e.g., Remark~4.3 of \cite{halko2011finding}. This extra orthonormalization
leads to more numerically stable outputs but has no impact on the theoretical results. 
For more discussion on randomized subspace iteration, see Section~4.5 of
\cite{halko2011finding}, Section~11.6 of \cite{martinsson_rnla}, and Section~4.3 of
\cite{woodruff2014sketching}, \cite{musco2015randomized}. 
Algorithm~\ref{RSRS} is known as the \texttt{PowerRangeFinder} and \texttt{SubspacePowerMethod} in \cite{halko2011finding,woodruff2014sketching,musco2015randomized}. 
\begin{algorithm}[tp]
\KwIn{$\hat{\M M}\in\RR^{m\times n}$, rank $k\geq 1$, sketching
  dimension $\sk \geq k$, power iterations $g \geq 1$.}
 Generate a $n \times k$ sketching matrix $\M G$ whose elements are iid standard normals\;
 If $\hat{\M M}$ is symmetric, compute $\hat{\M M}^{g} \M G$ by
iterating $\hat{\M M} \M G, \hat{\M M} (\hat{\M M}\M G), \dots, \hat{\M M} (\hat{\M M}^{g-1} \M G)$, otherwise, compute $\hat{\M M} (\hat{\M M}^{\top} \hat{\M M})^{g} \hat{\M G}$ by iterating $\hat{\M M} \M G, \hat{\M M}^{\top}(\hat{\M M} \M G), \dots,\hat{\M M} (\hat{\M M}^{\top} \hat{\M M})^{g} \hat{\M G}$\;  
Obtain the {\em exact} SVD of either $\hat{\M M}^{g} \M G$ or
$\hat{\M M} (\hat{\M M}^{\top} \hat{\M M})^{g} \M G$ and
let $\hat{\M U}_g^{(k)}$ be the $m\times k $ matrix whose columns are the $k$ leading left singular vectors\;
\KwOut{ Estimated singular vectors $\hat{\M U}_{g}^{(k)}$ and low-rank $\hat{\M U}_g^{(k)}\hat{\M U}_g^{(k)\T}\hat{\M M}$.}
 \caption{RSVD}\label{RSRS}
\end{algorithm}

\begin{remark}
If we set $g = 1$ in Algorithm~\ref{RSRS} then we
get the ``sketched SVD'' algorithm described in
\cite{lopes2020error,mahoney2011randomized}. Sketched
SVD is very useful when $\hat{\M M}$ is too large
to store in fast memory as the procedure only requires one pass
through the data. However, as we will show in Section~\ref{sec:thm} and
Section~\ref{sec:rgi}, 
setting $g = 1$ can lead to poor estimates of $\M U$ unless $\sk = \Omega(n) \gg k$. 
The choice $\sk = \Omega(n)$ has recently been
considered in the context of sketching PCA \citep{yang2021reduce}. 
However, in practice it is preferable to
choose $\sk$ as small as possible, because the Step~3 in Algorithm~\ref{RSRS} requires
$O(m \sk^2)$ flops.
\end{remark}

\section{Theoretical results}\label{sec:thm}
Let $\hat{\M M}$ be a $m \times n$ matrix and for any $k \leq \min\{m,n\}$ 
denote the SVD of $\hat{\M M}$ by 
\bee\label{decom:M:hatM}
\hat{\M M} :=  \hat{\M U}^{(k)} \hat{\bm{\Sigma}}^{(k)} \hat{\M V}^{(k)\top} + \hat{\M U}^{(k)}_{\perp} \hat{\bm{\Sigma}}^{(k)}_{\perp} \hat{\M V}^{(k)\top}_{\perp}
\ee
where $\hat{\bm{\Sigma}}^{(k)}$
is the diagonal matrix containing the $k$ largest singular values of $\hat{\M M}$, $\hat{\M U}^{(k)} \in \mathbb{R}^{m \times k}$ and $\hat{\M V}^{(k)} \in \mathbb{R}^{n \times k}$ are the corresponding left  and right singular vectors. 
We now present the general upper bounds for
$d_{2}(\hat{\M U}^{(k)}_g, \hat{\M U}^{(k)})$ and $d_{\twoinf}(\hat{\M U}^{(k)}_g, \hat{\M U}^{(k)})$.  

\begin{theorem}\label{thm2ss} 
Let $\hat{\M M}$ be given and compute $\hat{\M
U}^{(k)}_g$ via Algorithm \ref{RSRS} for some choices of $k, g$
and $\sk$ where $n \geq \sk \geq (1 - c_{\mathrm{gap}})^{-2}\{k + (8 k \log (1/\vartheta))^{1/2} + 2 \log (1/\vartheta) \}$;
here $c_{\mathrm{gap}} \in (0, 1)$ and $\vartheta > 0$ are both arbitrary. 
Denote by $\hat{\sigma}_i$ the $i$th largest singular value of $\hat{\M M}$ and let $\vertiii{\cdot}$ be any unitarily invariant norm. 
Then for all $g \geq 1$ we have
\begin{equation}
   \label{eq:ui_corollary_bound}
      \vertiii{\sin \Theta(\hat{\M U}_g^{(k)}, \hat{\M U}^{(k)})} = \vertiii{(\M I - \hat{\M U}_{g}^{(k)} \hat{\M U}_g^{(k)\top}) \hat{\M U}^{(k)}} \leq
      \frac{3 n^{1/2} \vertiii{\bigl(\hat{\bm{\Sigma}}_{\perp}^{(k)}\bigr)^{\tilde{g}}}}{c_{\mathrm{gap}} \sk^{1/2} \hat{\sigma}_{k}^{\tilde{g}} }
  \end{equation}
  with probability at least $1 - \vartheta - 2e^{-n/2}$, where $\tilde{g} = g$ if $\hat{\M M}$ is symmetric and $\tilde{g} = 2g + 1$ otherwise.
 Let $\hat{\zeta}_k = \hat{\sigma}_{k+1}/\hat{\sigma}_{k}$. Eq.~\eqref{eq:ui_corollary_bound} implies
  \begin{equation}
    \begin{split}
   \label{eq:d2_corollary_bound}
      d_2(\hat{\M U}_g^{(k)}, \hat{\M U}^{(k)}) \leq \sqrt{2} \|\sin \Theta(\hat{\M U}_g^{(k)}, \hat{\M U}^{(k)})\| \leq
      \frac{3\sqrt{2} n^{1/2}}{c_{\mathrm{gap}}\sk^{1/2}}
     \hat{\zeta}_k^{\tilde{g}}
      \end{split}
  \end{equation}
  with probability at least $1 - \vartheta - 2e^{-n/2}$. 
\end{theorem}
\begin{theorem}\label{thm3ss:general}
Consider the setting in Theorem~\ref{thm2ss}.
Let $\hat{\zeta}_k = \hat{\sigma}_{k+1}/\hat{\sigma}_{k}$. For any $\delta > 0$ such that ${\color{black}\tilde{k} \geq 2 \log \delta^{-1}}$, define
\begin{gather}
  \label{eq:r_twoinf}
  r_{\twoinf} =  \frac{\sqrt{128} e (k \log \delta^{-1})^{1/2} \hat{\zeta}_k^{\tilde{g}}}{c_{\mathrm{gap}}^2 \sk^{1/2}}  +
\frac{18 n  \|\hat{\M U}^{(k)}\|_{\twoinf}\hat{\zeta}_k^{2\tilde{g}}
}{c_{\mathrm{gap}}^2 \sk} + \frac{36 n
  (\log \delta^{-1})^{1/2} \hat{\zeta}_k^{3\tilde{g}}}{c_{\mathrm{gap}}^3 \sk}  
\end{gather}
where $\tilde{g} = g$ if $\hat{\M M}$ is symmetric and $\tilde{g} = 2g+1$ otherwise.   
  Then for all $g \geq 1$ we have
  \begin{gather}
    \label{eq:d2inf_corollary_bound}
    d_{\twoinf}(\hat{\M U}_g^{(k)}, \hat{\M U}^{(k)}) \leq 
    r_{\twoinf}
    \end{gather}
    with probability at least $1 - 4m \sk \delta - \vartheta - 2e^{-n/2}$, where $\vartheta$ appears
    in the lower bound condition for $\sk$ given in Theorem~\ref{thm2ss}. 
 Furthermore, for any $\gamma > 0$ such that $\tilde{k} \geq 2 \log \gamma^{-1}$, define
 \begin{equation}
   \label{eq:tilder2_inf}
    \tilde{r}_{\twoinf} =  \frac{\sqrt{128} e (k \log \gamma^{-1})^{1/2} \hat{\zeta}_k^{\tilde{g}}}{c_{\mathrm{gap}}^2 \sk^{1/2}}  +
   \frac{18 n  \|\hat{\M V}^{(k)}\|_{\twoinf}\hat{\zeta}_k^{2\tilde{g}}}{c_{\mathrm{gap}}^2 \sk} + \frac{36 n (\log \gamma^{-1})^{1/2} \hat{\zeta}_k^{3\tilde{g}}}{c_{\mathrm{gap}}^3 \sk},
 \end{equation}
 where $\tilde{g} = g$ if $\hat{\M M}$ is symmetric and $\tilde{g} = 2g+1$ otherwise. 
  Then for all $g \geq 1$, we have
  \begin{gather}
\label{eq:max_corollary_bound}
\|(\hat{\M U}_g^{(k)} \hat{\M U}_g^{(k)\top} - \hat{\M U}^{(k)} \hat{\M U}^{(k)\top}) \hat{\M M}\|_{\max}
   \leq \hat{\sigma}_1 (r_{\twoinf} \tilde{r}_{\twoinf}  + \|\hat{\M U}^{(k)}\|_{\twoinf} \tilde{r}_{\twoinf} + \|\hat{\M V}^{(k)}\|_{\twoinf} r_{\twoinf}),
   \end{gather}
  with probability at least $1 - 4\sk(m \delta + n \gamma) - \vartheta - 2e^{-n/2}$.
\end{theorem}
\begin{remark}[Technical ideas behind the $\ell_{\twoinf}$ bound]\label{rk:twoinf}\color{black}We briefly discuss the main technical ideas behind the
  \(\ell_{\twoinf}\) bound in Theorem~\ref{thm3ss:general}. For clarity we focus on the case where \(\hat{\M M}\) is an \(n\times n\)
  symmetric matrix with eigendecomposition \(\hat{\M U}\hat{\M\Lambda}\hat{\M U}^\T + \hat{\M U}_\perp\hat{\M\Lambda}_\perp\hat{\M U}_\perp^\T\). Note that
  the eigenvalues in \(\hat{\M\Lambda}\) and \(\hat{\M\Lambda}_\perp\) are ordered by decreasing magnitude, and the superscript ``\((k)\)'' is omitted for ease of notations.
  If $\sk - k > 0$ then (with probability one) the column space of \(\hat{\M U}\) coincides with that of \(\check{\M U}_g\) 
  whose columns are the $k$ leading left singular vectors of \(\hat{\M U} \hat{\bm\Lambda}^g \hat{\M U}^\top \M G\), and
  bounding \(d_{\twoinf}(\hat{\M U}_g,\hat{\M U})\) 
  reduces to bounding \(d_{\twoinf}(\check{\M U}_g,\hat{\M U})\). Next note that
\bee\label{decom:l2inf}
\hat{\M M}^{g} \M G = \hat{\M U} \hat{\bm\Lambda}^{g} \hat{\M U}^\top \M G + \hat{\M U}_\perp \hat{\bm\Lambda}_\perp^{g} \hat{\M U}_\perp^\top \M G,
\ee
where the second term on the right-hand side of Eq.~\eqref{decom:l2inf} can be viewed as an additive perturbation. 
Unlike standard settings, the matrices in Eq.~\eqref{decom:l2inf} are highly unbalanced when \(\sk \ll n\) and also strongly correlated (due to their dependency on a common \(\M G\)).
We address these issues as follows. First we
extend the deterministic Procrustes analysis in \citep{cape2019two} to
decompose the difference between \(\check{\M U}_g\) and \(\hat{\M U}_g\) into three terms \(\M T_1\), \(\M T_2\), and \(\M T_3\)
(see Eq.~\eqref{eq:Procrustean_1} in the Supplementary File). 
Next
we observe that \(\M T_1\) can be represented as the product of two {\em  independent} Gaussian matrices depending on $\M G$ from which, after some careful analysis
that exploit various properties of $\M G$, we obtain
a sharp concentration bound for \(\|\M T_1\|_{\twoinf}\).
Finally we bound $\|\M T_2\|_{\twoinf}$ and $\|\M T_3\|_{\twoinf}$ using a version of the Wedin $\sin$-$\Theta$ theorem for unbalanced matrices \citep{cai2018rate}. 
\end{remark}
\begin{remark}[Illustrations of the $r_{\twoinf}$ rate]\label{rk:illustrate}{\color{black}We simplify the \(\ell_{\twoinf}\) bound in~\eqref{eq:d2inf_corollary_bound} under different parameter regimes for illustration. Assume that \(\hat{\M M}\) is symmetric, \(k\) is fixed, \(\hat{\zeta}_k \asymp n^{-\psi}\), \(\sk \asymp n^{\phi}\) for some   \(\psi > 0\) and \(\phi \in (0,1]\), and \(\|\hat{\M U}^{(k)}\|_{\twoinf} \asymp n^{-1/2}\). Then~\eqref{eq:d2inf_corollary_bound} yields
\begin{gather}
\label{eq:d2inf_corollary_bound:2}
d_{\twoinf}(\hat{\M U}_g^{(k)}, \hat{\M U}^{(k)}) \precsim \sqrt{\log n} \cdot n^{-\psi g - \phi/2} + n^{1/2-2\psi g-\phi} + \sqrt{\log n} \cdot n^{1-3\psi g-\phi}
\end{gather}
with high probability. 
The dominant term in~\eqref{eq:d2inf_corollary_bound:2} depends on the
relative magnitudes of \(\phi\), \(\psi \) and $g$. For \(\phi \approx 0\), the first term dominates when
\(\psi g\) is large, i.e., the singular value gap and the number of
power iterations are relatively large, the second term dominates when
\(\psi g\) is moderate, and the third term dominates when \(\psi g\)
is small. For \(\phi \approx 1\), the second
term is always negligible compared to the first term, and furthermore the first term dominates for large \(\psi g\) and the third term dominates for small \(\psi g\).
Finally, for any choice of $\phi$, we can make 
$d_{\twoinf}(\hat{\M U}_g^{(k)},\hat{\M U}^{(k)})$ arbitrarily small by increasing $\psi$, or $g$, or both. 
See Section~\ref{sec:add:rg} of the Supplementary
File for concrete examples in the context of random graph inference.}
\end{remark}

  The choice of $c_{\mathrm{gap}}$ in Theorem~\ref{thm2ss} and Theorem~\ref{thm3ss:general} provides a trade-off between the magnitudes of $\sk - k$ and upper bounds for $\hat{\M U}_g^{(k)} - \hat{\M U}^{(k)}$ in $\ell_2$ and $\ell_{\twoinf}$ norms, i.e., if $c_{\mathrm{gap}}$ decreases then the ``oversampling'' dimension $\sk - k$ can be made smaller
but the upper bounds in Eq.~\eqref{eq:d2_corollary_bound}, Eq.~\eqref{eq:d2inf_corollary_bound}, and Eq.~\eqref{eq:max_corollary_bound} will increase. In our subsequent discussion we will, for ease of exposition, typically set $c_{\mathrm{gap}}$ to an arbitrarily small constant and thus
omit factors depending on $c_{\mathrm{gap}}$ from our bounds.
The condition for $\sk$ in
  Theorem~\ref{thm2ss} and Theorem~\ref{thm3ss:general} also depend on $\vartheta$, and this allows us to handle settings where $m$ and $n$ are of very different magnitudes.  
  For example, suppose $n \asymp e^{m^{\alpha}}$ for some $\alpha > 0$. 
  Then $\log n \gg \log m$ and thus the requirement for $\sk$ when setting $\vartheta = m^{-1}$
  is much less stringent compared to when setting $\vartheta = n^{-1}$. In the same vein,
  the bounds in Theorem~\ref{thm3ss:general} include two additional parameters $\delta$ and $\gamma$
  so that we can precisely control the magnitude of $r_{\twoinf}$ and $\tilde{r}_{\twoinf}$
  in Eq.~\eqref{eq:d2inf_corollary_bound} and Eq.~\eqref{eq:max_corollary_bound}, e.g.,
  if $m$ and $n$ are of very different magnitude then $\log \delta^{\color{black}-1} $ and $\log \gamma^{\color{black}-1} $
  can also be of different magnitude,
  while if $m \asymp n$ then we can choose $\log (\delta^{\color{black}-1}) = \log (\gamma^{\color{black}-1})  = \log (\vartheta^{\color{black}-1})  = (c+2) \log(m+n)$
  to guarantee that all of the
  bounds in Theorem~\ref{thm2ss} and Theorem~\ref{thm3ss:general} hold
  with probability at least $1 - O((m+n)^{-c})$ where $c > 0$ is any arbitrary constant. 
  Finally, for concreteness we had include
  explicit constants in our statement of Theorem~\ref{thm2ss} and Theorem~\ref{thm3ss:general} but
  their values are chosen mainly for ease of exposition and are thus possibly sub-optimal.

\begin{remark}
  For conciseness we only stated Theorem~\ref{thm2ss} and Theorem~\ref{thm3ss:general} for the approximate left singular vectors
  $\hat{\M U}_g^{(k)}$. Results for the approximate right singular vectors $\hat{\M V}_{g}^{(k)}$
  are obtained simply by applying the  algorithms and theorems to $\hat{\M M}^{\top}$, i.e., we replace
  $n$ with $m$ in Theorem~\ref{thm2ss} and swap the roles of $m$ and $n$ (as well as the roles of $\hat{\M U}^{(k)}$ and $\hat{\M V}^{(k)}$) in Theorem~\ref{thm3ss:general}. For example, we have
  \begin{equation}
    \label{eq:d2_hatV}
    d_2(\hat{\M V}_g^{(k)}, \hat{\M V}^{(k)}) \leq \frac{3 \sqrt{2} m^{1/2}}{c_{\mathrm{gap}} \sk^{1/2}} \Bigl(\frac{\hat{\sigma}_{k+1}}{\hat{\sigma}_{k}}\Bigr)^{\tilde{g}}
  \end{equation}
  with probability at least $1 - \vartheta - 2e^{-m/2}$ and
 \begin{equation}
   \label{eq:d2inf_hatV}
   d_{\twoinf}(\hat{\M V}_g^{(k)}, \hat{\M V}^{(k)}) \leq
 \frac{\sqrt{128} e (k \log \delta^{-1})^{1/2} \hat{\zeta}_k^{\tilde{g}}}{c_{\mathrm{gap}}^2 \sk^{1/2}}  +
\frac{18 m \hat{\zeta}_k^{2\tilde{g}} \|\hat{\M V}^{(k)}\|_{\twoinf}}{c_{\mathrm{gap}}^2 \sk} + \frac{36 m
  (\log \delta^{-1})^{1/2} \hat{\zeta}_k^{3\tilde{g}}}{c_{\mathrm{gap}}^3 \sk} 
 \end{equation}
 with probability at least $1 - 4n\sk \delta - \vartheta - 2e^{-m/2}$. 
\end{remark}

\begin{remark}
  \label{rk:smd}
  In this paper we only considered Gaussian sketching matrices as they are (1) most commonly used and
  (2) lead to simple and precise theoretical results. 
Other type of $\M G$ such as uniform and Rademacher have also been
studied in the literature \citep{mahoney2011randomized,woodruff2014sketching,kannan2017randomized}.
We note that any distribution for $\M G$ that satisfies Lemma~\ref{lm:basic} in the Supplementary,
will also lead to the same $\sin \Theta$ upper bounds (up to some multiplicative constants)
as those presented in Theorem~\ref{thm2ss}. In contrast, the analysis in  
Theorem~\ref{thm3ss:general} leverage several properties that are
intrinsic to normal random variables. We thus leave the extension of Theorem~\ref{thm3ss:general} for
non-Gaussian $\M G$ to future work. 
\end{remark}

Finally, let $\M M$ be a $m \times n$ matrix and, similar to Eq.\eqref{decom:M:hatM}, 
denote its SVD by
\[ \M M :=  \M U^{(k)} \bm{\Sigma}^{(k)} \M V^{(k)\top} + \M U^{(k)}_{\perp} \bm{\Sigma}^{(k)}_{\perp} \M V^{(k)\top}_{\perp},
\]
where $\bm{\Sigma}^{(k)}$
is the diagonal matrix containing the $k$ largest singular values of $\M M$ and $\M U^{(k)}, \M V^{(k)}$
are the corresponding left and singular vectors. Next suppose we only observe
$\hat{\M M} = \M M + \M E$ and want to use $\hat{\M U}_g^{(k)}$ as an estimate for $\M U^{(k)}$.  Then, by Weyl's inequality, $|\hat{\sigma}_i - \sigma_i| \leq \|\M E\|$
for all $i$, where $\sigma_1 \geq \sigma_2 \geq \dots$ are the singular values of $\M M$. Thus if
$\sigma_{k} > \|\M E\|$
then $\hat{\zeta}_k = \hat{\sigma}_{k+1}/\hat{\sigma}_{k} \leq (\sigma_{k+1} + \|\M E\|)/(\sigma_{k} - \|\M E\|)$.
Substituting this bound for $\hat{\zeta}_k$ into 
Theorem~\ref{thm2ss} and Theorem~\ref{thm3ss:general} we directly obtain upper bounds for
$d_2(\hat{\M U}_g^{(k)}, \M U^{(k)})$ and $d_{\twoinf}(\hat{\M U}_g^{(k)}, \M U^{(k)})$ that
depend only on $\|\M E\|$ and $\{\sigma_{k}, \sigma_{k+1}\}$. 
\begin{co}
  \label{co:arbitrary_noise}
  Consider the setting in Theorem~\ref{thm2ss} and suppose that $\hat{\M M} = \M M + \M E$ where $k$ is chosen
  such that $\sigma_{k} > \|\M E\| $. 
  Let 
$$
\zeta_k = (\sigma_{k+1} + \|\M E\|)/(\sigma_{k} - \|\M E\|),
$$ and define $r'_{\twoinf}$ and $\tilde{r}'_{\twoinf}$
  as $r_{\twoinf}$ and $\tilde{r}_{\twoinf}$ in Eq.~\eqref{eq:r_twoinf} and Eq.~(\ref{eq:tilder2_inf})
  but with $\zeta_k$, $u_{\twoinf} := \|\M U^{(k)}\|_{\twoinf} + d_{\twoinf}(\hat{\M U}^{(k)}, \M U^{(k)})$, and $v_{\twoinf} := \|\M V^{(k)}\|_{\twoinf} + d_{\twoinf}(\hat{\M V}^{(k)}, \M V^{(k)})$ in place of $\hat{\zeta}_k$,
  $\|\hat{\M U}^{(k)}\|_{\twoinf}$, and $\|\hat{\M V}^{(k)}\|_{\twoinf}$ respectively. 
  Then for all $g \geq 1$ we have
  \begin{equation}\nonumber
    \bigl|d_2(\hat{\M U}_g^{(k)}, \M U^{(k)}) - d_{2}(\hat{\M U}^{(k)}, \M U^{(k)})\bigr| \leq {\color{black}  d_2(\hat{\M U}_g^{(k)}, \hat{\M U}^{(k)})}\leq3 \sqrt{2} c_{\mathrm{gap}}^{-1}
    (n/\sk)^{1/2} \zeta_k^{\tilde{g}}
  \end{equation}
  with probability at least $1 - \vartheta - 2e^{-n/2}$, and
  \begin{equation}\nonumber
    \bigl|d_{\twoinf}(\hat{\M U}_g^{(k)}, \M U^{(k)}) - d_{\twoinf}(\hat{\M U}^{(k)}, \M U^{(k)})\bigr|\leq {\color{black}  d_{\twoinf}(\hat{\M U}_g^{(k)}, \hat{\M U}^{(k)})} \leq r'_{\twoinf}, 
  \end{equation}
  with probability at least $1 - 4m \sk \delta - \vartheta - 2e^{-n/2}$.
  Finally, denote $\hat{\bm{\Pi}}_g^{(k)} = \hat{\M U}_g^{(k)} \hat{\M U}_g^{(k)\top}$ and $\hat{\bm{\Pi}}^{(k)} = \hat{\M U}^{(k)} \hat{\M U}^{(k)\top}$. Then we also have
  \bee
    \label{eq:entrywise_co_general}
    \bigl|\|\hat{\bm{\Pi}}_g^{(k)} \hat{\M M} - \M M\|_{\max} -
    \|\hat{\bm{\Pi}}^{(k)} \hat{\M M} - \M M\|_{\max} \bigr|
    \leq (\sigma_{1} + \|\M E\|) \left(r'_{\twoinf} \tilde{r}'_{\twoinf} + 
    u_{\twoinf} \tilde{r}'_{\twoinf} 
    + v_{\twoinf} r'_{\twoinf}\right)
  \ee
  with probability at least $1 - 4\sk(m \delta + n \gamma) - \vartheta - 2e^{-n/2}$.
  In all of the above bounds we have $\tilde{g} = g$ if $\hat{\M M}$ is symmetric and $\tilde{g} = 2g+1$ otherwise. 
\end{co}

\subsection{Refined bounds for low-rank setting}
Theorems~\ref{thm2ss} and Theorem~\ref{thm3ss:general} 
hold for any $\hat{\M M}$, any $k$,
and any $g$.
In particular they implied that if $\hat{\M M} = \M M + \M E$ and
there exists a $c > 1$ such that $1/\zeta_k =  (\sigma_k - \|\M E\|)/(\sigma_{k+1} + \|\M E\|) \geq c$ 
then the upper bounds in Corollary~\ref{co:arbitrary_noise} converge to $0$ at rate $c^{-g}$, 
and are thus negligible for $g = \Omega(\log(n)/\log(c))$. 
If $\M M$ is low-rank with $\mathrm{rk}(\M M) = k_0$, then the bounds in Corollary~\ref{co:arbitrary_noise} can be further refined when we choose $k = k_0$, especially when $\|\M E\| \ll \sigma_{k_0}$. 
To reduce notational burden we only present results when
$\hat{\M M}$ is symmetric (this is
sufficient for our discussion in Section~\ref{sec:rgi})
so that the subsequent bounds and probabilities can be stated explicitly in terms of $n, \sk$ and $\sigma_{k_0}/\|\M E\|$. Similar results hold when $\hat{\M M}$ is asymmetric/rectangular and we leave them to the interested reader. 
\begin{co}\label{co:l2_noise}
  Assume the setting of Corollary~\ref{co:arbitrary_noise} where $\hat{\M M}$ is a symmetric $n \times n$ matrix
  with $\mathrm{rk}(\M M) = k_0$. Denote $E_n = \|\M E\|$ and
   let $\hat{\M U}^{(k_0)}_{g}$ be computed from $\hat{\M M}$ via Algorithm \ref{RSRS} for 
   some choice of $\sk \geq (1 - c_{\mathrm{gap}})^{-2} (k_0 + \sqrt{24 k_0 \log n} + 6 \log n)$ where
   $c_{\mathrm{gap}} \in (0,1)$ is a fixed but arbitrary constant. For any specified $\iota > 0$, let $g_\iota = \frac{\log\{3\sqrt{2}c_{\mathrm{gap}}^{-1}({n/\sk})^{1/2} \iota^{-1}\} }{ \log(1/\zeta_{k_0})}.$ 
   If $\sigma_{k_0} > 2E_n$ then we have
   \begin{equation}
   \label{eq:co_main_ws}
       d_{2}(\hat{\M U}^{(k_0)}_g, \M U^{(k_0)}) \leq 
\big(1 + \iota \cdot \zeta_{k_0}^{g   - g_\iota - 1} \big)\frac{E_n}{\sigma_{k_0}}, \text{ for any }g \geq 1, 
    \end{equation}
    with probability at least $1 - 2n^{-3}$; recall $\zeta_{k_0} = E_n/(\sigma_{k_0} -E_n) < 1$.
 Furthermore, if $\sigma_{k_0}/E_n \succsim n^{\epsilon}$ for some fixed but arbitrary $\epsilon > 0$ then
    Eq.~\eqref{eq:co_main_ws} can be refined  to
\bee\label{thm2*:resopt:main}
d_2(\hat{\M U}_{g}^{(k_0)},\M U^{(k_0)}) \leq  \begin{cases} 
 O\bigl((E_n/\sigma_{k_0})^{g - g_*}\bigr) & \text{if $g_* \leq g 
 \leq  1 + g_*$} \\
(1 + o(1))E_n/\sigma_{k_0}  & \text{if  $g \geq 1 + g_* + \Delta$}
\end{cases},
\ee
with probability at least $1 - 2n^{-3}$. Here $g_* =  \frac{\log(n/\sk)}{2 \log(1/\zeta_{k_0})}$ and $\Delta = \omega(\log^{-1}{n})$ is arbitrary; note that
$g_*\leq (2\epsilon)^{-1}$ for sufficiently large $n$. 
 
\end{co}

\begin{remark}\label{rk:tc:thm2}
We now discuss two special cases of Corollary~\ref{co:l2_noise}. 
  \begin{enumerate}
  \item 
    Suppose $cn \leq \sk \leq (1 - c)n$ for some constant 
    $c \in (0,1)$ and $\sigma_{k_0} =  \omega(E_n)$. 
    Then $\frac{\log(n/\sk)}{2 \log(1/\zeta_{k_0})} = o(1)$, and setting $\iota = 4\sqrt{2}c_{\mathrm{gap}}^{-1}c^{-1/2}$ in \eqref{eq:co_main_ws} we obtain 
    $
   d_2(\hat{\M U}^{(k_0)}_{g}, \M U^{(k_0)}) = (1 + o(1))E_n/\sigma_{k_0}  
  $ for all $g \geq 1$,  with probability at least $1 - 2n^{-3}$. 
    In other words, if we estimate $\M U^{(k_0)}$ using only the
    leading singular vectors of the sketched matrix $\hat{\M M} \M G$ then we need a
    sketching dimension of $\sk = \Omega(n)$.
    A similar phenomenon was observed for PCA using random projections \citep{yang2021reduce}. 
  \item Next suppose $k_0$ is bounded by a finite constant and $\sk \asymp \log n$.
    If $\sigma_{k_0}/E_n \succsim n^{\epsilon}$ for some fixed but arbitrary
    $\epsilon > 0$ then  \eqref{thm2*:resopt:main} implies
 \begin{equation}
 \label{eq:d2_logn}
     d_2(\hat{\M U}_g^{(k_0)}, \M U^{(k_0)}) \leq (1 + o(1)) E_n/\sigma_{k_0}  \text{ for all $g \geq 1 + (2 \epsilon)^{-1}$}, 
 \end{equation}
 with probability at least $1 - 2n^{-3}$, as $(2\epsilon)^{-1} -g_* = \omega(\log^{-1} n)$.
 \end{enumerate}
\end{remark}
\begin{co}\label{co:l2inf_noise}
  Assume the setting of Corollary~\ref{co:l2_noise} with $\sigma_{k_0} > 2 E_n$ where $E_n = \|\M E\|$.
Define
   \begin{gather*}
       g_* = \max\Bigl\{1 + \frac{\log(nk_0/\sk) + \log \log n}{2 \log(\tfrac{1}{2}\sigma_{k_0}/E_n)}, 
       \frac{1}{3} + \frac{\log(n^{3/2}/\sk) + \tfrac{1}{2} \log \log n}{3 \log(\tfrac{1}{2}\sigma_{k_0}/E_n)}, \frac{1}{2} + \frac{\log(n^{3/2} \|\M U^{(k_0)}\|_{\twoinf}/\sk)}{2 \log (\tfrac{1}{2} \sigma_{k_0}/E_n)}\Bigr\}, 
    \end{gather*}
Suppose also that $d_{\twoinf}(\hat{\M U}^{(k_0)}, \M U^{(k_0)}) \leq \|\M U^{(k_0)}\|_{\twoinf}$. 
We then have, for all $g \geq g_*$, 
    \begin{gather}
    \label{eq:twoinf_co_main2}
        d_{\twoinf}(\hat{\M U}_g^{(k_0)}, \M U^{(k_0)}) \leq 
        d_{\twoinf}(\hat{\M U}^{(k_0)}, \M U^{(k_0)}) + O\bigl(n^{-1/2}E_n/\sigma_{k_0}\bigr)
    \end{gather}
    with probability at least $1 - 5n^{-3} - 2e^{-n/2}$, where the hidden factor in $O(n^{-1/2} E_n/\sigma_{k_0})$
    only depend on $c_{\mathrm{gap}}$. Furthermore for all $g \geq g_*$ we have
    \begin{equation}
      \label{eq:entrywise_guarantee}
      \|\hat{\M U}^{(k_0)}_{g} \hat{\M U}^{(k_0)\top}_{g} \hat{\M M} - \M M\|_{\max}
      \leq \|\hat{\M U}^{(k_0)} \hat{\M U}^{(k_0)\top} \hat{\M M} - \M M \|_{\max} + O
      \bigl(\kappa n^{-1/2} E_n \|\M U\|_{\twoinf}\bigr),
    \end{equation}
    with probability at least $1 - 9n^{-3} - 2e^{-n/2}$,
    where $\kappa = \sigma_1/\sigma_{k_0}$ is the condition number for $\M M$ and
    the hidden factor in $O(n^{-1/2} \kappa E_n \|\M U^{(k_0)}\|_{\twoinf})$ only depends on $c_{\mathrm{gap}}$.  

    If $\sigma_{k} = \omega(E_n)$ then the terms $O(n^{-1/2} E_n /\sigma_{k_0})$ and $O(\kappa n^{-1/2} E_n \|\M U^{(k_0)}\|_{\twoinf})$ in Eq.~\eqref{eq:twoinf_co_main2} and Eq.~\eqref{eq:entrywise_guarantee}
    can be replaced by $o(n^{-1/2} E_n/\sigma_{k_0})$ and $o(\kappa n^{-1/2} E_n \|\M U^{(k_0)}\|_{\twoinf})$ provided that
    $g \geq g_* + \Delta$ where $\Delta > 0$ is any arbitrary constant. 
    Finally if $\sigma_{k_0}/E_n \succsim n^{\epsilon}$ for a fixed but arbitrary $\epsilon > 0$ then 
    the above expression for $g_*$   
    can be simplified to
   \begin{gather}
     \label{eq:g2_strong}
       g_* = \max\Bigl\{1 + \frac{\log(nk_0/\sk)}{2 \log(\sigma_{k_0}/E_n)}, 
       \frac{1}{3} + \frac{\log(n^{3/2}/\sk)}{3 \log(\sigma_{k_0}/E_n)}, \frac{1}{2} + \frac{\log(n^{3/2} \|\M U^{(k_0)}\|_{\twoinf}/\sk)}{2 \log(\sigma_{k_0}/E_n)} \Bigr\}.
    \end{gather}
 \end{co}

\subsection{Comparison with existing results}
\label{sec:related_works}
For ease of exposition we will fix some choice of $k$ and thus drop the index $k$ from the matrices $\hat{\M U}_g^{(k)}$ and $\hat{\M U}^{(k)}$. Also, we will implicitly assume that $\tilde{g} = g$ if $\hat{\M M}$ is
symmetric and $\tilde{g} = 2g+1$ otherwise. 
As we alluded to in the introduction, 
theoretical analysis of RSVD mostly focused on
spectral and Frobenius norms
upper bounds for $\sin \Theta(\hat{\M U}_g, \hat{\M U})$ and
$(\M I - \hat{\mathbf{U}}_g \hat{\M U}_g^{\top})
\hat{\M M}$ in
settings where $\hat{\M M}$ is assumed to be {\em noise-free}. 

Let $\vertiii{\cdot}$ be any unitarily invariant norm. 
By combining
Theorems~4 and 6 in \cite{saibaba}, we have
\begin{equation}
  \label{eq:saibaba}
  \vertiii{\sin \Theta(\hat{\M U}_g, \hat{\M U})} \leq 
  \frac{C \sk^{1/2} \vertiii{\hat{\bm{\Sigma}}_{\perp}^{\tilde{g}}}}
  {\hat{\sigma}_{k}^{\tilde{g}} (1 - \hat{\zeta}_{k})(\sk - k + 1) \eta^{1/(\sk - k + 1)}}
  \bigl((n - k)^{1/2} + \sk^{1/2} + (\log 1/\eta)^{1/2} \bigr) 
\end{equation}
with probability at least $1 - \eta$, where $C$ is a universal constant. Ignoring constant factors,
Eq.~\eqref{eq:ui_corollary_bound} is similar to Eq.~\eqref{eq:saibaba}, 
with the main difference being that Eq.~\eqref{eq:saibaba} include extra factors
$\eta^{1/(\tilde{k} - k +1)}$ and  
$(1 - \hat{\zeta}_k)$ in the denominator. 
Our bound in Eq.~\eqref{eq:ui_corollary_bound} is therefore sharper than Eq.~\eqref{eq:saibaba}
when $\hat{\zeta}_k = \hat{\sigma}_{k+1}/\hat{\sigma}_k \approx 1$ and/or $\eta \rightarrow 0$. 

Next, for $(\M I - \hat{\M U}_g \hat{\M U}_g^{\top}) \hat{\M M}$,
  let $\bm{\Pi}_{\M Y_g}$ be the orthogonal projection onto the column space of $\M Y_g$ where $\M Y_g =
\hat{\M M}^{g} \M G$ if $\hat{\M M}$ is symmetric and $\M Y_g = \hat{\M M}(\hat{\M M}^{\top} \hat{\M M})^{g} \M G$ otherwise. 
Then by Corollary~10.10 in \cite{halko2011finding} we
have for $\sk \geq k + 2$ that
\begin{equation}
  \label{eq:prior1}
  \begin{split}
\|(\M I - \bm{\Pi}_{{\M Y}_g}) \hat{\M M}\| &
\leq \Bigl(1 + \frac{k^{1/2}}{(\sk - k - 1)^{1/2}} + \frac{e \sk^{1/2} (\min\{m,n\} - k)^{1/2}}{\sk - k} \Bigr)^{1/\tilde{g}} \hat{\sigma}_{k+1}
    \end{split}
\end{equation}
with high probability. 
Meanwhile, from Theorem~\ref{thm2ss} we have by the triangle inequality that
\begin{equation}
  \label{eq:low_rank1_approximation}
  \begin{split}
  \|(\M I - \hat{\M U}_g \hat{\M U}_g^{\top}) \hat{\M M}\|\leq C (n/\sk)^{1/2} (\hat{\sigma}_{k+1}/\hat{\sigma}_{k})^{\tilde{g}} \hat{\sigma}_1 +
  \hat{\sigma}_{k+1}  
  \end{split}
\end{equation}
with high probability, where $C$ is a universal constant.
Comparing Eq.~\eqref{eq:prior1} and Eq.~\eqref{eq:low_rank1_approximation} we see that they 
both converge to $\hat{\sigma}_{k+1}$ as $\tilde{g}$ increases, but the manner in which they converge can be
significantly different.
More specifically Eq.~\eqref{eq:prior1} does not depend on the singular value gap $\hat{\sigma}_{k+1}/\hat{\sigma}_{k}$ but its rate of convergence becomes slower for larger values of $\tilde{g}$.
In contrast, Eq.~\eqref{eq:low_rank1_approximation} depends on $(\hat{\sigma}_{k+1}/\hat{\sigma}_k)^{\tilde{g}}$ and
converges to $\hat{\sigma}_{k+1}$ rather quickly if $\hat{\sigma}_{k+1} \ll \hat{\sigma}_{k}$ but
arbitrarily slowly when $\hat{\sigma}_{k+1} \approx \hat{\sigma}_{k}$. 
For many statistical applications including
those in Section~\ref{sec:rgi} and Sections~\ref{sec:mc}-\ref{sec:epca}, 
$\M M$ will be approximately low-rank so that 
bounds based on
Eq.~\eqref{eq:low_rank1_approximation} are sharper than those derived from Eq.~\eqref{eq:prior1}. 

We now discuss the upper bounds
for $d_{2 \to \infty}(\hat{\M U}_g, \hat{\M U})$. Perturbation bounds in $\ell_{\twoinf}$ norm
is rarely studied in the RSVD literature and among existing results, the one that is perhaps
most related to ours is from Proposition~1 in \cite{twoinf_subspace_damle} wherein
the authors approximate the leading eigenvectors of a symmetric $\hat{\M M}$ using a subspace iteration procedure similar to that in Algorithm~\ref{RSRS} with $\sk = k$, and they showed that
\begin{equation}
  \label{eq:damle_twoinf}
  d_{\twoinf}(\hat{\M U}_g, \hat{\M U}) \leq \frac{1}{\sqrt{1 - d_0^2}}
    \Bigl(\frac{\hat{\sigma}_{k+1}}{\hat{\sigma}_{k}}\Bigr)^{g}
    \Bigl[\sqrt{2} d_0 +  C_* (1 + \sqrt{n} \|\hat{\M U}\|_{\twoinf})
    d_{\twoinf}(\hat{\M U}_0, \hat{\M U})\Bigr]
\end{equation}
for all $g \geq 1$, where $\hat{\M U}_0$ denote any arbitrary initial estimate for $\hat{\M U}$ and 
$d_0 := d_{2}(\hat{\M U}_0, \hat{\M U})$. 
Eq.~\eqref{eq:damle_twoinf} is similar in spirit to
ours Eq.~\eqref{eq:d2inf_corollary_bound}, with the most obvious difference being that
Eq.~\eqref{eq:d2inf_corollary_bound} include terms with extra factors of $(n/\sk)^{1/2}$ and 
$(\hat{\sigma}_{k+1}/\hat{\sigma}_k)^{\tilde{g}}$. 
The main limitation of Eq.~\eqref{eq:damle_twoinf}, however, lies in the fact that it depends
on the assumption
\begin{equation}
  \label{eq:condition}
  \|\hat{\M U}_{\perp} \hat{\bm{\Sigma}}^{g} \hat{\M U}_{\perp}^{\top}\|_{\infty} \leq C_* \hat{\sigma}_{k+1}^{g} \|\M I - \hat{\M U} \hat{\M U}^{\top}\|_{\infty}
\end{equation}
for all $g \geq 1$, where $\|\cdot \|_{\infty}$ denote the maximum row-sum (after taking absolute values) of a matrix, and $C_*$ is the constant appearing in Eq.~\eqref{eq:damle_twoinf}; see Assumption~1 in \cite{twoinf_subspace_damle} for more details. The condition in Eq.~\eqref{eq:condition} can be   restrictive. 
In particular, for many statistical applications including random graph inference and matrix completion,
we have $\hat{\M M} = \M M + \M E$ where $\M M$ is low-rank,
$\M U$ has bounded coherence, and $\|\M E\|_{\infty} \gg \|\M E\|$. Letting $k = \mathrm{rk}(\M M)$, we typically have
$\|\hat{\M U} \hat{\bm{\Sigma}} \hat{\M U}^{\top} - \M M\|_{\infty} \ll \|\M E\|_{\infty}$ and
$\|\hat{\M U} \hat{\M U}^{\top}\|_{\infty} \asymp 1$ so that 
\begin{equation*}
  \|\hat{\M U}_{\perp} \hat{\bm{\Sigma}}_{\perp} \hat{\M U}_{\perp}^{\top}\|_{\infty} \geq \|\hat{\M M} - \M M\|_{\infty} - \|\hat{\M U} \hat{\bm{\Sigma}} \hat{\M U}^{\top}-\M M\| \succsim \|\M E\|_{\infty} \gg \|\M E\| \asymp \hat{\sigma}_{k+1} \|\M I - \hat{\M U} \hat{\M U}^{\top}\|_{\infty}, 
\end{equation*}
and hence Eq.~\eqref{eq:condition} does not hold for $g = 1$.
In contrast, the only assumption we need for Theorem~\ref{thm3ss:general} is that $\sk$ is lower bounded, i.e., $\sk \geq (1 - c_{\mathrm{gap}})^{-2}\{\sqrt{k} + \sqrt{2 \log \vartheta^{-1}}\}^2$
where $c_{\mathrm{gap}}$ and $\vartheta$ are both arbitrary. 
The discrepancies between the above sets of assumptions 
is mainly because the analysis in \cite{twoinf_subspace_damle} is
for subspace iteration and not
{\em randomized} subspace iteration.  More specifically, \cite{twoinf_subspace_damle} viewed
$\hat{\M U}_0$ as given but arbitrary and thus they cannot control
$\hat{\M U}_{\perp} \hat{\bm{\Sigma}}_{\perp}^{\tilde{g}} \hat{\M U}_{\perp}^{\top} \hat{\M U}_0$
for all $g$ and all $\hat{\M U}_0$ without making (possibly restrictive) assumptions on $\hat{\M M}$ itself.
In contrast, as we assume that $\hat{\M U}_0$ is a random matrix independent of $\hat{\M M}$ and only dependent on $\M G$, 
we can leverage the randomness in $\hat{\M U}_0$ to bound
$\hat{\M U}_{\perp} \hat{\bm{\Sigma}}_{\perp}^{\tilde{g}} \hat{\M U}_{\perp}^{\top} \hat{\M U}_0$
{\em conditional} on $\hat{\M M}$ and thus alleviate
the need to make assumptions on $\hat{\M M}$ (not to mention that, by carefully exploiting the Gaussianity of $\M G$, we also obtain
more precise $\ell_{2,\infty}$ error bounds compared to that for some given but arbitrary $\hat{\mu}_0$). 
{\color{black}See Remark~\ref{rk:twoinf} for more discussions.}
Finally we note that Eq.~\eqref{eq:max_corollary_bound}, Eq.~\eqref{eq:entrywise_co_general} and Eq.~\eqref{eq:entrywise_guarantee} are, to the best of our
knowledge, the first set of bounds for entrywise differences
between the RSVD-based low-rank
approximation $\hat{\M U}_g \hat{\M U}_g^{\top} \hat{\M M}$ and the truncated exact SVD $\hat{\M U} \hat{\M U}^{\top} \hat{\M M}$. If $\hat{\M M} = \M M + \M E$ then these
results also yield entrywise bounds
for $\M M - \hat{\M U}_g \hat{\M U}_g^{\top} \hat{\M M}$ as well 
as normal approximations and entrywise
confidence intervals for $\M M$; see Section~\ref{sec:mc} for more details.

\section{Random graph inference}\label{sec:rgi}
We now apply Corollary \ref{co:l2_noise} and Corollary~\ref{co:l2inf_noise}
to estimate the leading eigenvectors for
edge-independent random graphs with low-rank edge probabilities
matrices. Additional applications to matrix completion and PCA with missing data are presented in Section~\ref{sec:mc} and Section~\ref{sec:epca:supp} of the Supplementary File.
Let $\hat{\mathbf{M}} = [\hat{m}_{ij}]$ be the 
adjacency matrix of a random graph on $n$ vertices with edge probabilities
$\M M = [m_{ij}]$, i.e.,  
$\hat{\M M}$ is a symmetric, binary matrix  whose upper triangular entries are independent
Bernoulli random variables with $\mathbb{P}(\hat{m}_{ij} = 1) =
m_{ij}$. Suppose $\mathrm{rk}(\M M) = k_0$ for some constant $k_0$, 
$\|\M U^{(k_0)}\|_{2 \to \infty} \asymp k_0^{1/2} n^{-1/2}$, 
$\sigma_{k_0} \asymp n \rho_n$, and $\|\M E\| \precsim (n \rho_n)^{1/2}$ with probability at least $1 - n^{-3}$; here
$\rho_n \in [0,1]$ satisfies $n \rho_n = \Omega(\log n)$ as $n$ increases. 

The above assumption are quite mild and are
 satisfied by many random graph models
 including Erd\H{o}s--R\'enyi, stochastic blockmodels
 and its degree-corrected and/or mixed-membership variants
 \citep{holland1983stochastic,karrer2011stochastic,Airoldi2008},
 (generalized) random dot product graphs
 \citep{grdpg1}, as well as any edge-independent random graph whose edge
 probabilities are sufficiently homogeneous, i.e.,  $\max_{i} \sum_{j} m_{ij}
 \asymp \min_{i} \sum_{j} m_{ij}$.
The factor $\rho_n$ corresponds to the sparsity of
 $\hat{\M M}$, i.e., with high probability $\hat{\M
   M}$ has $\Theta(n^2 \rho_n)$ non-zero entries.
 We note that random graph inference using spectral methods generally require
 $n \rho_n \succsim \log{n}$ 
 as otherwise the leading eigenvalues and eigenvectors of $\hat{\M M}$
 may not yield consistent estimates for the corresponding eigenvalues and eigenvectors of $\M M$. 
 Finally 
 the assumption $\mathbb{P}(\|\M E\| \precsim (n \rho_n)^{1/2}) \geq 1 - n^{-3}$
 is chosen mainly for convenience as,
 under the above conditions, by using standard matrix concentration inequalities \citep{oliveira,tropp,lei2015consistency} we can show that for any finite constant $c > 0$ there exists a finite constant $C > 0$ depending only on $c$ such that
 $\mathbb{P}(\|\M E\| \leq C (n \rho_n)^{1/2}) \geq 1 - n^{-c}$.

 \subsection{Subspace perturbation error bounds}\label{sec:int:network}
Suppose $n$ is large and we are interested in computing the $k_0$ leading
singular vectors of $\hat{\M M}$ as an estimate for the $k_0$ leading
singular vectors of $\M M$. To save computational time, we will use
Algorithm \ref{RSRS} with some choices of $g \geq 1$ and $\sk \geq k_0$. 
Then by Corollary~\ref{co:l2_noise} we have
\begin{equation}
\label{thm2:rate:pre1}
  d_{2}(\hat{\M U}_g, \M U) =
    O((n \rho_n)^{-1/2}), \quad \text{ for all $g \geq 1 + \frac{(1 + o(1)) \log(n/\sk)}{\log(n \rho_n)}$}
\end{equation}
with high probability, 
where the hidden factor in $o(1)$ does not depend on $n, \sk$ and $\rho_n$. Note that, for ease of exposition, we omit the dependency on \( k_0 \) for all matrices.  

Furthermore, under the {\it strong signal regime} $n \rho_n = \Omega(n^{\beta})$ for some fixed but arbitrary
$\beta > 0$, by Corollary~\ref{co:l2_noise}, Eq.~\eqref{thm2:rate:pre1} can be strengthened to
\bee\label{thm2:rate:pre2}{
d_2(\hat{\M U}_g,\M U) = 
\begin{cases}
 O\left( (n\rho_n)^{-1/2}\right)   & \text{if $g > 1 + \alpha_*$}
\\
O\bigl((n \rho_n)^{-(g - \alpha_*)/2}\bigr) & \text{if  $\alpha_* \leq g \leq 1 + \alpha_*$},
\end{cases}}
\ee
with probability at least $1 - 3n^{-3}$, where $\alpha_* = \frac{\log(n/\sk)}{\log n \rho_n} \leq \beta^{-1} + {\color{black}o(1)}$. 
Eq.~\eqref{thm2:rate:pre1} and Eq.~\eqref{thm2:rate:pre2} imply a
phase transition for $d_2(\hat{\M U}_g,\M U)$. In particular,
the assumptions at the beginning of this section together with the Davis-Kahan
theorem imply $d_2(\hat{\M U}, \M U) \leq \frac{E_n}{\sigma_{k_0}}$ with
high probability. 
Therefore, under the strong signal regime, if $g\geq 1 + \alpha_*$ then $d_2(\hat{\M U}_g, \M
U)$ is asymptotically equivalent to $d_2(\hat{\M U}, \M U)$. In contrast,  
if $\alpha_* < g < 1 + \alpha_*$ then
$d_2(\hat{\M U}_g,\M U)$ converges to $0$ at the slower rates of
$n^{-(g - \alpha_*)/2}$. Convergence
of $d_2(\hat{\M U}_g, \M U)$ to $0$ is not guaranteed when $g \leq
\alpha_*$. Finallly, under the {\it weak signal regime} $n \rho_n \succsim \log n$,
we can similarly show that $d_{2}(\hat{\M U}_g, \M U)$ is asymptotically equivalent to $d_{2}(\hat{\M U}, \M U)$ when $g \geq
1 + \frac{(1 + o(1)) \log(n/\sk)}{\log n \rho_n}$, while
convergence of $d_{2}(\hat{\M U}_g, \M U)$ to $0$ is not guaranteed when
$g \leq \frac{(1 - o(1)) \log(n/\sk)}{\log n \rho_n}$.  

We now discuss the behavior of $d_{\twoinf}(\hat{\M U}_g, \M U)$. For random graphs satisfying the assumptions in this section, the sharpest known upper bound for $d_{\twoinf}(\hat{\M U}, \M U)$ is
\begin{equation}\nonumber
  d_{\twoinf}(\hat{\M U}, \M U) = O\Big((\log n)^{1/2} \frac{E_n}{\sigma_{k_0}} \|\M U\|_{\twoinf} \Big)= O\Big( \frac{(\log n)^{1/2}}{n \rho_n^{1/2}}\Big)
\end{equation}
with high probability; see e.g., \cite{cape2019signal,grdpg1,abbe2020entrywise,mao2020estimating,xie2021entrywise}
for more details. Define
\[g_* = \max\Bigl\{1 + \frac{ (1 + o(1)) (\log(n/\sk) + \log \log n)}{\log n \rho_n}, \frac{1}{3} + \frac{(1 + o(1))(\log(n^3/\sk^2) + \log \log n)}{3 \log n \rho_n}\Bigr\}.\] 
where the hidden factor in $o(1)$ does not depend on $n, \sk$ and $\rho_n$. 
Then by Corollary~\ref{co:l2inf_noise}, for $g \geq g_*$ we have
\begin{equation}
  \label{eq:d_twoinf_rg}
  d_{\twoinf}(\hat{\M U}_g, \M U) = 
    O\Bigl(\frac{(\log n)^{1/2}}{n \rho_n^{1/2}}\Bigr)  
\end{equation}
with probability at least $1 - 6n^{-3}$.  
Furthermore, under the strong signal regime $n \rho_n = \Omega(n^{\beta})$,
  {a careful investigation of the proof of Corollary~\ref{co:l2inf_noise}} can strengthen Eq.~\eqref{eq:d_twoinf_rg} to yield
\begin{equation}
  \label{eq:d_twoinf_rg_part2}
  d_{\twoinf}(\hat{\M U}_g, \M U) = \begin{cases} O\bigl(\frac{(\log n)^{1/2}}{n \rho_n^{1/2}}\bigr) & \text{if $g \geq \max\{1 + \alpha_*^{(1)}, 1/3 + \alpha_*^{(2)}\}$} \\
    O\bigl(\frac{(\log n)^{1/2}}{n^{1/2} (n \rho_n)^{\xi_g}}
    \bigr) &
\text{if $\max\{\alpha_*^{(1)}, \alpha_*^{(2)}\} < g < \max\{1 + \alpha_*^{(1)}, 1/3 + \alpha_*^{(2)}\}$} 
  \end{cases}
\end{equation}
where $\alpha_*^{(1)} = \frac{\log(n/\sk)}{\log n \rho_n} \leq \beta^{-1}$, $\alpha_*^{(2)} =
  \frac{\log(n^{3}/\sk^2)}{3 \log n \rho_n} \leq \beta^{-1}$, and
  $\xi_g = \tfrac{1}{2} \min\{g - \alpha_*^{(1)}, 3(g - \alpha_*^{(2)})\} \leq \tfrac{1}{2}$.
Eq.~(\ref{eq:d_twoinf_rg}) and Eq.~(\ref{eq:d_twoinf_rg_part2}) also imply a
phase transition for $d_{\twoinf}(\hat{\M U}_g, \M U)$ as $g$ increases, e.g., 
 if $g \geq \max\{1 + \alpha_*^{(1)}, 1/3 + \alpha_*^{(2)}\}$ then
$d_{\twoinf}(\hat{\M U}_g, \M U)$ converges to $0$ at the
same rate as $d_{\twoinf}(\hat{\M U}, \M U)$, while if $\max\{\alpha_*^{(1)}, \alpha_*^{(2)}\} < g 
< \max\{1 + \alpha_*^{(1)}, 1/3 + \alpha_*^{(2)}\}$ then $d_{\twoinf}(\hat{\M U}_g, \M U)$ converges to $0$ at a slower (degenerate) rate of $n^{-1/2} (n \rho_n)^{-\xi_g}$; convergences
of $d_{\twoinf}(\hat{\M U}_g, \M U)$ is not guaranteed when $g \leq \max\{\alpha_*^{(1)}, \alpha_*^{(2)}\}$.
\par
\begin{figure}[t]
\includegraphics[width=0.8\linewidth]{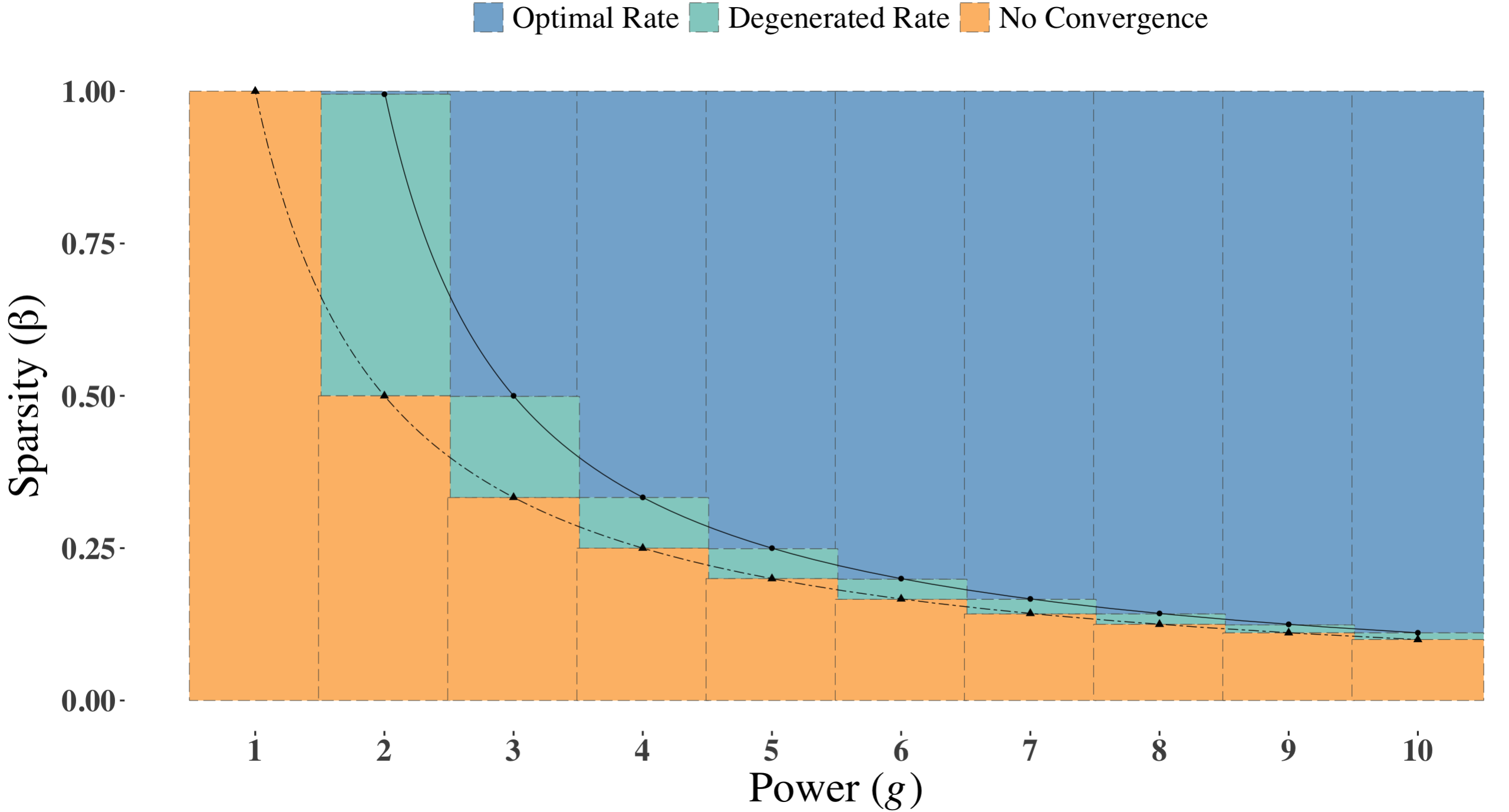} \centering
\caption{Phase transition diagram of error rates under the random graph setting with \(\sk \asymp \log n\) and $n \rho_n \asymp n^{\beta}$. Different regions correspond to different convergence rates of $d_2(\hat{\mathbf{U}}_g, \mathbf{U})$ or $d_{2 \to \infty}(\hat{\mathbf{U}}_g, \mathbf{U})$; see Eq.~\eqref{rate:simple} for details. The X- and Y-axes represent $g$ and $\beta$, respectively, and dashed and solid lines represent the thresholds $g = \beta^{-1}$ and $g = 1 + \beta^{-1}$.}
  \label{fig:phase}
\end{figure}

 
{\color{black}To further illustrate our   results, we focus on a practical and computationally efficient scenario where the RSVD sketching dimension is $\sk \asymp \log n$
 under the strong signal regime $n\rho_n\asymp n^\beta$. Eq.~\eqref{thm2:rate:pre2} and Eq.~\eqref{eq:d_twoinf_rg_part2} then become
\begin{gather}\label{rate:simple}
 d_2(\hat{\M U}_g,\M U) = \begin{cases}
O( n^{-\beta/2} ) & \text{if $g \geq 1 + \beta^{-1}$}
\\
  O(n^{(-\beta g + 1)/2})  & \text{if  $\beta^{-1} \leq g \leq 1 + \beta^{-1}$}
\end{cases},
\\
\label{rate:simple2}
 d_{\twoinf}(\hat{\M U}_g, \M U) = (\log n)^{1/2}\begin{cases}
 n^{-\beta/2-1/2}   & \text{if $g \geq 1 + \beta^{-1}$}
\\
  n^{(-\beta g + 1)/2-1/2}  & \text{if  $\beta^{-1} \leq g \leq 1 + \beta^{-1}$}
\end{cases}. 
\end{gather}
Similar to our previous discussions, the above rates are optimal when $g \geq 1 +\beta^{-1}$, degenerate when $\beta^{-1} \leq g \leq 1 + \beta^{-1} $, and no convergence
can be guaranteed when $g \leq \beta^{-1}$. 
See Section~\ref{sec:lower_bound} for matching lower bounds (up to logarithmic factors), Section~\ref{sec:ptv} for
further numerical evidence, and Figure \ref{fig:phase} for a summary of these phase transitions. 
In summary the convergence rates in Eq.~\eqref{rate:simple} and Eq.~\eqref{rate:simple2}
support the well-known recommendation that selecting \(\sk\) slightly larger than $k_0$ and increasing \(g\)
is crucial for the practical success of RSVD \citep{martinsson_rnla}.  }

\subsection{Lower bound and phase transition sharpness}
\label{sec:lower_bound} 
We now study the sharpness of the phase transition thresholds
described in Section~\ref{sec:int:network}. 
More specifically, we derive lower bounds for
$d_2(\hat{\M U}_g, \M U)$ and $d_{\twoinf}(\hat{\M U}_g, \M U)$ under the same regime as the upper bounds in Eq.~\eqref{rate:simple}, namely with
$n \rho_n \asymp n^{\beta}$ and $\sk \asymp \log n$. 
Our lower bounds depend on the following
assumption for the growth rate of
 $\mathrm{tr} \, \hat{\M M}^{2g}$. \begin{assumption}
  \label{ass:lower_bound}
For any $g \geq 1$ there exists a constant $c_g > 0$ depending on $g$ such that
\bee\label{network:setting3}
\mathbb{E}( \mathrm{tr} \,\hat{\M M}^{2g}) \geq c_g \left\{n^{g+1} \rho_n^{g} + (n \rho_n)^{2g}\right\}.
\ee
\end{assumption}
Assumption~\ref{ass:lower_bound}
is satisfied by any edge-independent random graphs with homogeneous variances including 
stochastic blockmodel graphs, their degrees-corrected and mixed
membership variants  \citep{karrer2011stochastic,Airoldi2008},
and (generalized) random dot product graphs \citep{grdpg1}.
More specifically,
let $\M M$ be such that $m_{ij} \asymp \rho_n$ for all $i,j$.
Then by Theorem~1 and Lemma~S.7 in \cite{maugis}, we have 
$$\mathbb{E}(\mathrm{tr} \, \hat{\M M}^{2g}) \asymp_g (n \rho_n)^{2g} + n^{g+1} \mathbb{P}\bigl(\{v_1, v_2, \dots, v_{g+1}\} \text{ forms a tree}\bigr) \asymp_g n^{g+1} \rho_n^{g} + (n \rho_n)^{2g}$$
where $\{v_1, \dots, v_{g+1}\}$ is any collection of $g+1$ distinct indices in $\hat{\M M}$.
\begin{theorem}\label{thm:lower}
  Assume the setting of Corollary~\ref{co:l2_noise} where $k_0$ is bounded by a finite constant not depending
  on $n$. 
  Furthermore suppose that (i) $\sigma_{1} \asymp n \rho_n \asymp n^{\beta}$ and $\sigma_{k_0}/E_n \succsim
  n^{\beta/2}$ 
  for some constant $\beta> 0$ and (ii) $\hat{\M M}$ 
  satisfies Assumption~\ref{ass:lower_bound}. 
Choose $\sk = 2 c_{\mathrm{gap}} \log n$ and let $p_0 \in (0,1)$ be fixed but
arbitrary. Then there exists a constant $C_{p_0}>0$
depending only on $p_0$ such that for all $g < \beta^{-1}$ we have with probability at least $1 - 2p_0$,
\bee\nonumber
d_2(\hat{\M U}_g,\M U) 
\geq C_{p_0} c_g, \quad   
d_{\twoinf}(\hat{\M U}_g,\M U) 
\geq C_{p_0} c_g \cdot n^{-1/2},
\ee
and for all $\beta^{-1}\leq g < 1 +\beta$, we have with probability at least $1 - 2p_0$,
$$
d_2(\hat{\M U}_g,\M U) 
\geq C_{LB}c_g\cdot (\log n)^{-1/2}n^{  (- g\beta + 1)/2}, \quad  
d_{\twoinf}(\hat{\M U}_g,\M U) \geq C_{LB}c_g\cdot(\log n)^{-1/2}n^{  (- g\beta + 1)/2 - 1/2} .
$$
\end{theorem}
Theorem~\ref{thm:lower} implies that if \(g < \beta^{-1}\) then \(\hat{\mathbf{U}}_g\) is not a consistent estimator for \(\mathbf{U}\). In contrast, if 
\(g \geq 1 + \beta^{-1}\) then \(d_2(\hat{\mathbf{U}}_g, \mathbf{U})\) and \(d_{2 \to \infty}(\hat{\mathbf{U}}_g, \mathbf{U})\) attain the same convergence rates as 
\(d_2(\hat{\mathbf{U}}, \mathbf{U})\) and \(d_{2 \to \infty}(\hat{\mathbf{U}}, \mathbf{U})\) (see Eq.~\eqref{rate:simple} and Eq.~\eqref{rate:simple2}), 
so that \(\hat{\mathbf{U}}_g\) is rate-optimal. 
\subsection{Exact recovery for stochastic blockmodels}\label{sec:sbm}
We first recall the notion of stochastic blockmodel graphs \citep{holland1983stochastic}, one of the
most widely used generative model for networks with an  intrinsic
community structure. 
\begin{definition}[SBM]\label{def:sbm}
  Let $K \geq 1$ be a positive integer and let $\bm{\pi} \in
  \mathbb{R}^{K}$ be a probability vector. Let $\mathbf{B}$ be a
  symmetric $K \times K$ matrix whose entries are in $[0,1]$. We say that $(\mathbf{A}, \bm{\tau}) \sim
  \mathrm{SBM}(\mathbf{B}, \bm{\pi})$ is a $K$-blocks stochastic
  blockmodel (SBM) graph  with parameters $\mathbf{B}$ and $\bm{\pi}$, and
  sparsity factor $\rho_n$, if the following holds.
  First, sample $\bm{\tau} = (\tau_1, \dots, \tau_n)$ where the $\tau_i$ are
  iid with $\mathbb{P}(\tau_i = \ell) = \pi_{\ell}$ for all $\ell \in [K]$. Then
  sample $\mathbf{A}$ as a $n \times n$ symmetric binary matrix 
  whose upper triangular entries $\{A_{ij}\}_{i \leq j}$ 
  are independent Bernoulli random variables with
  $\mathbb{P}(A_{ij} = 1) = \rho_n B_{\tau_i, \tau_j}$.
 \end{definition}
Community detection is a well-known problem with many available
techniques; see \citet{abbe2017community} for a survey.
  We focus on spectral clustering,
a simple and popular 
algorithm wherein, given $\mathbf{A}$, 
we first choose an embedding dimension $d$ and compute the matrix $\hat{\M U}$ of eigenvectors corresponding to the
$d$ largest (in modulus) eigenvalues of $\M A$. Next we cluster the rows
of $\hat{\M U}$ into $K$ cluster using either the $K$-means or
$K$-medians algorithms and let $\hat{\tau}_i$ be the resulting cluster
membership for the $i$th row of $\hat{\M U}$. 
Statistical properties of spectral clustering had been widely studied
in recent years. See e.g., \cite{rohe2011spectral, lei2015consistency,abbe2020entrywise} for an incomplete list of references. 

In   real-world applications like social or biological networks,
the number of nodes can be on the order of $10^6$ \citep{gopalan2013efficient}. As a result,
obtaining $\hat{\M U}$ using standard
SVD algorithms can be prohibitively demanding in terms of both
the computational time and memory requirement.
We thus consider a RSVD-based 
spectral clustering procedure  
that replaces $\hat{\M U}$
by its approximation $\hat{\M U}_g$, which  
results in a procedure with running time
$O(g \sk \times \mathrm{nnz}(\mathbf{A}))$  
and memory consumption $O(\mathrm{nnz}(\mathbf{A}))$ where $\mathrm{nnz}(\M A)$ denote the number of non-zero
entries of $\M A$. 
We leverage our bounds for
$d_{\twoinf}(\hat{\M U}_g, \M U)$  to show that
$K$-means clustering over $\hat{\M U}_g$ yields an exact recovery of $\bm{\tau}$  with high probability, i.e., it
yields a $\hat{\bm{\tau}}$ such that there exists a
permutation $\varsigma$ of $\{1,\dots,K\}$ for which $\hat{\tau}_i =
\varsigma(\tau_i)$ for all $i \in \{1,\dots,n\}$. 
\begin{theorem}\label{thm:sbm}
Let $(\mathbf{A}, \bm{\tau}) \sim \mathrm{SBM}(\M
B,\bds\pi)$ be a $K$-blocks SBM with sparse parameter $\rho_n$ where $\pi_{\ell} >0$ for all $\mathrm{\ell}\in[K]$ 
and $n\rho_n = \omega(\log n)$. Let $d =
\mathrm{rk}(\mathbf{B})$ and suppose that
$\hat{\bds\tau}$ is the $K$-means clustering of the rows of $\hat{\M U}_g$ where
$\sk \asymp \log n$ and $g \geq \frac{2 \log n}{\log n \rho_n}$.
Then for sufficiently large $n$, $\hat{\bds\tau}$ exactly recovers $\bds\tau$
with probability at least $1 - 2n^{-3}$. If $n \rho_n \succsim n^{\beta}$ for any fixed but arbitrary
$\beta > 0$ then the
threshold for exact recovery of $\hat{\bds{\tau}}$ can be sharpened to
$g \geq 1 + \beta^{-1}$. 
\end{theorem}
\begin{remark}
  \label{rem:earlier_RSVD}
  Community detection using RSVD was also studied
  in \cite{randomized_sc}. In particular, Theorem~1 in
  \cite{randomized_sc} shows that if $n \rho_n = \omega(\log n)$ then weak
  recovery is possible using $\hat{\M U}_g$,
  provided that one choose $\tilde{k} \geq k + 4$ and $g =
  \Omega(n^{\delta})$ for any fixed but arbitrary $\delta >
  0$. Recall that weak
  recovery only requires the {\em proportion} of mis-clustered vertices
  to converge to $0$. Comparing the results in \cite{randomized_sc} to Theorem~\ref{thm:sbm}, we see that while they have
  the same sparsity requirement, the exact recovery in 
  Theorem~\ref{thm:sbm} is a much stronger guarantee than the weak recovery in \cite{randomized_sc}.
  Furthermore Theorem~\ref{thm:sbm} only requires at most $g = O(\log n)$ power iterations (which can be reduced further to $O(1)$ iterations whenever $n \rho_n = \Omega(n^{\beta})$ for some $\beta > 0$)
  while \cite{randomized_sc} require $g = \Omega(n^{\delta})$ power iterations for some
  fixed but arbitrary $\delta > 0$.  
\end{remark}
 
\subsection{Row-wise normal approximation}\label{sec:rownorm}
We now consider the row-wise fluctuations of $\hat{\M U}_g$.
For random graphs satisfying the assumptions at the beginning of this section, there exists a sequence of orthogonal matrices $\M W_n$ such that for any index $i \in [n]$ we have
\begin{equation}
  \label{eq:clt}
  \bm{\Gamma}_i^{-1/2} \bigl(\M W_n [\hat{\M U}]_i - [\M U]_i\bigr) \overset{\mathrm{d}}{\longrightarrow} \mathcal{N}(0, \mathbf{I})\end{equation}
as $n \rightarrow \infty$. Here $[\hat{\M U}]_{i}$ (resp. $[\M U]_i$) denote the $i$th row of $\hat{\M U}$ (resp. $\M U$) and $\bm{\Gamma}_i$ is defined as
\begin{equation}
  \label{eq:gamma_def}
 \bm{\Gamma}_i =
 \bm{\Lambda}^{-1} \Bigl(\sum_{j} m_{ij}(1 - m_{ij}) [\M U]_j ([\M U]_j)^{\top}\Bigr) \bm{\Lambda}^{-1}
 \end{equation}
where $\{m_{ij}\}$ are the entries of $\M M$ and $\bm{\Lambda}$ are the non-zero eigenvalues of $\M M$;
we note that $\|\bm{\Gamma}_i^{-1/2}\| \asymp n \rho_n^{1/2}$ for all $i$.
Eq.~\eqref{eq:clt} provides normal approximations for the rows of $\hat{\M U}$
when centered around the corresponding rows of $\M U$; see e.g., 
\cite{grdpg1,cape2019signal,xie2021entrywise} for more details.  
By Corollary~\ref{co:l2inf_noise}
there exists a choice of $g = O(\log n)$ and orthogonal $\M W_g$ depending on $g$ and $n$
such that for any $i \in [n]$ we have
\begin{equation*}
  \begin{split}
  \M W_n \M W_g [\hat{\M U}_g]_i - [\M U]_i &= \M W_n (\M W_g [\hat{\M U}_g]_i - [\hat{\M U}]_i) +
                                              \M W_n [\hat{\M U}]_i - [\M U]_{i} 
                                               \\
                                              &= \M W_n
                                   [\hat{\M U}]_i - [\M U]_i + o\bigl(n^{-1} \rho_n^{-1/2}\bigr)
  \end{split}
\end{equation*}
with probability at least $1 - 5n^{-3}$. 
Combining Eq.~\eqref{eq:clt} with the above bound we obtain  the following
normal approximations for the rows of $\hat{\M U}_g$. 
\begin{theorem}\label{thm:clt}
  Let $\hat{\M M}$ be the adjacency matrix for a random graph with edge probabilities matrix $\M M$ where $\M M$
  satisfies the assumption at the beginning of this section. 
Let $\hat{\M U}_{g}$ be generated via Algorithm \ref{RSRS} with $\sk \asymp \log n$ and $g \geq \frac{2 \log n}{\log n \rho_n}$. 
Suppose $n \rho_n = \omega(\log n)$. Then there exist a sequence of orthogonal matrices $\M W_n^*$ such that for any $i\in[n]$ we have
\begin{equation}
  \label{eq:clt_graphs}
\bm{\Gamma}_i^{-1/2} \bigl(\M W_n^* \bigl[\hat{\M U}_{g}\bigr]_i - \bigl[\M U\bigr]_i
\bigr) \overset{\mathrm{d}}{\longrightarrow} {\bf\mathcal{{N}}}(\bds 0,\M I).
\end{equation}
Moreover, if $n \rho_n \succsim n^{\beta}$ for some $\beta > 0$ then Eq.~\eqref{eq:clt_graphs} holds for all $g \geq 2 + \beta^{-1}$. 
\end{theorem}
Under the strong signal regime $n \rho_n \succsim n^{\beta}$, 
the condition $g \geq 2 + \beta^{-1}$ in Theorem~\ref{thm:clt} is slightly more stringent than
$g \geq 1 +  \beta^{-1}$ in Eq.~\eqref{eq:d_twoinf_rg_part2}, and the
reason
for this discrepancy is that while $g \geq 1 + \beta^{-1}$ is sufficient for $d_{2
  \to \infty}(\hat{\M U}_g, \M U)$ to achieve the optimal error rate,
it does not guarantee that the fluctuations of
$ \M W_n^* [\hat{\M U}_g]_i - [\M U]_i$ is asymptotically
equivalent to that of $\M W_n [\hat{\M U}]_i - [\M U]_i$. See Section
\ref{sec:RLD:emp} for a numerical experiment supporting this claim. 
\section{Numerical studies}\label{sec:simu} We conduct simulations to
support our theoretical results for RSVD-based random graph inference,
and present a real-data analysis using RSVD-based PCA.
Section~\ref{sec:ptv} illustrates the subspace perturbation error
rates from Section~\ref{sec:int:network}, and
Section~\ref{sec:RLD:emp} illustrate the Gaussian approximation from
Section~\ref{sec:rownorm}.  Section~\ref{pca:realdata} analyzes an
scRNA-seq data using the RSVD-based PCA in
Section~\ref{sec:epca}. Additional numerical results, including 
simulations for missing-data PCA, exact recovery in SBM, and a real-data analysis for distance matrix completion,
are further provided in the Supplementary File.

For the general simulation setting, we consider two-blocks SBM graphs with equal sized blocks and block
probabilities matrix 
$\M B_0 =
\rho_n\left(\begin{smallmatrix} 0.8 & 0.3 \\ 0.3 & 0.8
\end{smallmatrix}\right).
$
Recall that $\rho_n$ is the sparsity scaling parameter. 
We consider three different regimes for $\rho_n$, namely $\rho_n = 1$
(dense setting), $\rho_n = 3n^{-1/3}$ (semi-sparse setting I), $\rho_n
= 4n^{-1/2}$ (semi-sparse setting II). As $\mathrm{rk}(\mathbf{B}_0) = 2$, we also set $k = 2$. 
\subsection{Phase transition}
\label{sec:ptv}
\begin{figure}[t]
  \centering
\includegraphics[width=0.45\linewidth]{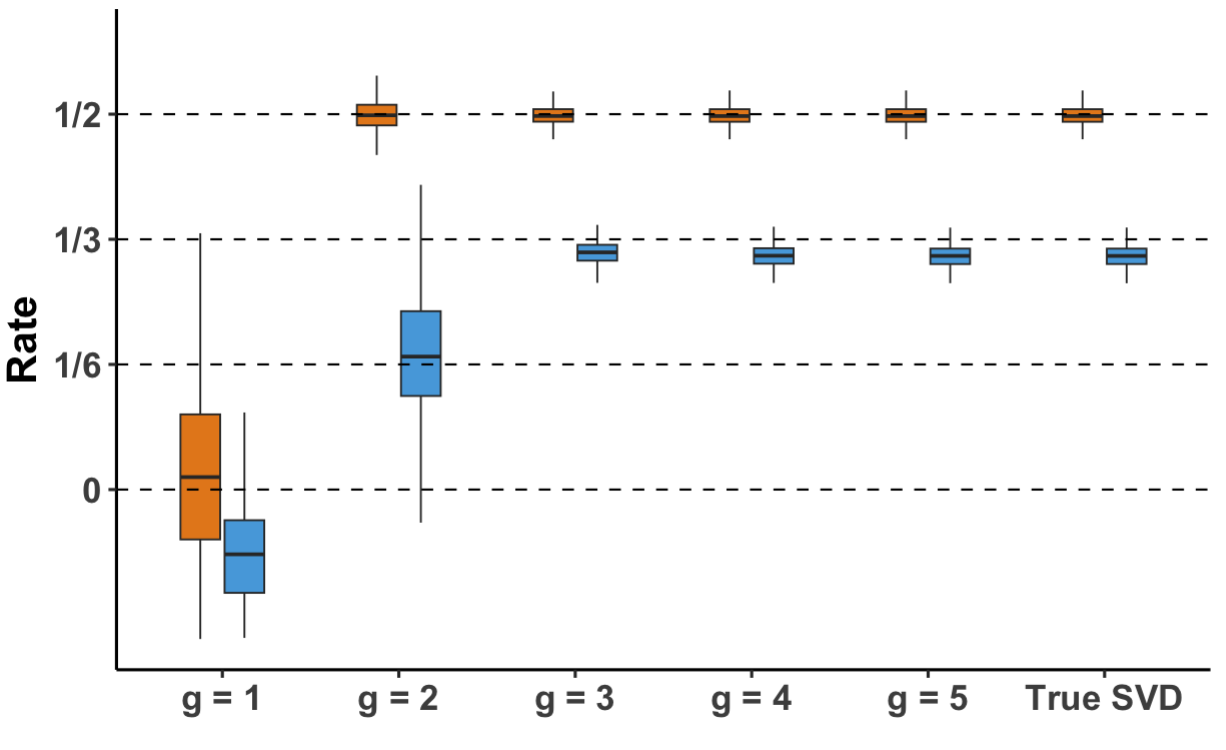} \quad\quad
\includegraphics[width=0.45\linewidth]{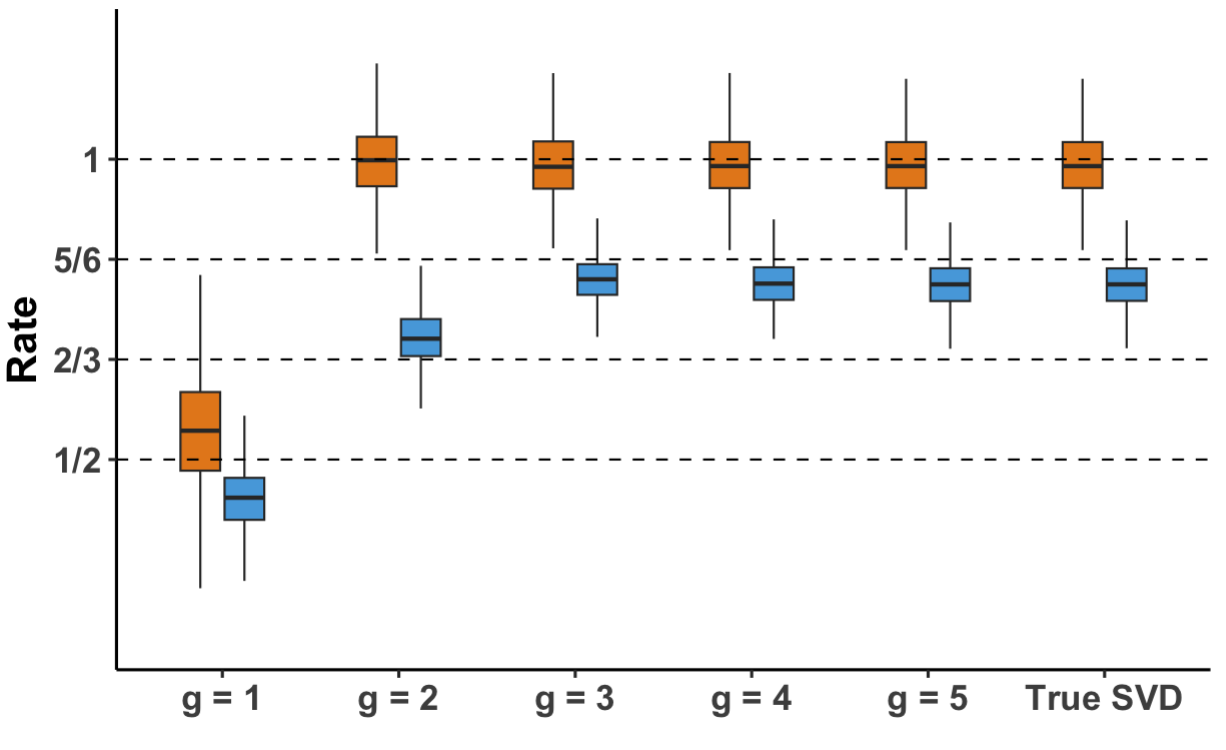}
\caption{Box plots of error rate for $d_{2}(\hat{\M U}_g, \M U)$ (top panel) and $d_{\twoinf}(\hat{\M U}_g, \M U)$ (bottom panel), where
  $\sk = 5 \log n$ and $1 \leq g \leq 5$. The colors denote different sparsity level, with $\rho_n = 1$ being blue and $\rho_n = 3n^{-1/3}$ being yellow. }\label{fig:rp}
\end{figure}

We first verify the convergence rate for $d_2(\hat{\M U}_g, \M U)$ and
$d_{2 \to \infty}(\hat{\M U}_g, \M U)$ as we vary
$g$ and $\rho_n$.
 For simplicity we only consider $\rho_n = 1$ and $\rho_n = 3n^{-1/3}$, 
and we ignore any potential logarithmic factors in the
convergence rate for $d_2(\hat{\M U}_g,\M U)$ and $d_{\twoinf}(\hat{\M
  U}_g,\M U)$ as $n$ increases. 
For each choice of $\rho_n$ we numerically
estimate the convergence rates of $d_2(\hat{\M U}_g,\M U)$
and $d_{\twoinf}(\hat{\M U}_g,\M U)$ as follows.
We first generate
a realization of $\M A$ with $n \in \{2000, 3000, \dots, 7000\}$ vertices from $\mathbf{B}_0$ with 
equal block sizes, and then
compute $\hat{\M U}_g$ via Algorithm~\ref{RSRS} with $\sk = 5 \log n$
and $1 \leq g \leq 5$. We then evaluate $d_2(\hat{\M U}_g, \M U)$
and $d_{2 \to \infty}(\hat{\M U}_g, \M U)$ and compare it against the
theoretical error rate given in Section~\ref{sec:int:network} by
running a simple linear regression between the (negative logarithm) of
the empirical error as the response variable against $\log n$ as
the predictor variable. 
The estimated coefficient $\hat{\beta}$ are then recorded. We repeat the above
steps for $500$ Monte Carlo (MC) iterations to get an empirical distribution for $\hat{\beta}$ as $g$ varies, 
and present
their box-plots 
in Figure~\ref{fig:rp}. For comparison we had also included the
(estimated) convergence rate for $d_2(\hat{\M U}, \M U)$ and
$d_{2\to\infty}(\hat{\M U}, \M U)$ (these are labeled as ``true
SVD'').  The empirical results in Figure~\ref{fig:rp}  match the theoretical rates
presented in Theorem~\ref{thm:lower} and Eq.\eqref{eq:d_twoinf_rg_part2}
exactly (see Section~\ref{rate:example}): in particular we see no convergence when $g = 1$ and $\rho_n
= 1$, slow or no convergence when $g \leq 2$ and $\rho_n = 3n^{-1/3}$, 
and asymptotically optimal convergence when
$g \geq 2$ and $\rho_n =1 $, or $g \geq 3$ and $\rho_n = 3n^{-1/3}$. 

\subsection{Limiting distribution}
\label{sec:RLD:emp} 
\begin{figure}[t]
\begin{subfigure}{0.16\textwidth}
\includegraphics[width=1\linewidth]{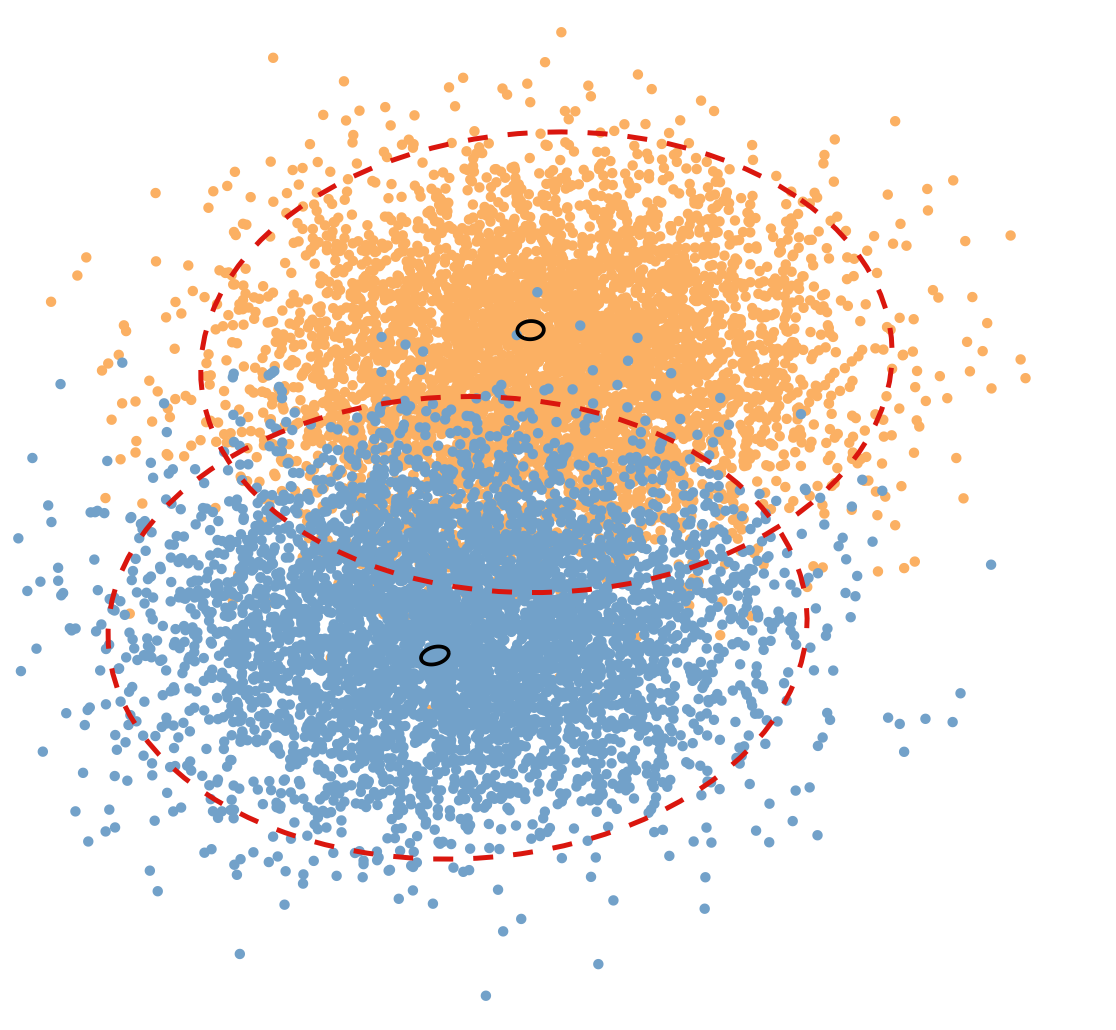}  
\includegraphics[width=1\linewidth]{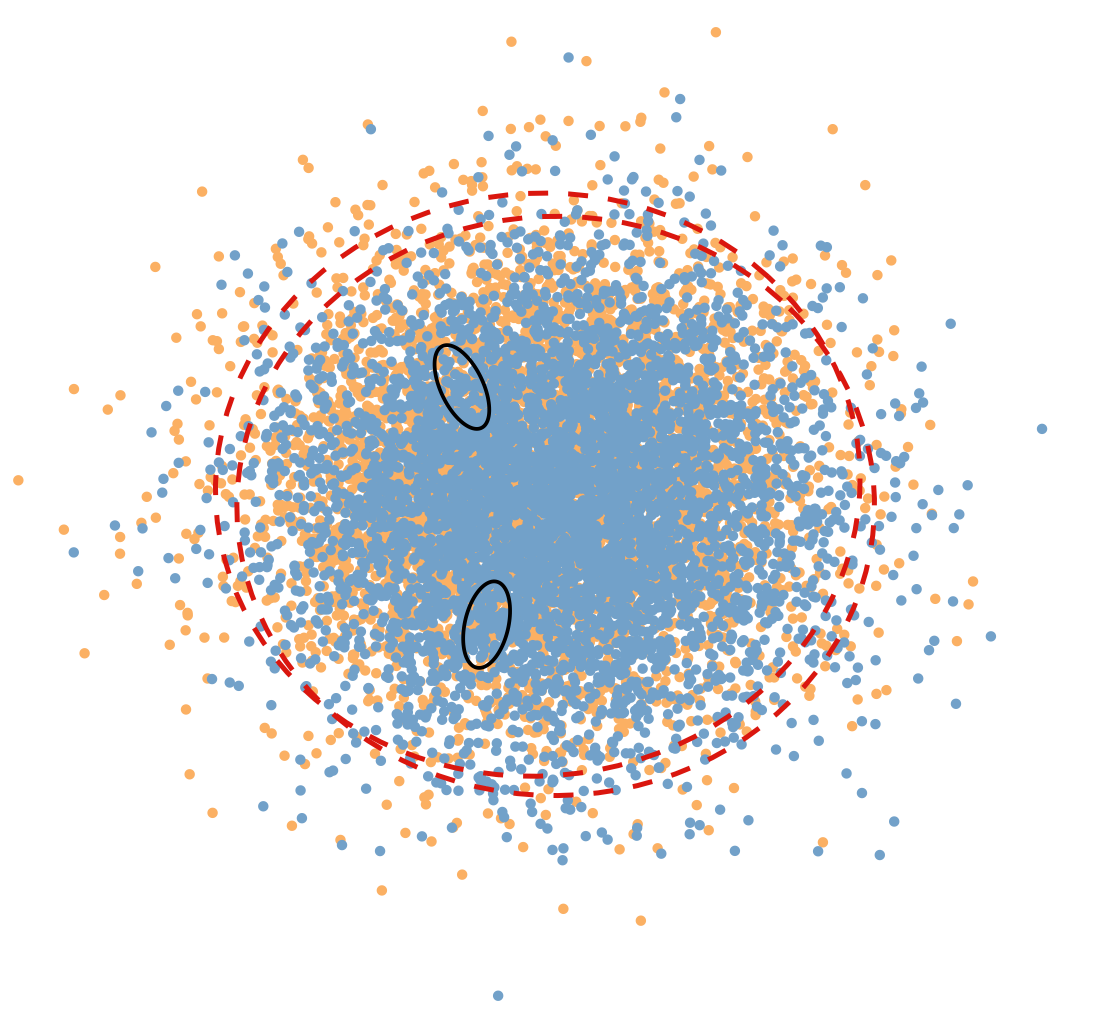} 
\caption{$g = 1$}
\end{subfigure}\begin{subfigure}{0.16\textwidth}
\includegraphics[width=1\linewidth]{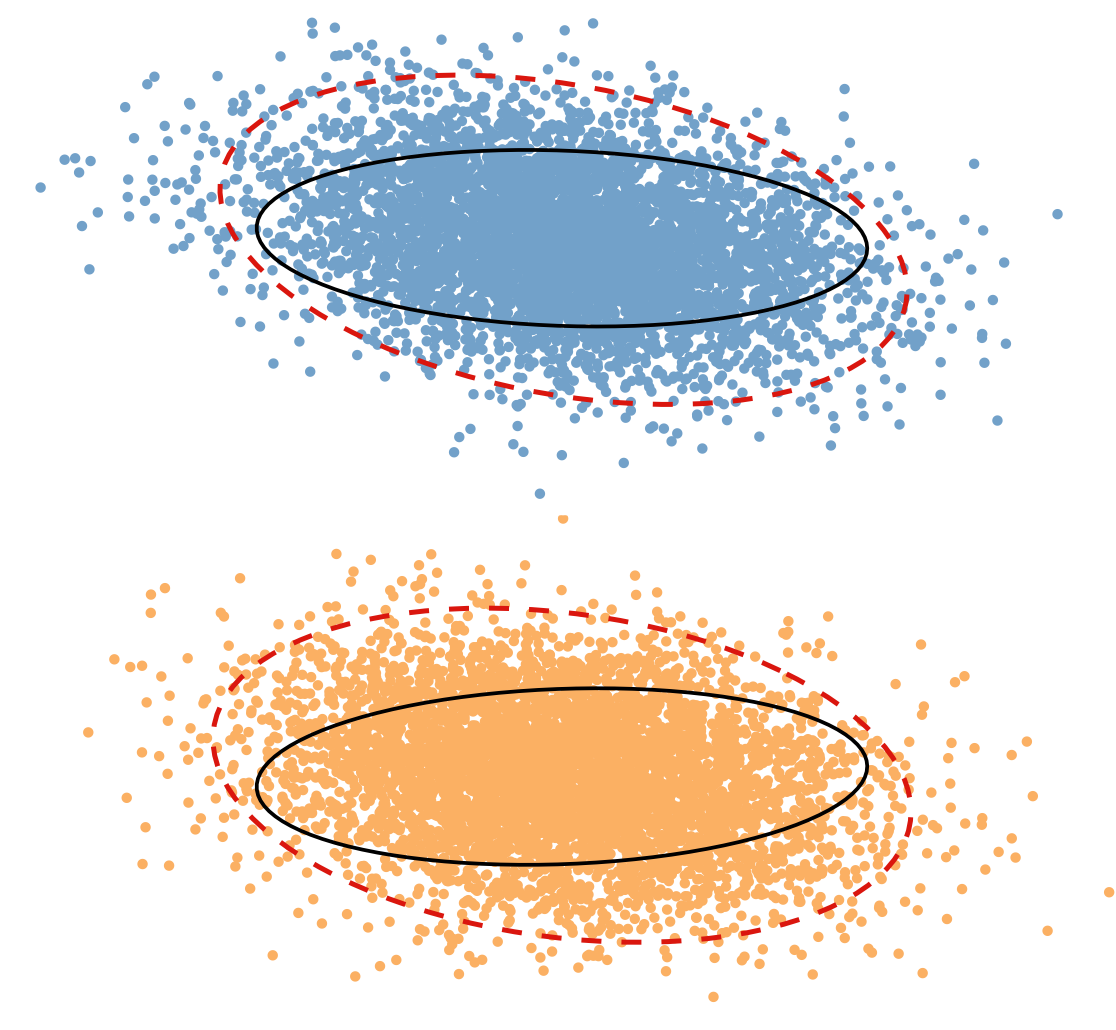} 
\includegraphics[width=1\linewidth]{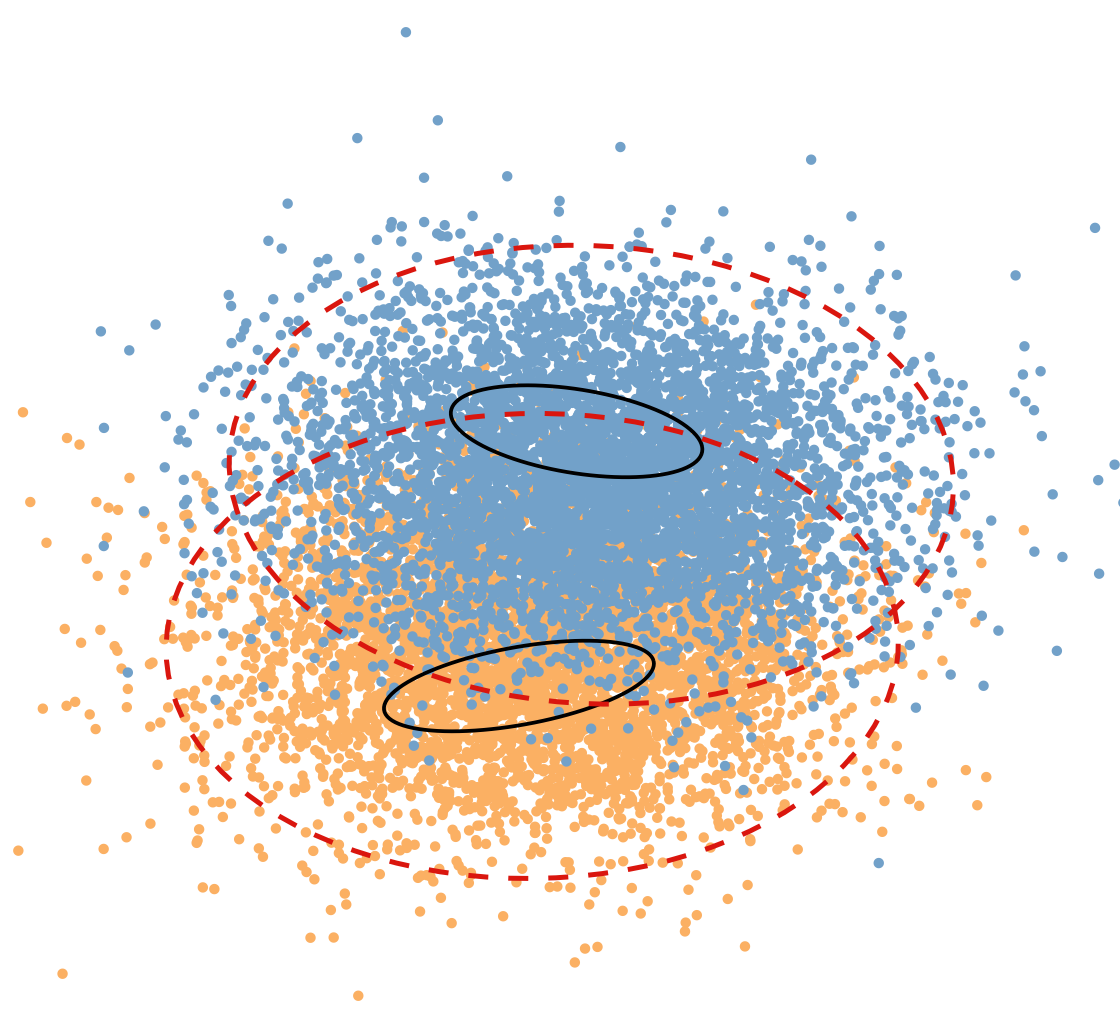} 
\caption{$g = 2$}
\end{subfigure}\begin{subfigure}{0.16\textwidth}
\includegraphics[width=1\linewidth]{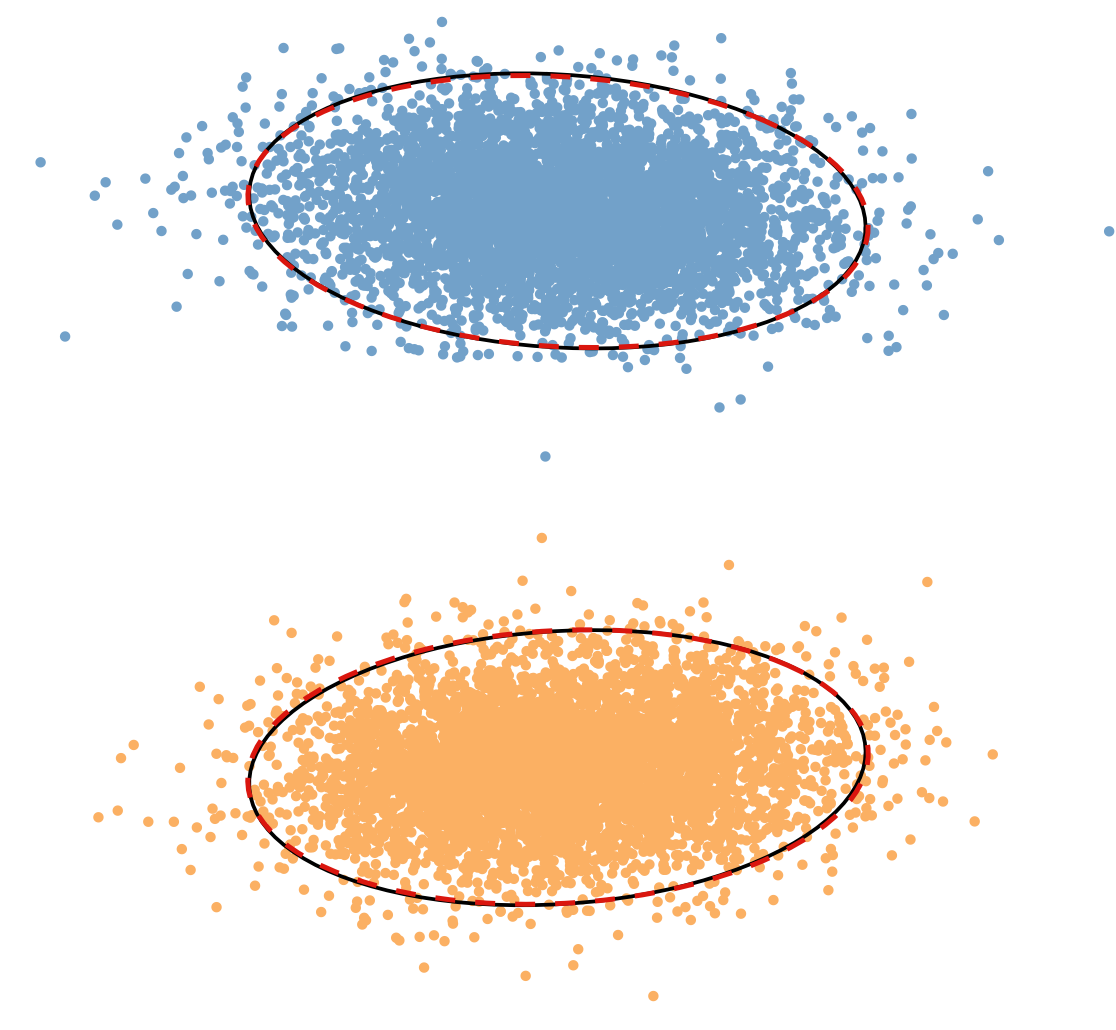} 
\includegraphics[width=1\linewidth]{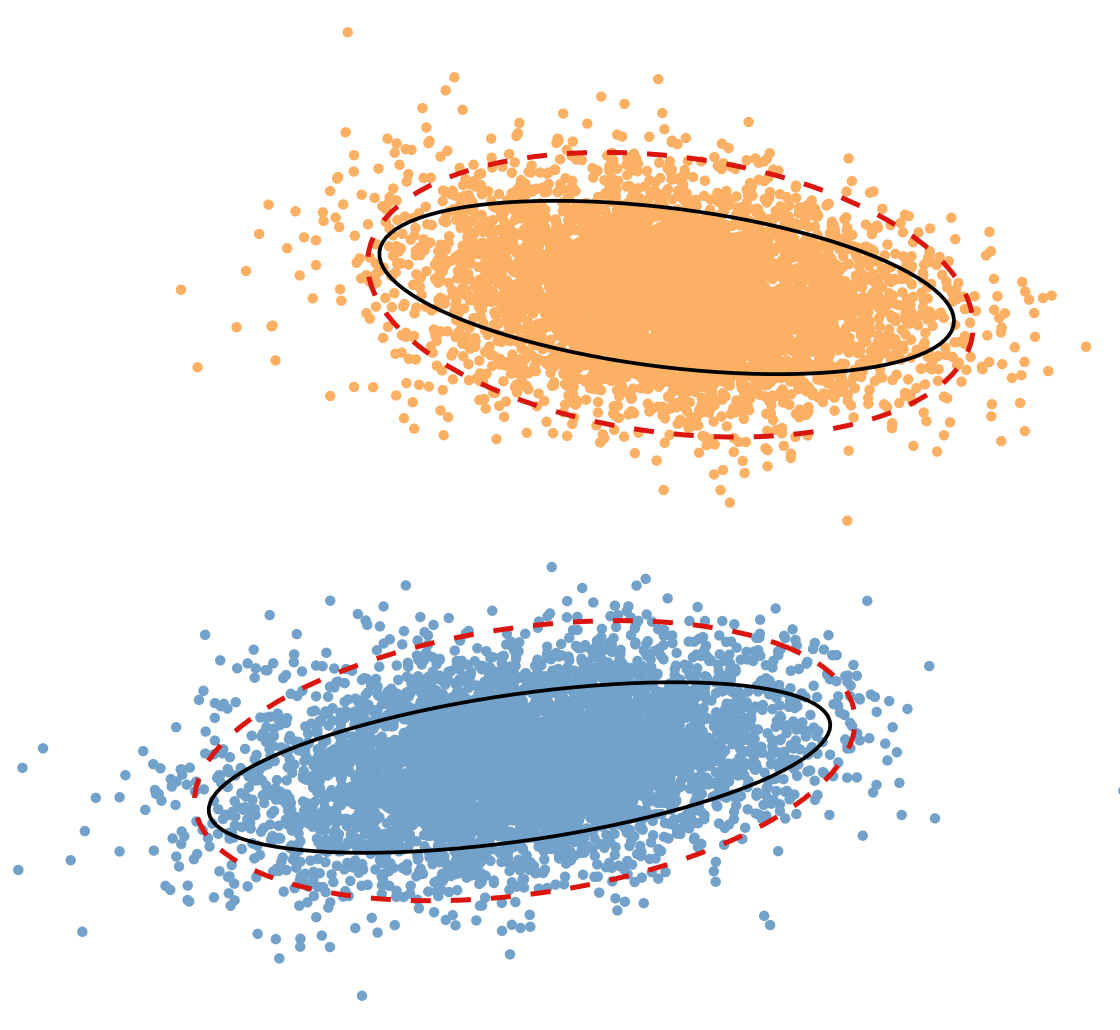}  
\caption{$g = 3$} 
\end{subfigure}\begin{subfigure}{0.16\textwidth}  
\includegraphics[width=1\linewidth]{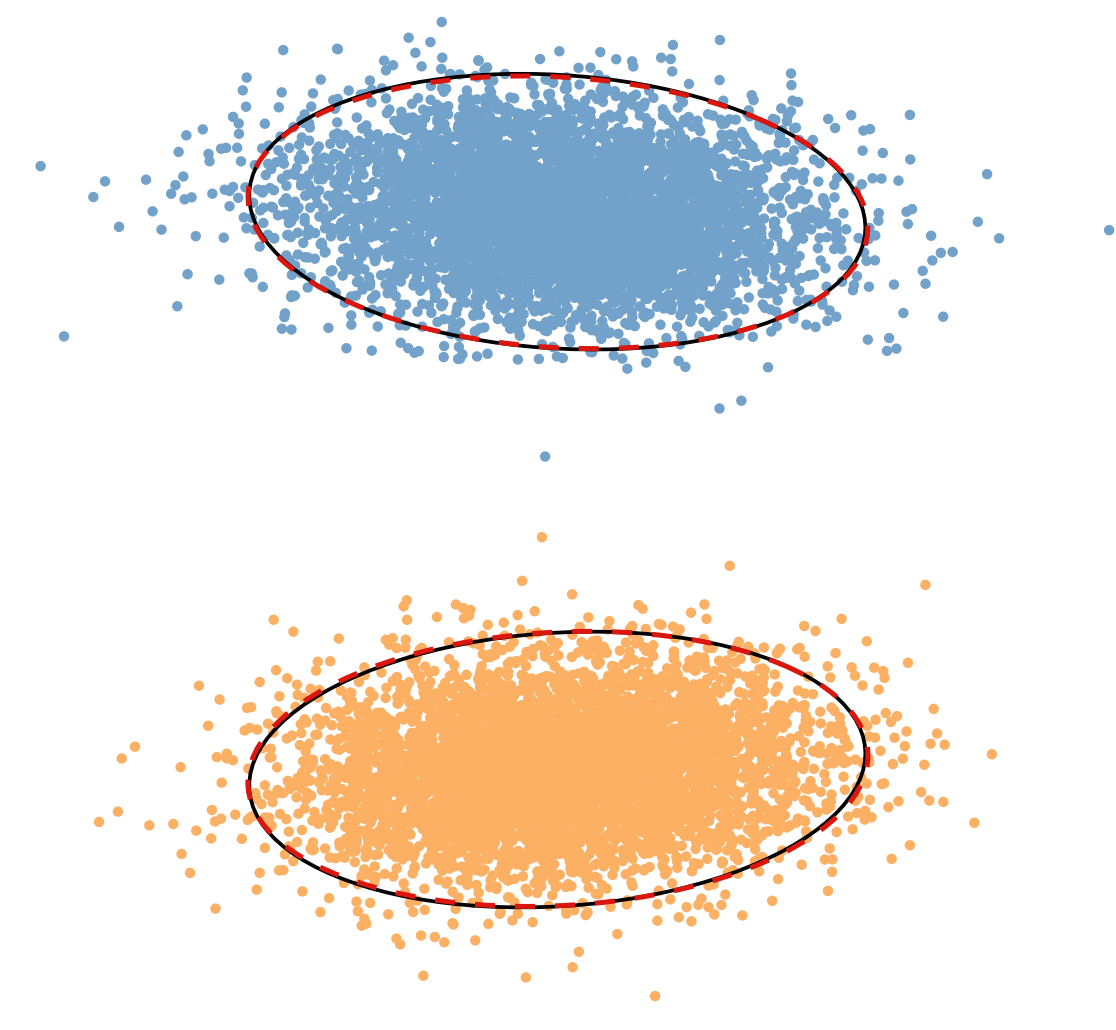} 
\includegraphics[width=1\linewidth]{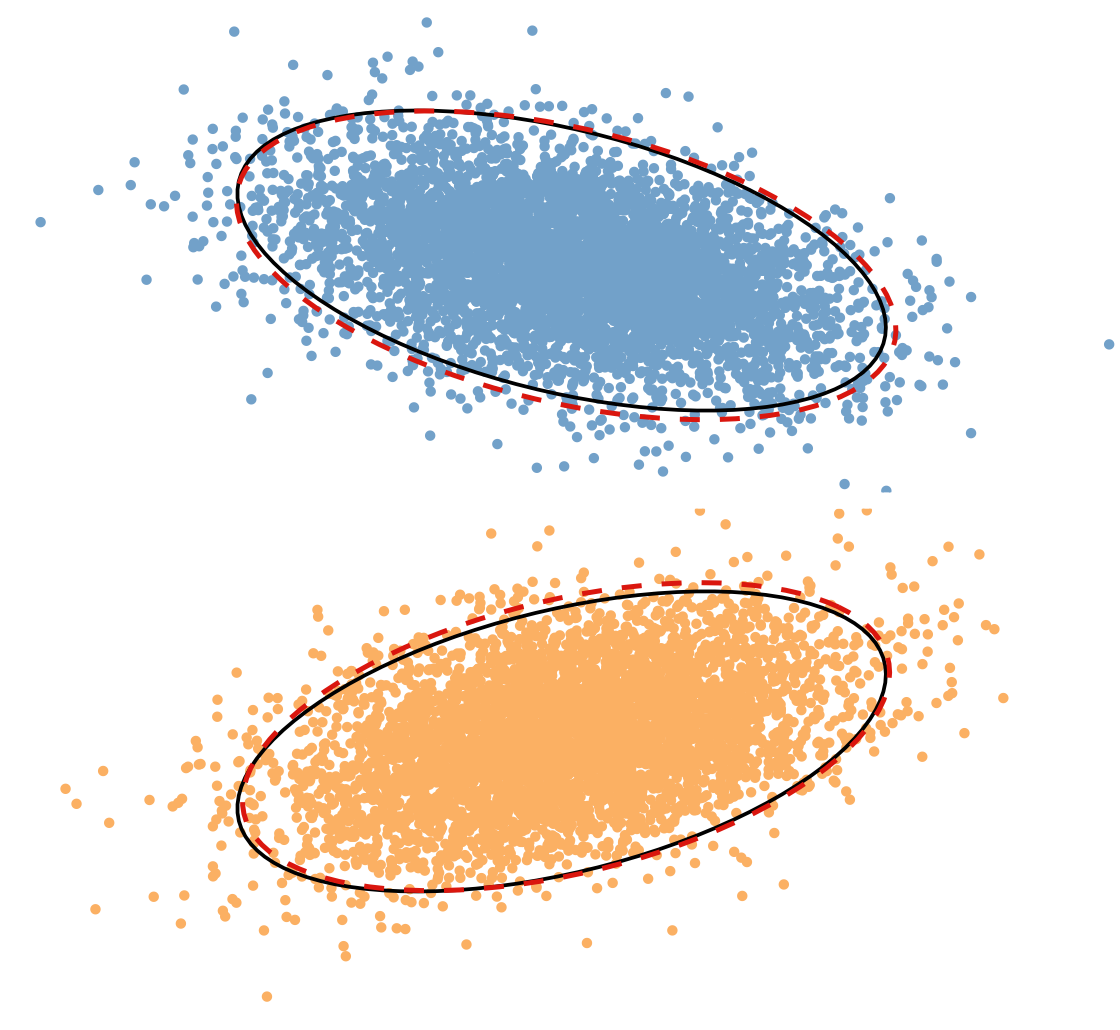} 
\caption{$g = 4$}
\end{subfigure}\begin{subfigure}{0.16\textwidth}
\includegraphics[width=1\linewidth]{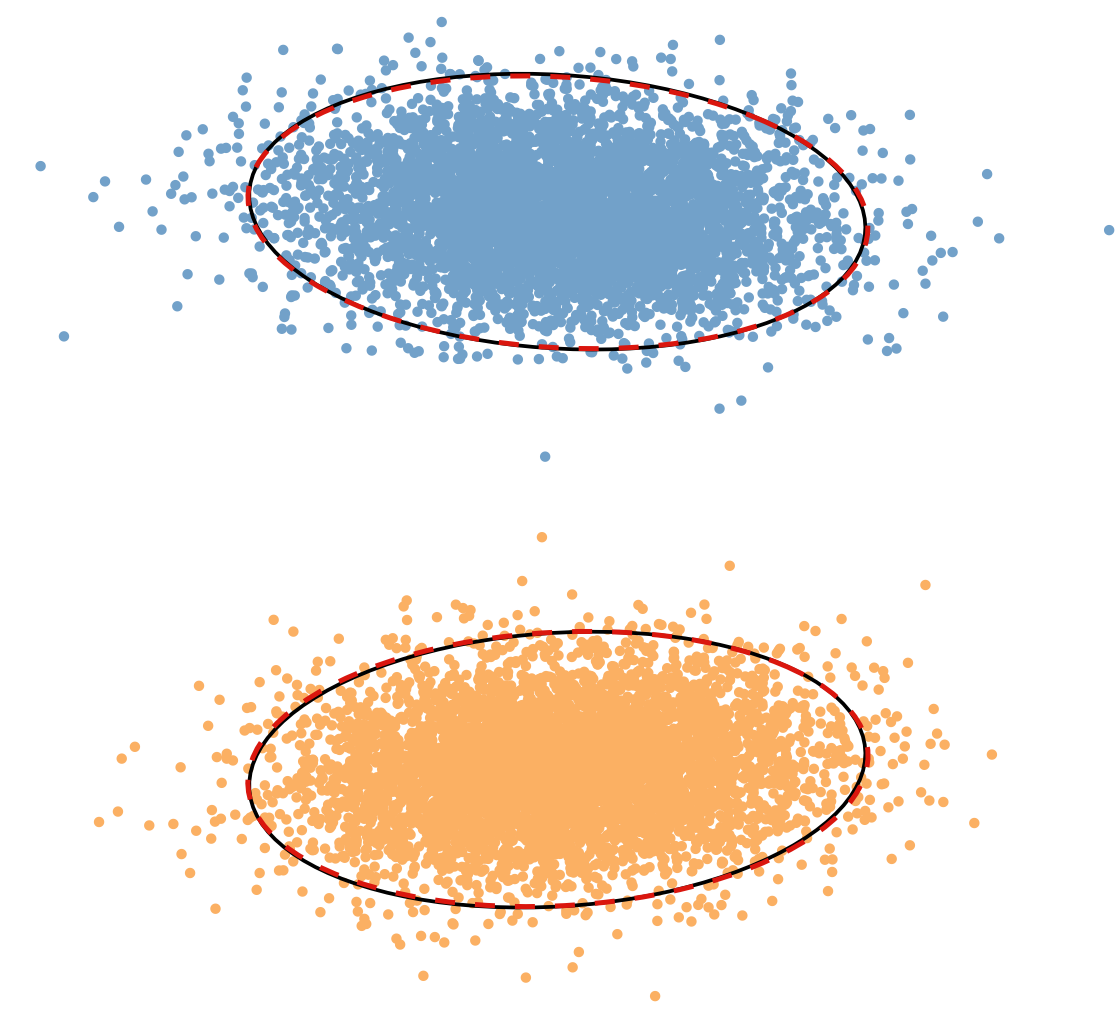} 
\includegraphics[width=1\linewidth]{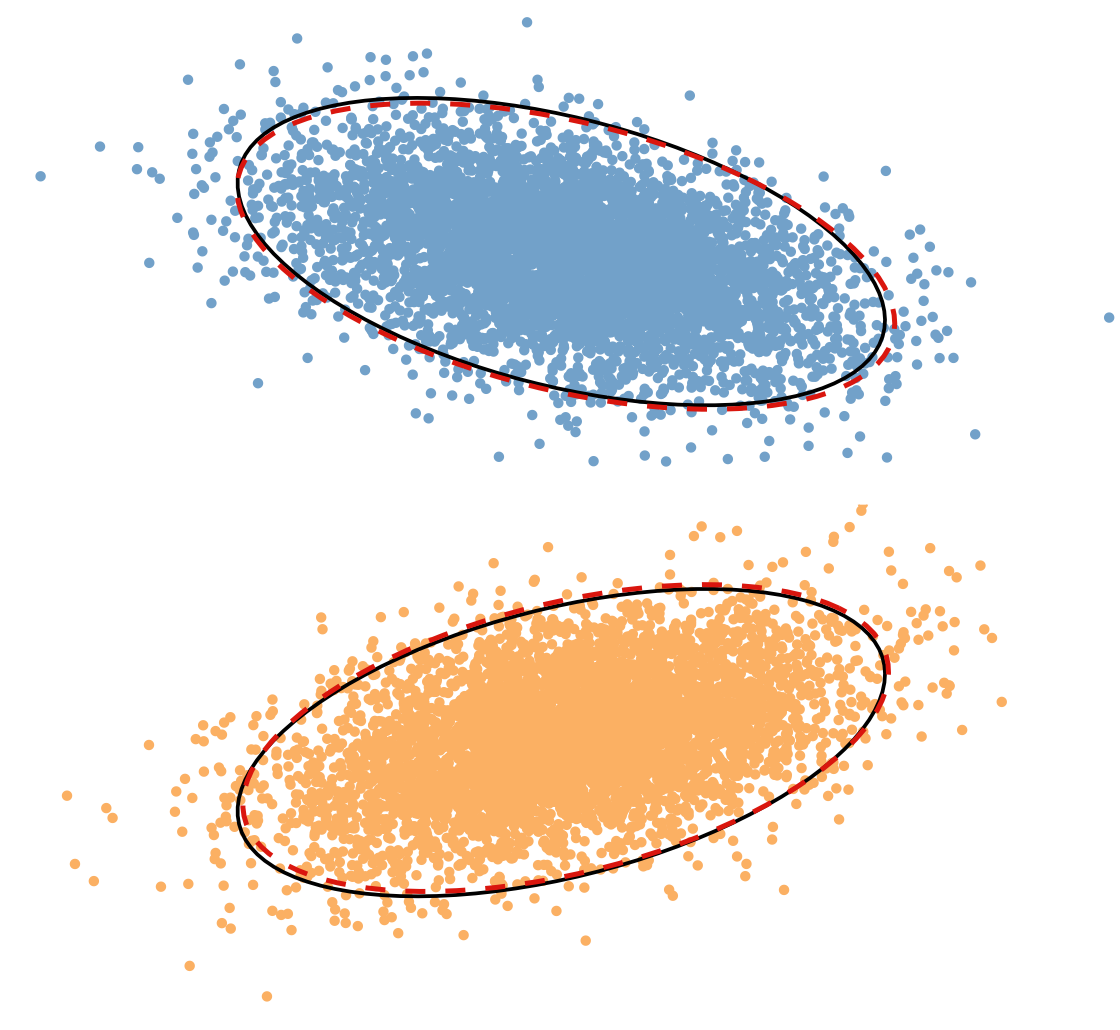} 
\caption{$g = 5$}
\end{subfigure}\begin{subfigure}{0.16\textwidth}
\includegraphics[width=1\linewidth]{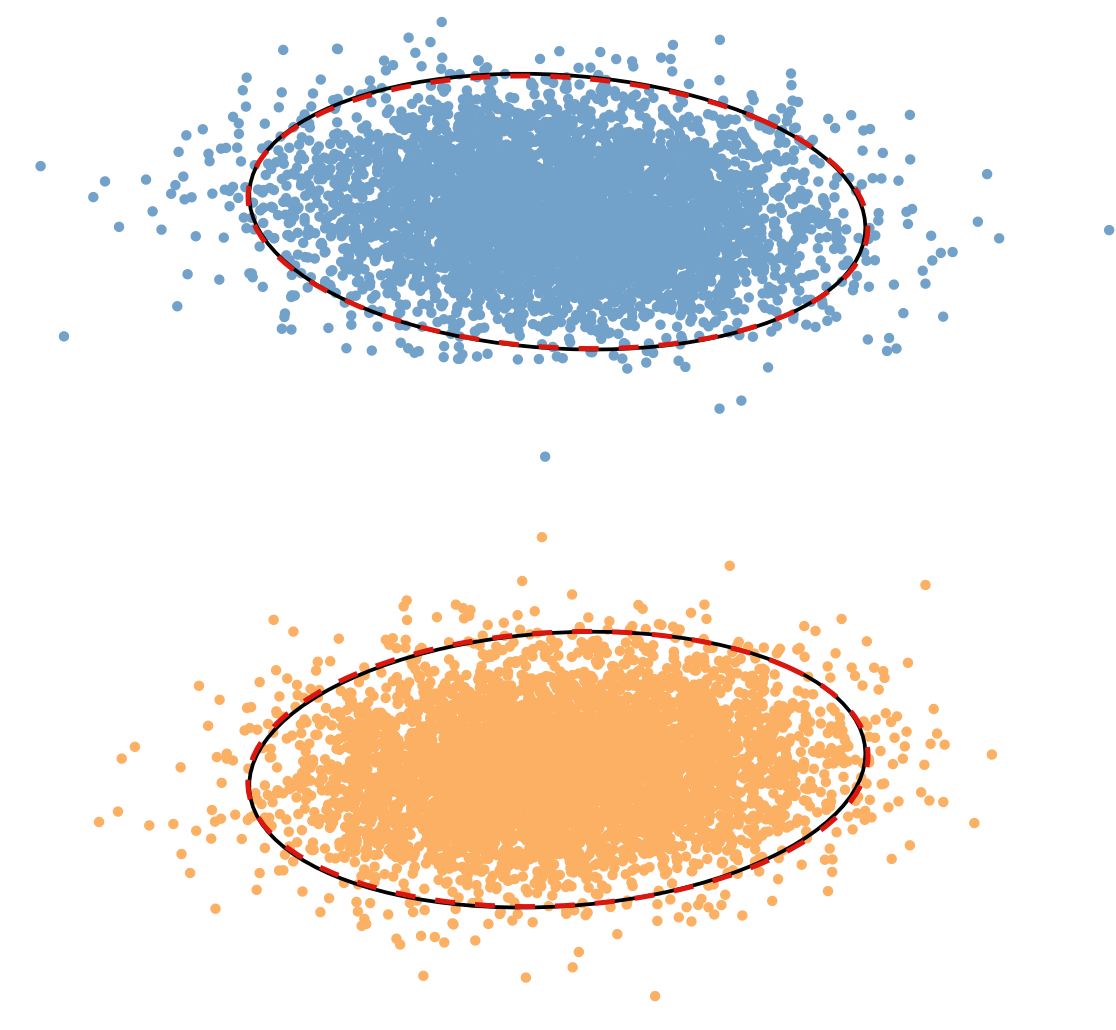} 
\includegraphics[width=1\linewidth]{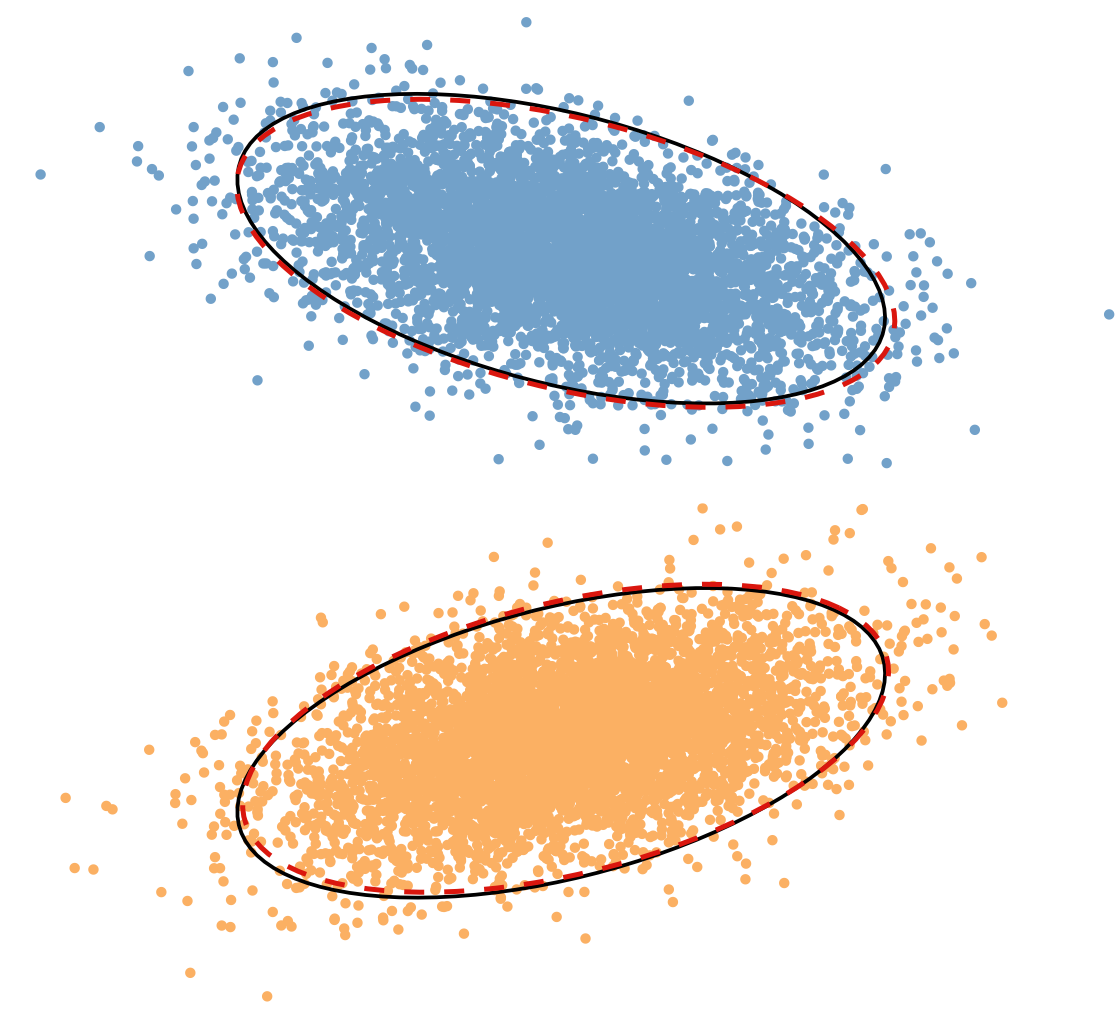} 
\caption{Exact SVD}
\end{subfigure}

\caption{Row-wise fluctuations of \(\hat{\mathbf{U}}_g\) and \(\hat{\mathbf{U}}\) for the two-block SBM with \(n = 2000\) and either \(\rho_n \asymp 1\) (top panels) or \(\rho_n \asymp n^{-1/2}\) (bottom panels). From left to right, scatter plots show the row vectors of \(\hat{\mathbf{U}}_g\) for \(g = 1, \dots, 5\) and \(\hat{\mathbf{U}}\). Points are colored according to true community memberships. Red dashed curves represent the 95\% empirical confidence ellipses; solid black curves represent the 95\% theoretical confidence ellipses.}\label{fig:cltsparse3a}
\end{figure}

We now illustrate the normal approximations for
the row-wise fluctuations of $\hat{\M U}_g$ as implied by
Theorem~\ref{thm:clt}. For brevity we set $\sk = \log n$.  We generate
a single realization of $\M A$ on $n=2000$ vertices and $\rho_n$
follows either the dense setting ($\rho_n = 1$) or the semi-sparse
setting II ($\rho_n = 4n^{-1/2}$), with equal sized blocks and block
probabilities $\mathbf{B}_0$.  Scatter plots of the rows of $\hat{\M
U}_g$ for $g \in \{1,2,\dots,5\}$ as well as the rows of $\hat{\M U}$
are presented in Figure \ref{fig:cltsparse3a}. The top panels of
Figure~\ref{fig:cltsparse3a} show that under the dense regime \(\rho_n
\asymp 1\), \(g \geq 2\) suffices for the rows of \(\hat{\M U}_g\) to
form two clusters, while \(g \geq 3\) is required for the empirical
95\% confidence ellipses to align with the theoretical ellipses from
Theorem~\ref{thm:clt}. The bottom panels show that under the
semi-sparse regime \(\rho_n \asymp n^{-1/2}\), \(g \geq 3\) suffices
for exact recovery, and \(g \geq 4\) is necessary for confidence ellipse
alignment. Moreover, for \(g \geq 4\), the scatter plots are nearly
indistinguishable from those of \(\hat{\M U}\).  These observations
confirm that the threshold \(g \geq 2 + \beta^{-1}\) in
Theorem~\ref{thm:clt} is both necessary and sufficient for the normal
approximation of \(\hat{\M U}_g\); note that \(\beta^{-1} = 1\) for
\(\rho_n \asymp 1\) and \(\beta^{-1} = 2\) for \(\rho_n \asymp
n^{-1/2}\).
 
\subsection{Subpopulation discovery of large  immune population}\label{pca:realdata}

\begin{figure}[t] 
 \includegraphics[width=1\linewidth]{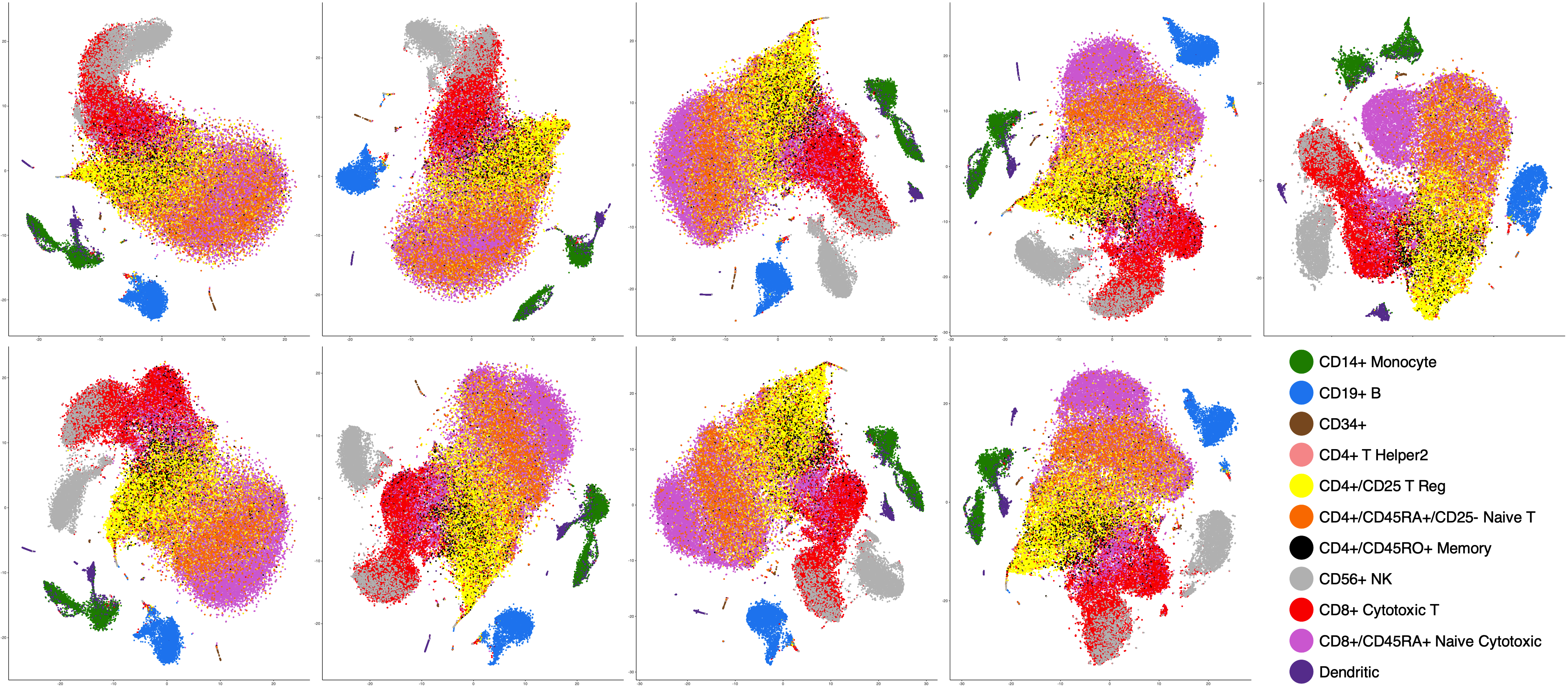}   
\caption{The $t$-SNE embeddings of the 68k PBMC gene expressions, projected on the top 50 PCs obtained via RSVD-based PCA (Algorithm~\ref{RSRS:pca}). Top row (left to right): results for $g = 1$ with $\sk = 55$, $100$, $300$, $1000$, and the reference embeddings from \citet{zheng2017massively}. Bottom row (left to right): results for $g = 2$ with the same sequence of $\sk$ values.}\label{fig:pca:real}
\end{figure} We evaluate the performance of RSVD-based PCA with
missing data (cf.~Section~\ref{sec:epca:supp}) on a widely used
single-cell RNA-seq dataset introduced by \citet{zheng2017massively},
which includes the gene expressions of approximately $6.8 \times 10^4$
peripheral blood mononuclear cells (PBMCs)
from a single donor.\footnote{The dataset is available at
\url{https://github.com/10XGenomics/single-cell-3prime-paper}.} Each
cell contains expression measurements for $\approx 2 \times 10^4$
genes, resulting in a data matrix of size $\approx (6.8
\times 10^4) \times (2 \times 10^4)$. The analysis in
\citet{zheng2017massively} first filtered the data by selecting the top
1,000 genes ranked by normalized dispersion \citep{macosko2015highly},
producing a matrix of size $\approx (6.8 \times 10^4) \times 10^3$. PCA
was then applied to this filtered and thinned matrix to obtain projections of the
cell expressions data onto the top 50 principal components (PCs). A
two-dimensional $t$-distributed stochastic neighbor embedding
($t$-SNE) of the projected data, colored by inferred cell types
based on the reference profiles from 11 purified PBMC subpopulations
\citep[Supplementary Figure 7]{zheng2017massively}, is shown in the
top-right panel of Figure~\ref{fig:pca:real}.

As the original dataset is quite large of size $\approx (6.8 \times 10^4) \times (2 \times 10^4)$, 
direct computation of its SVD (as required by traditional PCA or missing-data PCA of
\citet{cai2021subspace}) is infeasible using the memory capacity of a standard
laptop (e.g., a MacBook with 36GB of RAM). The filtering
step in \citet{zheng2017massively} directly reduces the data
dimension by dropping a large amount of the gene expression
measurements, thereby easing the memory and computational burden for
PCA. As an alternative, we apply our RSVD-based PCA with missing data
(Algorithm~\ref{RSRS:pca}) to the {\em whole} data matrix,
enabling efficient and scalable approximations to the top PCs. 
More specifically we project the whole data onto the leading $k = 50$
approximated PCs (computed for various choices of $\sk$ and $g$) 
and visualize the resulting low-dimensional
embeddings via $t$-SNE. The results are presented in Figure~\ref{fig:pca:real} and
Figure~\ref{fig:pca:real:supp}, where cell types are labeled using 
the same classification scheme as that in \citet{zheng2017massively}. We
observe that if $\sk$ and $g$ are both small then 
the $t$-SNE embeddings exhibit weaker separation among subpopulations. For
instance, when $g = 1$ and $\sk = 55$, CD56$^+$ NK cells (grey) are
not well-separated from CD8$^+$ cytotoxic T cells (red). When $g = 1$
and $\sk = 55$ or $100$, CD8$^+$/CD45RA$^+$ T cells
(light purple) are indistinguishable from CD4$^+$/CD45RA$^+$/CD25$^-$
T cells (orange). In contrast, as $\sk$ and $g$ increase, the
RSVD-based PCA closely resembles the original analysis from
\citet{zheng2017massively}, recovering meaningful cell
subpopulations. Figure~\ref{fig:pca:real:supp} shows that while the
embeddings do vary for different values of $g$ and $\sk$
when $g\geq 3$, the overall subpopulation clustering patterns remain
consistent. 
These findings align with our theoretical insights, and
demonstrate the computational and statistical efficiency of RSVD   for 
large-scale scRNA-seq analysis.

\bibliographystyle{chicago}
\bibliography{ref}
\newpage
\renewcommand{\thetheorem}{S\arabic{theorem}}
\renewcommand{\thelemma}{S\arabic{lemma}}
\renewcommand{\theremark}{S\arabic{remark}}
\renewcommand{\thealgocf}{S\arabic{algocf}}
\renewcommand{\thesection}{S\arabic{section}}
\renewcommand{\theco}{S\arabic{co}}
\renewcommand{\thetable}{S\arabic{table}}
\renewcommand{\thefigure}{S\arabic{figure}}
{{\bf\Large Supplementary File for ``Perturbation Analysis of Randomized SVD and its Applications  to Statistics''}}
\bigskip
 \setcounter{section}{0}

\def\spacingset#1{\renewcommand{\baselinestretch}%
{#1}\small\normalsize} \spacingset{1}
\spacingset{1.83}
 Section~\ref{sec:add:rg} complements Section~\ref{sec:rgi} by providing additional discussion and numerical results for random graph inference. 
In Sections~\ref{sec:mc} and~\ref{sec:epca:supp}, we apply our general theoretical framework to two additional inference problems: matrix completion with noise and PCA with missing data, and present the corresponding theoretical results, computational refinement, and numerical experiments.
Section~\ref{sec:dist} introduces the application of RSVD to distributed estimation in multi-layer networks, along with the associated theoretical guarantees. Section~\ref{sec:discussion} discusses several directions for future research.
Sections~\ref{sec:proof_thm} and~\ref{sec:pfa} contain all technical proofs. 
 \section{Additional results for   random graph inference}\label{sec:add:rg}
 \subsection{Examples of error rates under different parameter regimes}\label{rate:example}
We provide a few additional examples of the error rates presented in
Section~\ref{sec:int:network} for random graph inference under
different parameter regimes. These examples illustrate the impact of
$g$, $n, \sk$, and $n \rho_n$ on the convergence rates for
$d_{2}(\hat{\M U}_g, \M U)$ and $d_{\twoinf}(\hat{\M U}_g, \M U)$. 
see also the visual summary in
Figure \ref{fig:phase}.  For ease of exposition we will ignore
all factors depending on $\log n$ when discussing the convergence rates for 
$d_{2}(\hat{\M U}_g, \M U)$ and $d_{\twoinf}(\hat{\M U}_g, \M U)$. 
We also assume, unless stated
otherwise, that $n \rho_n \asymp n^{\beta}$ for some $\beta \in
(0,1]$. This corresponds to graphs where the average degree is of
order $\Theta(n^{\beta})$.
\begin{enumerate}
\item[(i)] Suppose $\sk \asymp \log n$. Then the threshold for convergence of $d_{2}(\hat{\M U}_g, \M U)$ and 
$d_{\twoinf}(\hat{\M U}_, g \M U)$ are the same, and are given by $\alpha_* = \frac{\log n}{\log (n \rho_n)}$. In particular: 
\begin{itemize}
\item if $n \rho_n \asymp n^{1/2}$ then $g \geq 2$
  is sufficient for $d_2(\hat{\M U}_g, \M U)$ and
$d_{\twoinf}(\hat{\M U}_g, \M U)$
to attain the optimal rate of $n^{-1/2}$ and $n^{-1}$, respectively;
no convergence is guaranteed when $g = 1$.
\item if $n \rho_n \asymp n^{2/3}$ then for $d_{2}(\hat{\M U}_g,\M U)$, the optimal rate
 $n^{-1/3}$ is attained when $g\geq
3$ and a sub-optimal rate  $n^{-1/6}$ is attained when $g = 2$. For $d_{\twoinf}(\hat{\M U}_g,\M U)$, the optimal rate
 $n^{-5/6}$ is attained when $g\geq
3$ and a sub-optimal rate $n^{-2/3}$ is attained when $g = 2$. No convergence is guaranteed when $g = 1$. 
\item If $n \rho_n \asymp n^{1/2}$ then for $d_{2}(\hat{\M U}_g,\M U)$, the optimal rate
 $n^{-1/4}$ is attained when $g \geq 3$. For $d_{\twoinf}(\hat{\M U}_g,\M U)$, the optimal rate
 $n^{-3/4}$ is attained when $g \geq 3$. No convergence is guaranteed when $g \leq 2$. 
\end{itemize}
\item[(ii)] Next suppose $\sk = \Omega(n)$.
  Then $d_{2}(\hat{\M U}_g, \M U)$ attains the optimal rate of $n^{-\beta/2}$  for all $g \geq 1$ and
  for any $\beta > 0$. In contrast, $d_{\twoinf}(\hat{\M U}_g, \M U)$
  attains the optimal rate of $n^{-(1+\beta)/2}$ for $g \geq 1$ and any $\beta \geq 1/2$,
  but a possibly slower rate when $g = 1$ and $\beta < 1/2$, i.e., 
  \begin{itemize}
  \item if $3a+2 < \beta^{-1} < 3a+3$ for some integer $a \geq 0$ then $d_{\twoinf}(\hat{\M U}_g, \M U)$
  attain the optimal rate of $n^{-(1 + \beta)/2}$ when $g \geq a+2$, a slower rate of $n^{-3(a+1)\beta/2}$ when $g = a+1$, and might not converge when $g = a$.
  \item if $3a \leq \beta^{-1} \leq 3a+2$ for some integer $a \geq 0$ then $d_{\twoinf}(\hat{\M U}_g, \M U)$ attain the optimal rate of $n^{-(1 + \beta)/2}$ when  $g \geq a+1$ but might not converge when $g \leq a$. 
  \end{itemize}
\end{enumerate} 
In summary, the above discussions provide further evidence to the
well-known advice that choosing a slightly larger $g$ and $\sk$ are
essential to the success of RSVD in practical applications \citep{martinsson_rnla}.
\subsection{Additional numerical results for exact recovery} We
used the same simulation setting as
that described in Section \ref{sec:ptv} in the main paper, with $n \in
\{1000,\dots,$ $5000\}$ and $\rho_n \in \{1, 3n^{-1/3}, 4n^{-1/2}\}$.  
For each combination of $n$ and $\rho_n$, we generate
an adjacency matrix $\M A$ with equal sized blocks and block
probabilities $\mathbf{B}_0$. 
We then perform RSVD-based spectral clustering
with $g \leq 3$, where we use the \texttt{Gmedian}
library in $\mathbf{R}$ to perform fast (approximate) $K$-medians clustering on
the rows of $\hat{\M U}_g$. For comparison, we also use
\texttt{Gmedian} to perform $K$-medians clustering on the rows of
$\hat{\M U}$. We repeat the above steps for $500$ Monte Carlo
replicates. The proportion of times (among these $500$ replicates) in
which the $K$-medians clustering of either $\hat{\M U}_g$ or $\hat{\M
  U}$ correctly recover the community memberships for {\em all} nodes are reported in 
Table~\ref{tb:2}.
\begin{table}[t]\linespread{1} 
\setlength{\tabcolsep}{14pt}
  \footnotesize
\caption{Proportion of times that RSVD-based spectral clustering (g =
  1,2,3) or the original spectral clustering (OSC) exactly recover the
  memberships of all nodes, among $500$ MC rounds; here
  $\tilde{k} = 5 \log n$. Standard errors are reported in parentheses.}
\label{tb:2}
\centering
\begin{tabular}{c|ccccc}
\hline
Sparsity&$n$&$g = 1$&
$g = 2$&
$g = 3$&
OSC\\
\hline
 & 1000 & 0.536 {(0.022)}  &  $0.384$  {(0.000)} &  $1.000$ {(0.000)} &  $1.000$  {(0.000)}
\\
& 2000 & 0.616 {(0.022)} &  $1.000$ {(0.000)}  &  $1.000$ {(0.000)}  &  $1.000$ {(0.000)}  
\\
 $\rho_n \asymp 1$& 3000 & 0.674 {(0.021)}   &  $1.000$ {(0.000)}  &  $1.000$ {(0.000)} &  $1.000$ {(0.000)}   
 \\
&  4000 & 0.600 {(0.022)}   &  $1.000$ {(0.000)} &  $1.000$ {(0.000)} &   $1.000$ {(0.000)}  
 \\
 & 5000 & 0.626 {(0.022)}   &  $1.000$ {(0.000)}&  $1.000$ {(0.000)} &  $1.000$ {(0.000)}    
 \\
\hline
 & 1000 & 0.000 {(0.000)}  &  $0.982$ {(0.012)} &   $1.000$ {(0.000)} &  $1.000$ {(0.000)}
\\
&  2000 & 0.000 {(0.000)} &  $0.996$ {($0.006$)}  &  $1.000$ {(0.000)}  &  $1.000$ {(0.000)}    
\\
$\rho_n\asymp n^{-1/3}$ & 3000 & 0.000 {(0.000)}   &  $1.000$ {($0.000$)}  &  $1.000$ {(0.000)}  &  $1.000$ {(0.000)}   
 \\
&  4000 &0.000 {(0.000)}   &  $1.000$ {($0.000$)} &  $1.000$ {(0.000)} &  $1.000$ {(0.000)}  
 \\
&  5000 & 0.000 {(0.000)}   &  $1.000$ {($0.000$)}&  $1.000$ {(0.000)} &  $1.000$ {(0.000)}   
 \\
\hline
& 1000 & 0.000 {(0.000)}  &  0.000 {(0.000)} &   $0.604$ {(0.022)}&  $0.914$  {(0.013)}
\\
 & 2000 & 0.000 {(0.000)} &  0.000 {(0.000)}  &  $0.892$ {(0.014)}  &  $1.000$ {(0.000)}   
\\
$\rho_n\asymp n^{-1/2}$   &3000 & 0.000 {(0.000)}   &  0.000 {(0.000)}  &  $0.974$ {(0.007)}  &  $1.000$ {(0.000)}   
 \\
  &4000 &0.000 {(0.000)}   &  0.000 {(0.000)} &  $0.996$ {(0.003)} &  $1.000$ {(0.000)}  
 \\
 & 5000 & 0.000 {(0.000)}   &  0.000 {(0.000)}&  $0.998$ {(0.002)} &  $1.000$ {(0.000)}   
  \\
\hline
\end{tabular}
\end{table}
From 
Table \ref{tb:2} we see that 
$K$-means clustering on $\hat{\M U}_g$ exactly recover the underlying community
assignment with probabilities converging to $1$ as $n$ increases,
provided that $g \geq 2$ when $\rho_n \in \{1,3n^{-1/3}\}$ and $g \geq
3$ when $\rho_n = 4n^{-1/2}$. 
These empirical results are consistent with the theoretical results
presented in Theorem~\ref{thm:sbm},
i.e., exact recover is
guaranteed if and only if $g > \beta^{-1}$. Note that $\beta^{-1} = 1$
when $\rho_n = 1$, $\beta^{-1} = \tfrac{3}{2}$ when $\rho_n \asymp
n^{-1/3}$, and $\beta^{-1} = 2$ when $\rho_n \asymp n^{-1/2}$.
Finally we note that if $g = 1$ and $\rho_n = 1$ then RSVD-based spectral clustering
  occasionally recovers the community membership for all nodes. 
This is due to the fact that although $d_{\twoinf}(\hat{\M U}_1,\M U) =
\Omega(n^{-1/2})$ with high probability, it is still possible that
$d_{\twoinf}(\hat{\M U}_1,\M U) \leq \zeta_*$ where $\zeta_*$ is the
minimum $\ell_2$ distance between any two nodes $i$ and $j$ belonging to different
communities; exact recovery is certainly expected if $d_2(\hat{\M
  U}_1, \M U) \leq \zeta_*$. 
 \section{Matrix completion with noises}\label{sec:mc}%
\subsection{Theoretical results}
Let $\M T\in\RR^{n\times n}$ be a matrix whose entries are only partially and
noisily observed. Such matrix occurs in many real-world
applications, including the well-known Netflix challenge. 
As another example, if $\mathbf{T}$ is an Euclidean distance matrix (EDM) between $n$
points in $\mathbb{R}^{d}$ then $\mathrm{rk}(\M T) \leq d+2$ and it is commonly the case that
$\mathbf{T}$ is noisily observed \citep{Javanmard2013}; similarly, if
$\mathbf{T}$ is a signal correlation matrix between multiple remote sensors
then $\mathbf{T}$ is partially observed
due to power constraints \citep{cheng2012stcdg}. 
Assume, for the current discussion, that $\M T$ is symmetric and we observed
\bee\label{T:miss}
\hat{\M T} = {\mathcal{P}}_{\M \Omega}(\M T + \M N) := \bm{\Omega} \circ (\mathbf{T} + \M N),
\ee
where $\M N$ denote an unobserved symmetric $n\times n$ noise matrix, 
$\bm{\Omega}$ is a symmetric matrix with $\{0,1\}$ entries, and $\circ$ denote the
Hadamard product.
We shall assume, for ease of exposition, that the (upper triangle)
entries of $\M N$ are iid $\mathcal{N}(0, \sigma^2)$ random
variables while the (upper triangular) entries of $\M \Omega$ are iid
Bernoulli random variables with success probability $p$.  
As $\E[p^{-1}\hat{\M T}] = \M T$, one simple and widely used
estimate for $\M T$ is given by $p^{-1} \hat{\M T}^{(k)}$ where
$\hat{\M T}^{(k)}$ is the truncated rank-$k$ SVD of $\hat{\M T}$ for some choice of $k$; see \cite{abbe2020entrywise,
  chen2021spectral,chatterjee2015matrix} and the
references therein. 

In many real-world applications, the dimensions of $\hat{\M T}$ can be rather large and yet $\hat{\M
  T}$ can be quite sparse compared to $\M T$, i.e., the number of
non-zero entries of $\hat{\M T}$ is much smaller than $n^2$. It is
thus computationally attractive to approximate the left singular
vectors of $\hat{\M T}$ using randomized SVD.
More specifically, let $\hat{\M U}_g$ be the output of Algorithm
\ref{RSRS} with $\hat{\M M} = \hat{\M T}$ for some 
choices of $k, \sk$ and $g$.
Given $\hat{\M U}_g$ we compute a rank-$k$ approximation for $\hat{\M
  T}$ via $\hat{\M U}_g \hat{\M U}_g^{\top} \hat{\M T}$. We 
can then take $\hat{\M T}_g  := p^{-1}\hat{\M U}_g \hat{\M U}_g^{\top} \hat{\M T}$ or $2^{-1} (\hat{\M T}_g + \hat{\M T}_g^{\top})$ as
an estimate for $\M T$. 

We now combine the $\ell_{\twoinf}$ perturbation and entrywise
concentration bounds in Corollary~\ref{co:l2inf_noise} of our paper
with Theorem~3.4 of \cite{abbe2020entrywise} 
to obtain error bounds for $\hat{\M U}_g$ and
$\hat{\M T}_g$ as estimates for $\M U$ and $\M M$, respectively. 
For ease of exposition we shall assume that $p$ is known.
If $p$ is unknown then, as the entries
of $\M T$ are assumed to be missing completely at random, 
it can be consistently estimated from the proportion
of observed entries in $\hat{\M T}$. The resulting $\hat{p}$ converges to $p$ at rate $n^{-1}p^{-1/2}$ and 
has no effect on the theoretical results. 
\begin{theorem}\label{thm:mc} Let $\M T$ be a symmetric $n \times n$
  matrix and denote $k_0:= \mathrm{rk}(\M T)$.
  Let $\hat{\M T}$ be a noisily observed version of $\M T$ sampled according to \Eq\eqref{T:miss} for some known value of
  $p \in (0,1)$. Let $\lambda_i(\M T)$ denote the $i$th largest eigenvalue (in modulus) of $\M T$. 
  Define $E_n = (n/p)^{1/2} \{\|\M T\|_{\max} + \sigma\}$ and suppose that
   \begin{eqnarray}
     \label{eq:cond_matrix_comp1}
     np \succsim \log n, \quad \text{and} \quad
     |\lambda_{k_0}(\M T)|/E_n \succsim \kappa  (\log n)^{1/2},
   \end{eqnarray}
   where $\kappa = |\lambda_{1}(\M T)/\lambda_{k_0}(\M T)|$ is the condition number for $\M T$.
   Let $\hat{\M U}_g$ be obtained from Algorithm~\ref{RSRS}
  with
  $\sk = (1 - c_{\mathrm{gap}})^{-2} \{k_0 + \sqrt{12 k_0 \log n}
  + 6 \log n\}$ 
  and $g \geq g_*:= \frac{\log(n/\sk)}{\log(|\lambda_{k_0}(\M T)|/E_n)}.$
We then have, with probability at
  least $1 - 2n^{-3}$, that 
  \begin{gather}
\label{eq:matrix_completion}
d_{2}(\hat{\M U}_g,\M U)\precsim  \frac{(n/p)^{1/2}\{\|\M T\|_{\max} + \sigma\}}{|\lambda_{k_0}(\M T)|},
\\
\label{eq:matrix_completion2}
d_{\twoinf}(\hat{\M U}_g,\M U)\precsim \frac{\kappa^2 (n/p)^{1/2} (\log n)^{1/2} \{\|\M T\|_{\max} +
  \sigma\}\|\M U\|_{\twoinf}}{|\lambda_{k_0}(\M T)|}, 
\\
\label{eq:matrix_completion3}
\|\hat{\M T}_g - \M T\|_{\max} \precsim \kappa^4 (n/p)^{1/2} (\log n)^{1/2} \{\|\M T\|_{\max} + \sigma\}
\|\M U\|_{\twoinf}^2
\end{gather}
simultaneously. 
If $|\lambda_{k_0}(\M T)|/ E_n
\succsim n^{\epsilon}$ for a fixed but arbitary $\epsilon > 0$ then
Eq.~\eqref{eq:matrix_completion} through Eq.~\eqref{eq:matrix_completion3} holds for all $g \geq 1 + (2\epsilon)^{-1}$.
\end{theorem}
The uniform entrywise bound for $\hat{\M T}_g - \M T$ in Eq.~\eqref{eq:matrix_completion3}
can be further refined to yield entrywise limiting distributions. For ease of exposition we only consider the case where
$\M T$ is homogeneous, has finite rank, and bounded condition number, as these assumptions lead to
results that are simple to state while containing all key features
of more general results.
\begin{co}
  \label{co:entrywise}
  Consider the setting in Theorem~\ref{thm:mc} with $\M T$ satisfying following assumptions
\begin{enumerate}
\item[C1.] $\min_{k \ell} |T_{k \ell}| \asymp \|\M T\|_{\max}$ and
  $\mathrm{rk}(\M T) = k_0$ for some finite constant $k_0$. 
\item[C2.] $\kappa:= |\lambda_1(\M T)/\lambda_{k_0}(\M T)| \leq C_{\kappa}$ for some finite constant $C_{\kappa}$ not depending on $n$. 
\item[C3.] $np \succsim \log^{6}{n}$ and $p \leq 1 - \delta$ for some constant $\delta$ not depending on $n$.
\item[C4.] $\lambda_{k_0}(\M T) \succsim (n/p)^{1/2} \sigma \log^{3}{n}$. 
\end{enumerate}
  Let
$\zeta_{k \ell} = [\M U \M U^{\top}]_{k \ell}$  
  and denote the variance of $[\M U\M U^\T {\M E}]_{ij} + [{\M E}\M U\M U^\T ]_{ij}$ by
  \begin{equation}
    \label{eq:def_vstar}
    \begin{split}
v_{ij}^* &:= 
\frac{1}{p} \sum_{\ell \neq j} \bigl\{(1-p) T_{i \ell}^2 + \sigma^2\bigr\}
\zeta_{\ell j}^2  + \frac{1}{p} \sum_{\ell \neq i} \bigl\{(1-p)
T_{\ell j}^2 + \sigma^2\bigr\} \zeta_{i \ell}^2 
\end{split}
\end{equation}
Then for $\sk$ and $g \geq g_*$ as specified in Theorem~\ref{thm:mc}, 
and for any indices pair $(i,j)$, we have
\begin{equation}
  \label{eq:mc_clt1}
  (v^*_{ij})^{-1/2}[\hat{\M T}_g - \M T]_{ij} \rightsquigarrow
  \mathcal{N}(0,1) \quad \text{as $n \rightarrow \infty$}.
\end{equation}
\noindent\textbf{(Entrywise confidence interval)}
Let $\hat{\zeta}_{k \ell} = [\hat{\M U}_g \hat{\M U}_g^{\top}]_{k
  \ell}$ 
and define
\bee\label{def:empv}
\hat{v}_{ij} := \sum_{\ell \neq j}[\hat{\M E}]_{i \ell}^2
\hat{\zeta}_{\ell j}^2 + \sum_{\ell \neq i}[\hat{\M
  E}]_{\ell j}^2\hat{\zeta}_{i \ell}^2  + [\hat{\M
  E}]_{ij}^2\big\{\hat{\zeta}_{ii} + \hat{\zeta}_{jj} \big\}^2
\ee
where $\hat{\M E}:=\hat{\M T}_g - p^{-1}\hat{\M T}$. Then for any
indices pair $(i,j)$, we have 
\begin{equation}
  \label{eq:mc_clt2}
(\hat{v}_{ij})^{-1/2}[\hat{\M T}_g - \M T]_{ij}\rightsquigarrow
\mathcal{N}(0,1), \quad \text{as $n \rightarrow \infty$}.  
\end{equation}
If $|\lambda_{k_0}(\M T)|/E_n \succsim n^{\epsilon}$
for any fixed $\epsilon > 0$
then Eqs.~\eqref{eq:mc_clt1} and \eqref{eq:mc_clt2} holds for $g \geq 2 + (2 \epsilon)^{-1}$. 
\end{co}
Corollary~\ref{co:entrywise} provides more precise control of the entrywise fluctuations
for $\hat{\M T}_g - \M T$ compared to Theorem~\ref{thm:mc} and thus require slightly stronger conditions for
$np$ and $\lambda_{k_0}(\M T)$ compared to that in Eq.~\eqref{eq:cond_matrix_comp1}. These 
same conditions  for $np$ and $\lambda_{k_0}(\M T)$ as well as the assumption that $\M T$ is homogeneous 
were also used in the proof of Theorem~4.12 of \cite{chen2021spectral} for the estimator $p^{-1} \hat{\M T}^{(k)}$; homogeneity 
of $\M T$ guarantees 
that the entrywise noise levels for $\hat{\M T}$ are roughly on the same
order and leads to a convenient lower bound for $v_{ij}^*$. The assumption that $p$ is bounded away from $1$ is a mild assumption (it is used implicitly in the proof of Theorem~4.12 in \cite{chen2021spectral}) 
as the typical setting for matrix completion is that $p = o(1)$ as $n$ increases. The assumptions of finite rank and bounded condition numbers are also commonly seen in the literature. These assumptions can be
relaxed with substantially more involved book-keeping;
see Eqs.~(4.80) and (4.169) of \cite{chen2021spectral} for examples of conditions
where $\kappa$ and $k_0$ are allowed to vary with $n$. 
Finally, $\hat{v}_{ij}$ in Eq.~\eqref{def:empv}
is computable using only the RSVD output $\hat{\M U}_g$. 
 \subsection{Real data application: Distance matrix completion}\label{sec:RDA}
For this section we apply Algorithm \ref{RSRS} to recover the
missing entries of a partially observed Euclidean
distance matrix. In particular we use the 
\texttt{world\_cities} dataset containing the locations of the $4428$ most
populous cities around the world; this dataset is part of the
\texttt{mdsr} library in $\mathrm{R}$ \citep{baumer2017modern}. We
first construct the $4428 \times 4428$ matrix $\M D = [D_{ij}]$ whose
elements are \bee\nonumber D_{ij} = (\mathtt{Lon}_i - \mathtt{Lon}_j)^2 +
(\mathtt{Lat}_i - \mathtt{Lat}_j)^2. \ee
Here $\mathtt{Lon}_i$ and
$\mathtt{Lat}_i$ represent the longitude and latitude of the $i$th city,
respectively. 
We then sample a matrix $\M D_{0.8}$ (resp. $\M D_{0.4})$ by keeping roughly $80\%$ (resp. $40\%$) of the
entries in $\M D$, i.e.,
$\M D_{0.8} = \bm{\Omega}
\circ \M D$ where $\bm{\Omega}$ is a symmetric matrix whose upper
triangular entries are iid $\mathrm{Bernoulli}(0.8)$. 
 \begin{figure}[t]\centering
\includegraphics[width=0.23\linewidth]{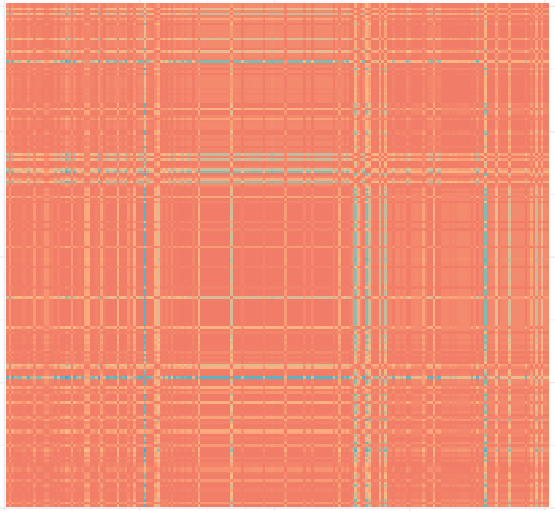} \quad\quad\quad
\includegraphics[width=0.23\linewidth]{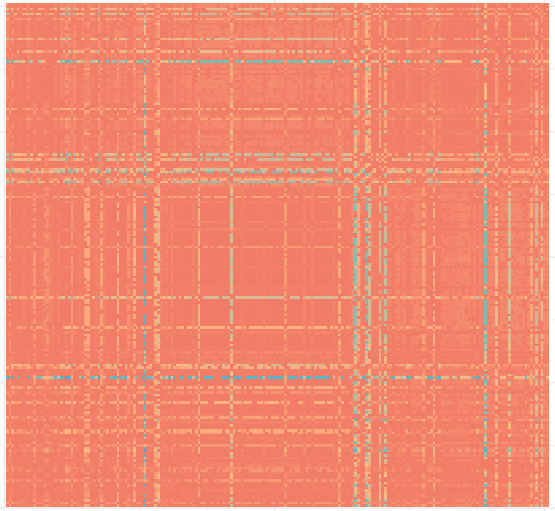} \quad\quad\quad
\includegraphics[width=0.273\linewidth]{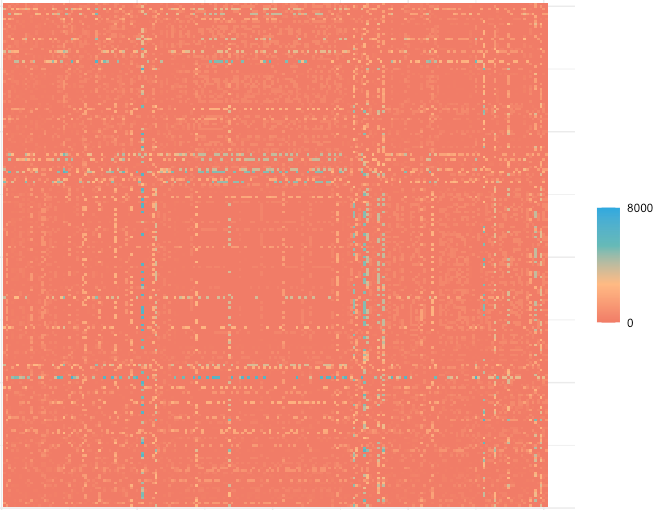} 
\caption{Matrix plots for the true $\mathbf{D}$ (Left), partially observed
  $\mathbf{D}_{0.8}$ (Middle), and partially observed $\mathbf{D}_{0.4}$ (Right).}\label{fig:MC1}
\end{figure}
 \begin{figure}[t]\centering
\includegraphics[width=0.2\linewidth]{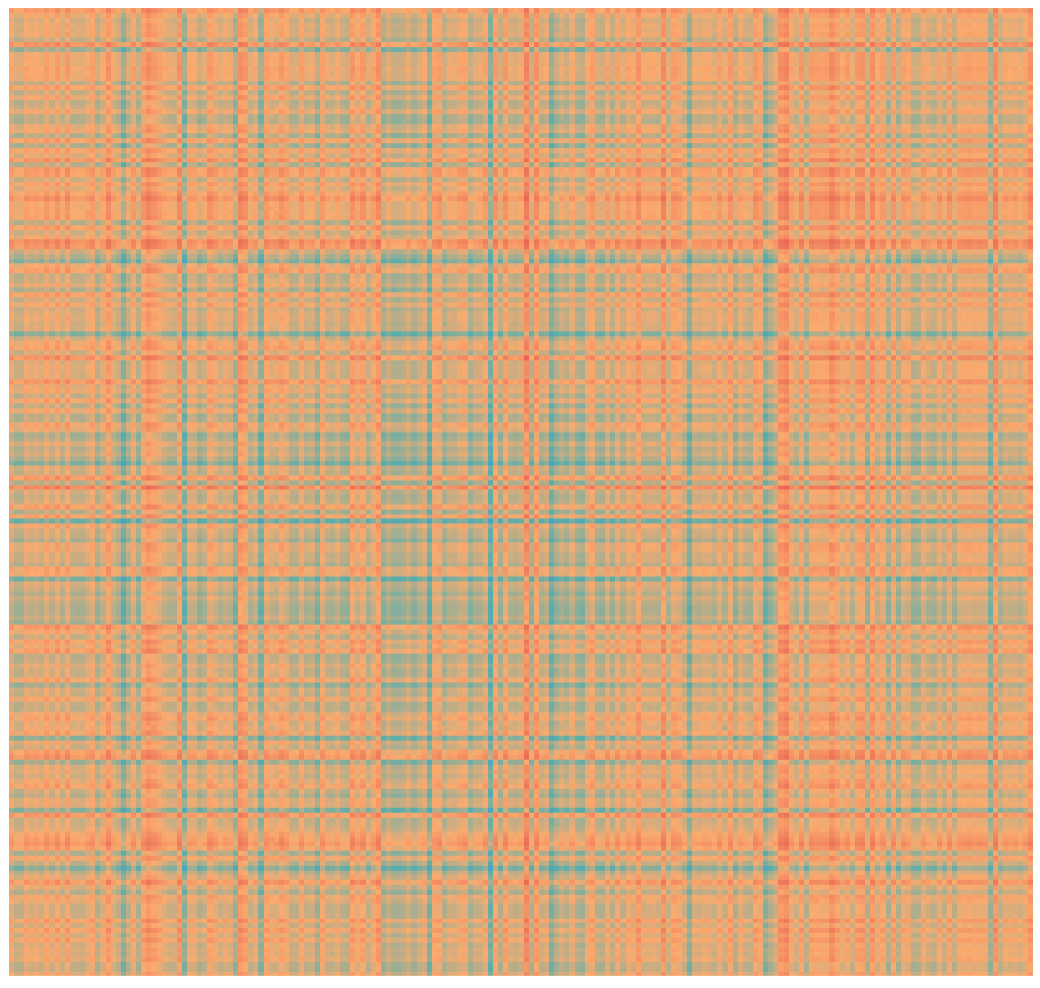} \quad
\includegraphics[width=0.2\linewidth]{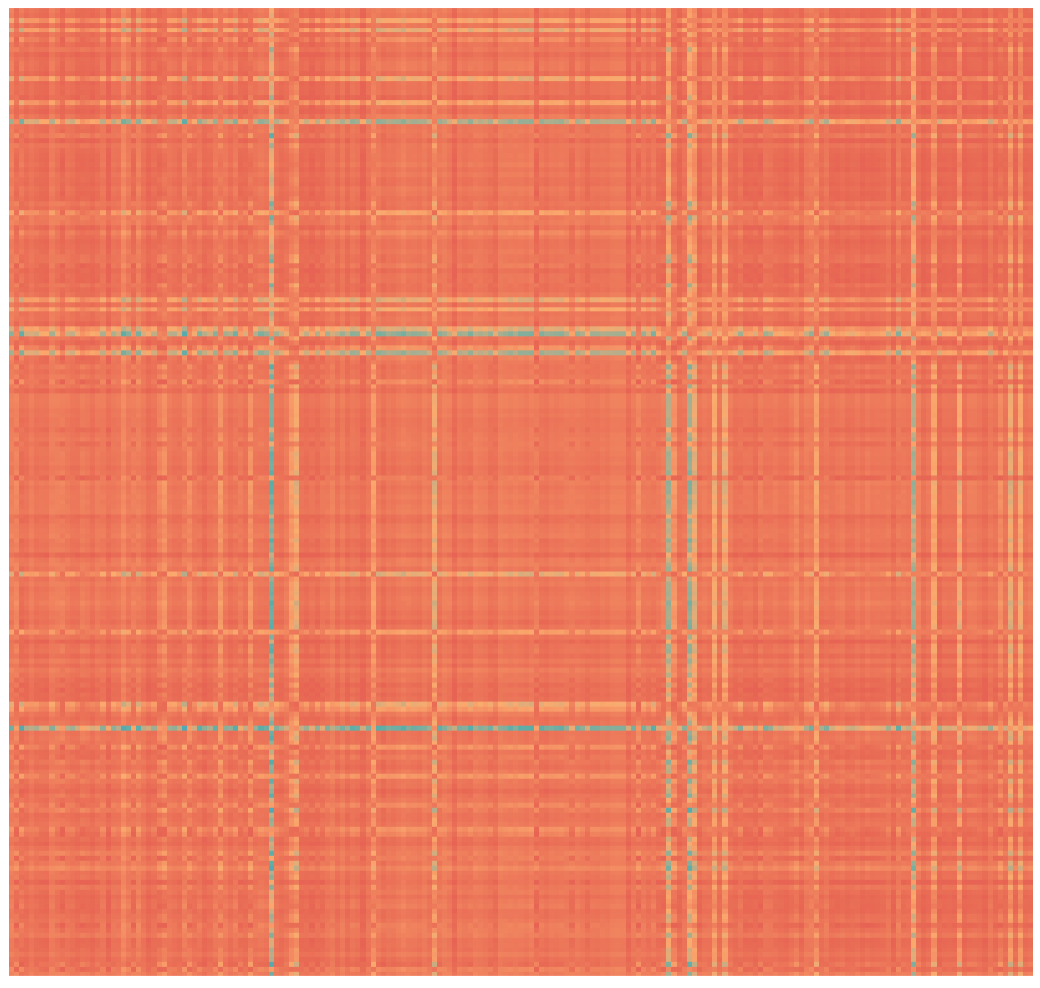} \quad
\includegraphics[width=0.2\linewidth]{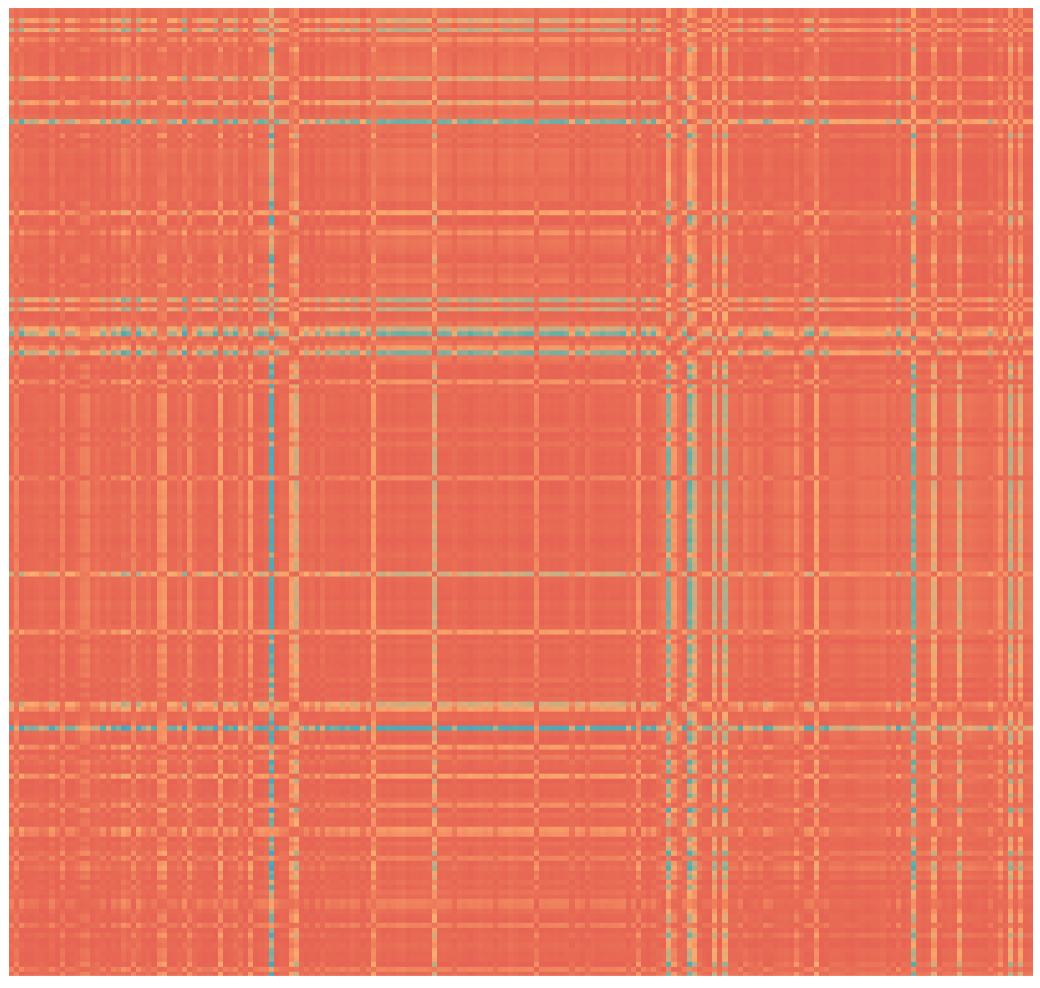} \quad
\includegraphics[width=0.244\linewidth]{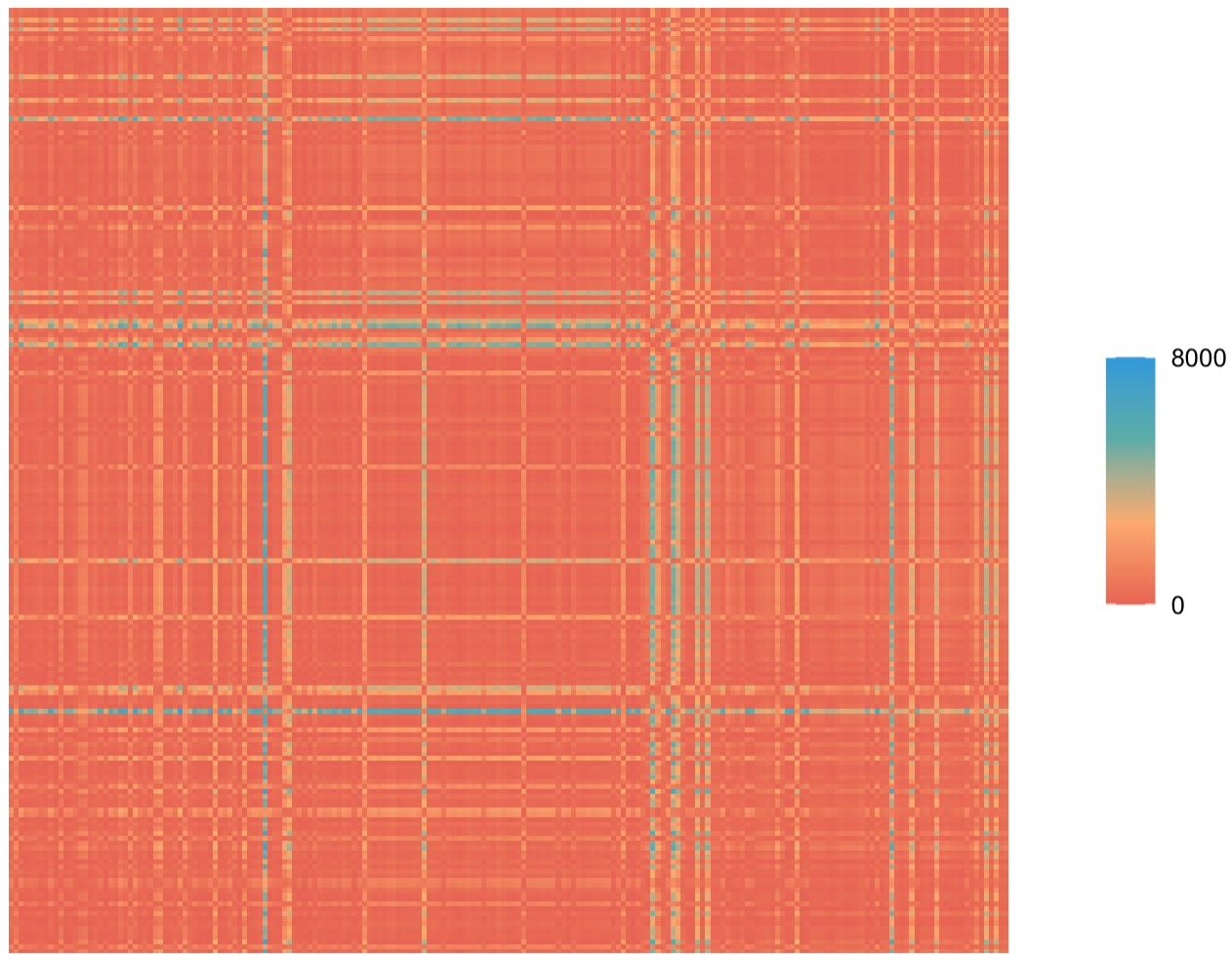} 
\caption{From left to right: matrix plots of RSVD-based estimates
  $\hat{\mathbf{D}}_{0.8}^{(1)}$, $\hat{\mathbf{D}}_{0.8}^{(2)}$, $\hat{\mathbf{D}}_{0.8}^{(5)}$, and exact SVD estimate $\hat{\mathbf{D}}_{0.8}$}\label{fig:MC2}
\end{figure}
 \begin{figure}[t]\centering
\includegraphics[width=0.2\linewidth]{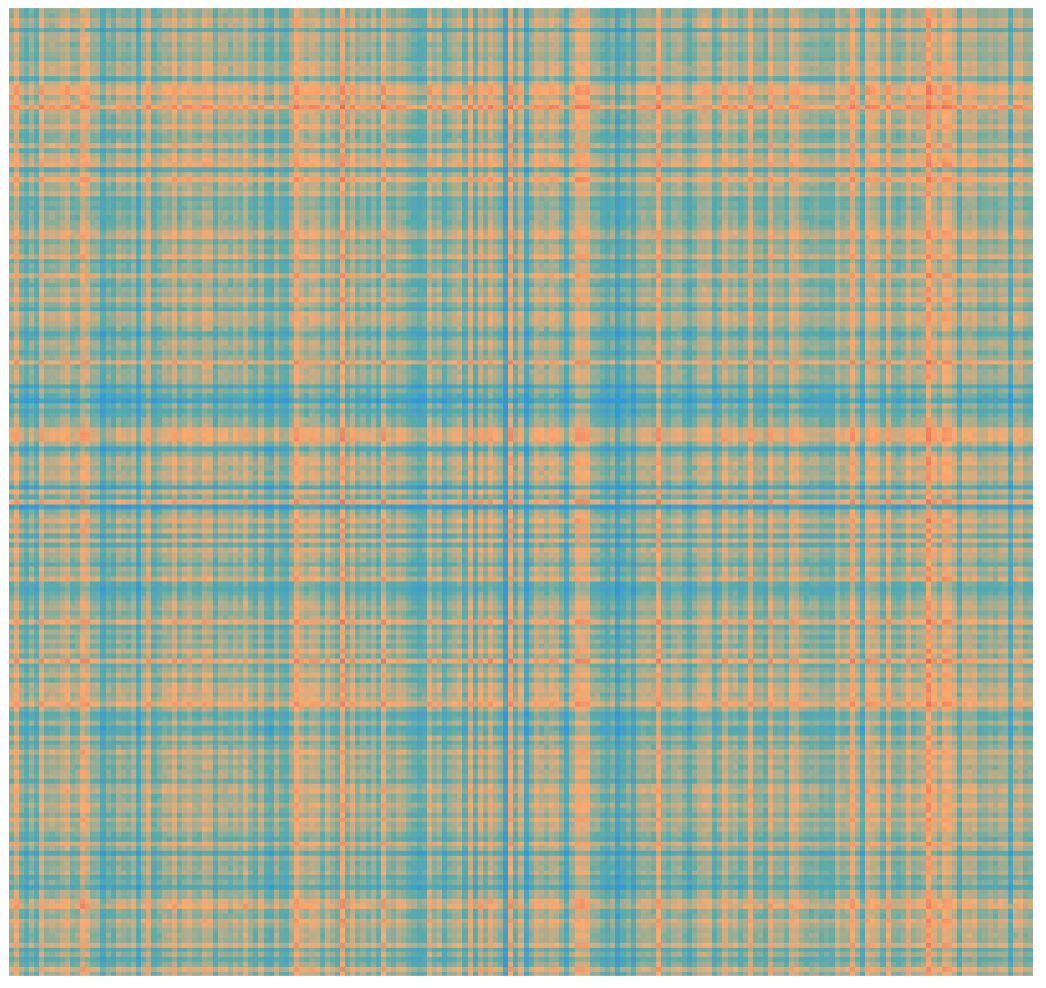}\quad 
\includegraphics[width=0.2\linewidth]{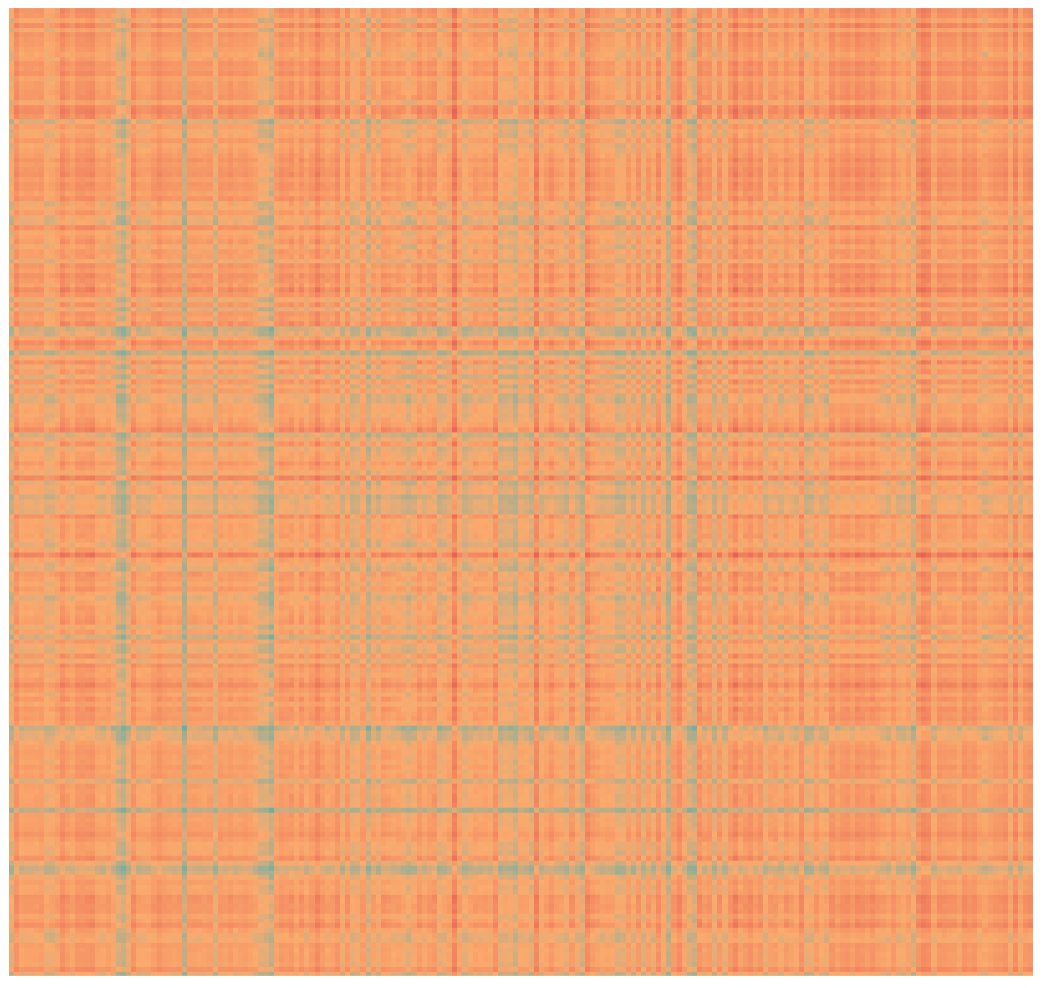} \quad
\includegraphics[width=0.2\linewidth]{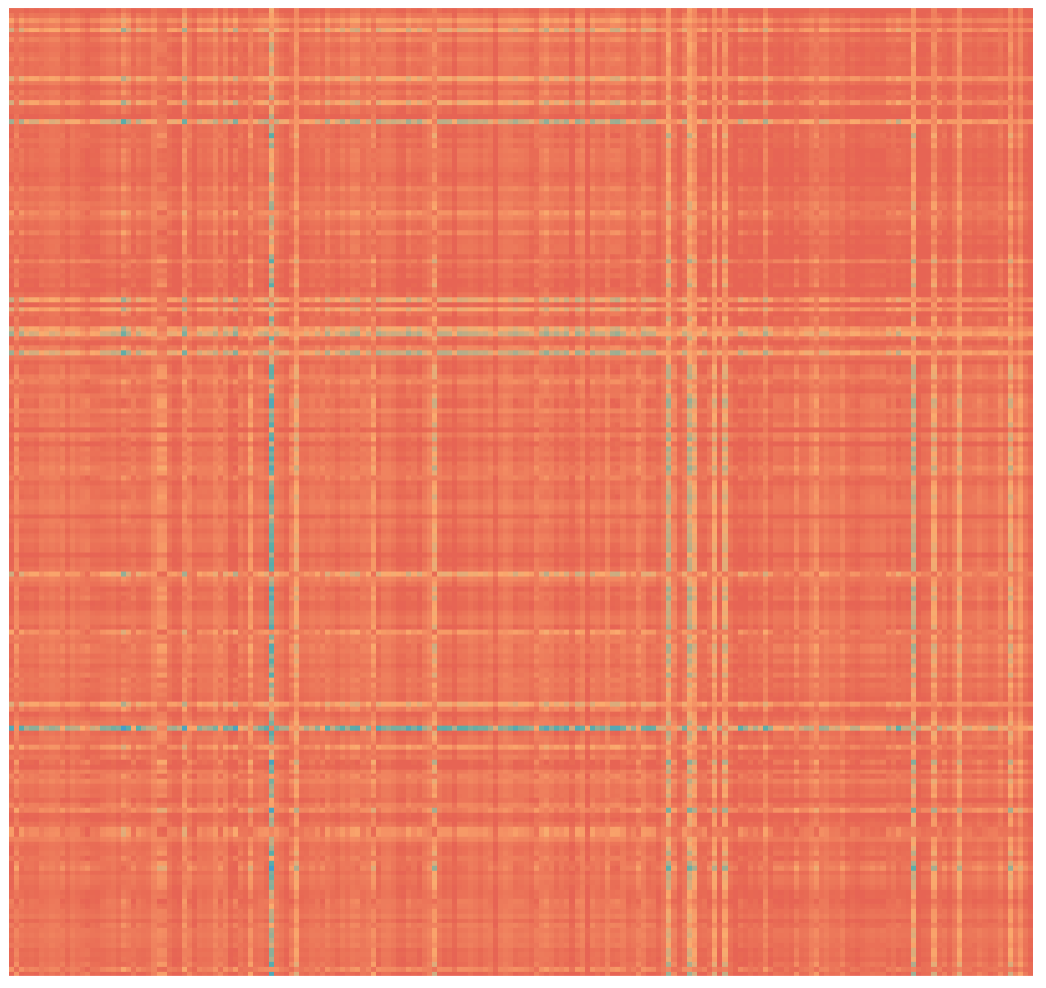} \quad
\includegraphics[width=0.244\linewidth]{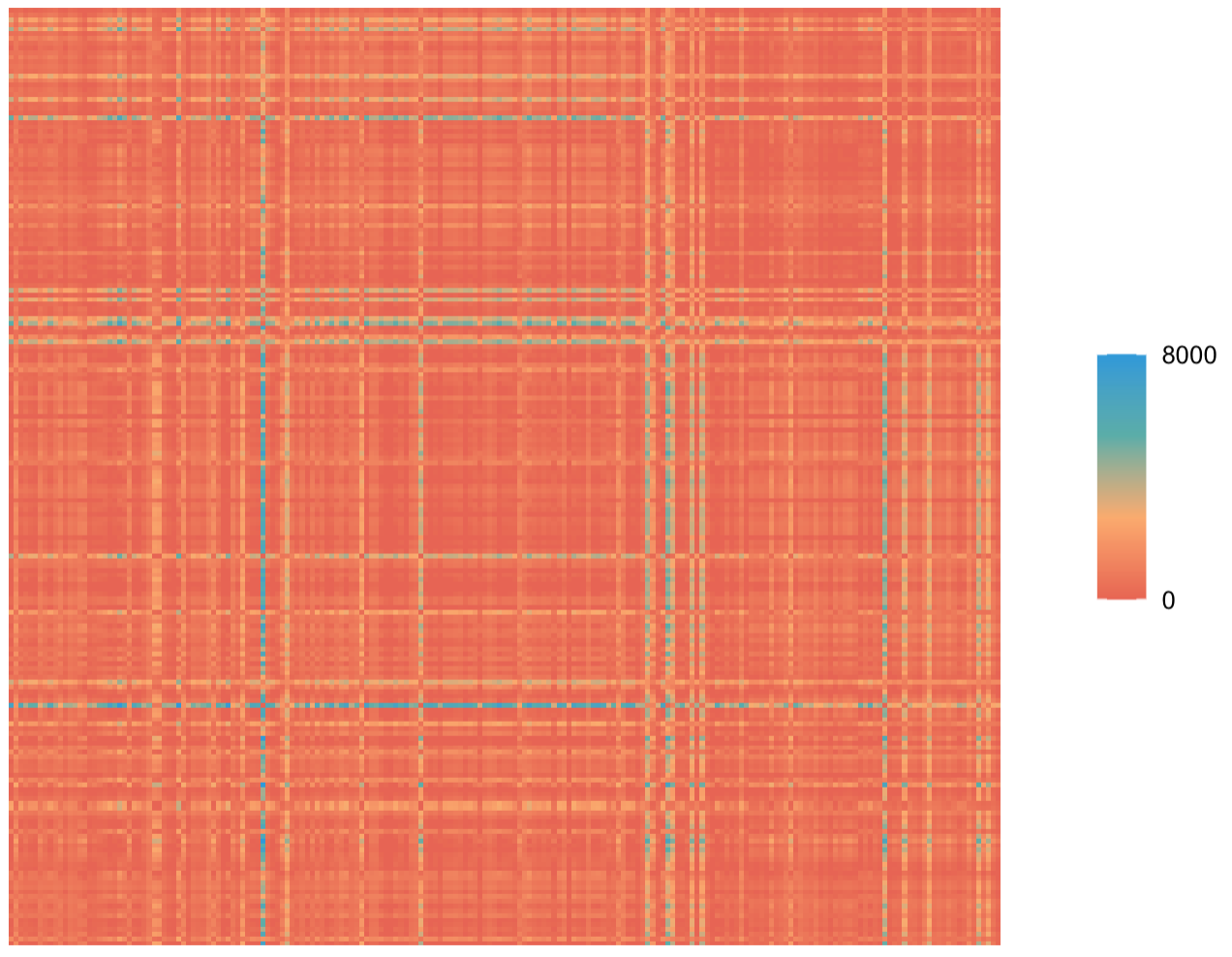} 
\caption{From left to right: matrix plots of RSVD-based estimates
  $\hat{\mathbf{D}}_{0.4}^{(1)}$, $\hat{\mathbf{D}}_{0.4}^{(2)}$, $\hat{\mathbf{D}}_{0.4}^{(5)}$, and truncated (exact) SVD estimate $\hat{\mathbf{D}}_{0.4}$}\label{fig:MC3}
\end{figure}
\begin{table}[!]
\setlength{\tabcolsep}{13pt}
\caption{The $\alpha$th--quantiles for the relative entrywise
  errors of the RSVD-based estimates $\hat{\mathbf{D}}_{p}^{(g)}$ and truncated (exact) SVD estimate
  $\hat{\mathbf{D}}_{p}$.
} \footnotesize
\centering{ 
\begin{tabular}{ c|ccccccc }
\hline
$\alpha$&$0.05$&$0.15$&
$0.35$&
$0.5$&
$0.65$ &$0.85$ &$0.95$ \\
\hline
 $\hat{\M D}_{0.8}^{(1)}$ & 0.0169& 0.0517& 0.1320& 0.2210& 0.4199 &2.4288& 15.3182 \\
  $\hat{\M D}_{0.8}^{(2)}$ & 0.0038 &0.0115 &0.0298& 0.0504 &0.0899 &0.3170 & 1.8230  \\
  $\hat{\M D}_{0.8}^{(5)}$ & 0.0024 &0.0073& 0.0182 &0.0289& 0.0465 &0.1532&  0.7829  \\
  $\hat{\M D}_{0.8}$& 0.0024& 0.0073& 0.0182& 0.0289& 0.0464& 0.1530& 0.7816\\
  \hline
 $\hat{\M D}_{0.4}^{(1)}$ & 0.0381& 0.1140& 0.2701 &0.4120& 0.6794 &3.9035 &24.6562  \\
  $\hat{\M D}_{0.4}^{(2)}$& 0.0106& 0.0326& 0.0858& 0.1488& 0.2750 &1.3274&  9.0400  \\
  $\hat{\M D}_{0.4}^{(5)}$ & 0.0073&  0.0221&  0.0558&  0.0911 & 0.1570&  0.5468 &  3.3533 \\
  $\hat{\M D}_{0.4}$ &0.0064& 0.0194& 0.0484& 0.0774& 0.1269& 0.4439& 2.4614\\
\hline
\end{tabular}}
\label{tb:3}
\end{table}
We now recover $\M D$ from $\M D_{0.8}$ (resp. $\M D_{0.4}$) using
Algorithm~\ref{RSRS}. As the entries of $\M D$ are  Euclidean
distances between points in $\mathbb{R}^2$, we have $\mathrm{rk}(\M D) \leq
4$. We therefore choose $k = 4$, $\sk = 20$ and $g \in \{1,2,5\}$, and let $\hat{\M D}_{0.8}^{(g)}$ 
(resp. $\hat{\M D}_{0.4}^{(g)}$) be the resulting estimate of $\M
D$. For comparison we also consider the spectral 
estimate $\hat{\M D}_{0.8}$ (resp. $\hat{\M D}_{0.4}$) obtained by
truncating the exact SVD of $\M D_{0.8}$ (resp. $\M D_{0.4}$); see
\cite{keshavan2009gradient} for more details. 

A plot of the true distance matrix $\M D$ and one random realization
of the partially observed
$\M D_{0.8}$ and $\M D_{0.4}$ are presented in 
Figure \ref{fig:MC1}.
The corresponding RSVD-based estimates and
exact SVD based estimates of $\M D$ are then shown in Figure
\ref{fig:MC2}--\ref{fig:MC3}. 
Figure~\ref{fig:MC2} shows that $\hat{\M
  D}_{0.8}^{(2)}$ and $\hat{\M D}_{0.8}^{(5)}$ both have comparable
accuracy to $\hat{\M D}_{0.8}$ while Figure~\ref{fig:MC3} shows that $\hat{\M D}_{0.4}^{(2)}$ has much
worse accuracy compared to $\hat{\M D}_{0.4}^{(5)}$ and $\hat{\M D}_{0.4}$. 

We also record the entrywise
relative errors between the RSVD-based estimates $\hat{\M
  D}_{p}^{(g)}$ (resp.~the exact SVD-based estimate $\hat{\M
  D}_{p}$) against that of $\M D$. A summary of the quantile levels for these
relative errors are presented in 
Tables \ref{tb:3}. The $\alpha$th--quantile of the entrywise relative errors
between an estimate $\M Z$ and the true distance $\M D$ is defined as the $\alpha$th quantile
of $\big\{|[\M Z - \M D]_{ij}/[\M D]_{ij}| \colon
(i,j)\in[4428]\times [4428]\big\}$; for example the median relative
error for $\hat{\M D}_{0.8}^{(5)}$ and $\hat{\M D}_{0.4}^{(5)}$ are
$\approx 0.029$ and $\approx 0.086$, respectively. Note that the numbers in
Table~\ref{tb:3} are averaged over $200$ independent random samples of
either $\M D_{0.8}$ or $\M D_{0.4}$. 
From Table~\ref{tb:3} we 
see that the relative error decreases as $g$ increases with $p$
fixed; indeed, the relative errors of $\hat{\M D}_{0.8}^{(5)}$ are nearly identical to those of $\hat{\M D}_{0.8}$. Table~\ref{tb:3}
also indicates that as $p$ decrease we need to increase $g$ 
to achieve a recovery rate close to that of $\hat{\M D}$. These observations
are consistent with the theoretical results in  
Theorem \ref{thm:mc}.

\section{PCA with missing data}\label{sec:epca:supp}
\subsection{Computational refinement and theoretical results}\label{sec:epca}
We now consider principal components estimation
with missing data. In particular, we focus on the following factor model from \cite{cai2021subspace}: 
\begin{equation}
  \label{eq:pca_setup}
\M X = \M B \M F + \M  N.
\end{equation}
Here $\M
B$ is a $d \times k_0$ matrix, 
$\M F = [\mathbf{f}_1,\dots,\mathbf{f}_m]\in\RR^{k\times m}$ is a $k_0
\times m$ matrix whose entries are iid $\mathcal{N}(0,1)$, and $\M N$ is a $d \times m$ matrix whose entries are iid $\mathcal{N}(0, \sigma^2)$; note that $\M B$, $\M N$ and $\M F$ are assumed to be mutually
independent. The columns of $\M X$ 
are then iid random vectors with mean $\bm{0}$ and covariance matrix
$\mathbf{B} \mathbf{B}^{\top} + \sigma^2
\mathbf{I}_d$. The columns of $\M B\M F$ are   iid signal random vectors with mean $\M 0$ and a low-rank covariance matrix $\M B\M B^\T$. Denote eigendecomposition $\M B \M
B^{\top} = \M U \bm{\Lambda} {\M U}^{\T}$ where $\mathbf{U} \in
\mathbb{O}_{d \times k}$ and $\M \Lambda =
\diag(\lambda_1,\dots,\lambda_k)$;
the columns of $\M U$ are the leading principal components.

Due to sampling issues and/or privacy-preserving intention,
it is often the case that only a partial subset of the entries in $\M
X$ are observed. More specifically let $\M \Omega$  be a $d \times m$ binary
matrix whose entries are iid Bernoulli random variables with success
probability $p$. Then, instead of
observing $\M X$, we only observe 
$$
\M Y = \mathcal{P}_{\M\Omega}(\M X) =
\bm{\Omega} \circ \mathbf{X} = \bm{\Omega} \circ (\M B \M F + \M  N),
$$ 
where $\circ
$ denotes the Hadamard
product between matrices.
Given observed 
$\M Y$,
\cite{cai2021subspace} propose the following spectral procedure for 
recovering the principle components $\M U$: form the matrix
\begin{equation}
  \label{eq:J_pca}
\M Q = \frac{1}{mp^2} \mathcal{P}_{\text{off-diag}}\left( \M Y\M Y^\T\right) ,
\end{equation}
where $\mathcal{P}_{\text{off-diag}}(\cdot)$ sets
 diagonal entries of the corresponding matrix to zero, and compute the
$d \times k_0$ matrix $\hat{\M U}$ whose columns
are the leading singular vectors of $\M Q$.
As the dimension $d$ can be reasonably
large compared to the number of samples $m$ while the number of non-zero entries
in $\M Y$ can be much smaller than $md$,  
we can replace the singular vectors
$\hat{\M U}$ of $\M Q$ by the approximate singular vectors $\hat{\M U}_g$
computed using RSVD.

\begin{algorithm}[tp]
\KwIn{$\M Y\in\RR^{d\times m}$, rank $k\geq 1$, sketching
  dimension $\sk \geq k$, power iterations $g \geq 1$.}
 Generate a $d \times k$ sketching matrix $\M G$ whose elements are iid standard normals\;
  Compute the diagonal entries of $\M Y\M Y^\T$ and obtain $\mathcal{P}_{\mathrm{diag}}(\M Y\M Y^\T)$\;   

  \For{$s = 1$ \KwTo $g$}{
Compute:
\bee\nonumber
\M Y^\T  \M Q^{s-1}\M G& \gets \M Y^\T ( \M Q^{s-1}\M G);
\\ 
\M Q^{s}\M G&\gets\frac{1}{mp^2}\M Y \left(\M Y^\T  \M Q^{s-1} \M G\right) - \frac{1}{mp^2}\mathcal{P}_{\text{diag}}\left(\M Y\M Y^\T\right) \M Q^{s-1} \M G\text{\;}
\ee}

Obtain the {\em exact} SVD of $\M Q^g\M G$. Let $\hat{\M U}_g$ be the $d\times k $ matrix  whose columns are the $k$ leading left singular vectors of $\M Q^g\M G$\;
 Compute the RSVD-based rank-$k$ approximation of ${\M Q}$. First,  compute 
$
 \hat{\M U}_g^\T   \M Y \M Y^\T  
$
 and $\hat{\M U}_g^\T  \mathcal{P}_{\text{diag}}(\M Y\M Y^\T) $ where the matrix multiplications are done from left to right. Then obtain
$$
\hat{\M Q}_g \gets  \frac{1}{mp^2}\hat{\M U}_g\times(\hat{\M U}_g^\T   \M Y \M Y^\T ) -  \frac{1}{mp^2}\hat{\M U}_g\times\left(\hat{\M U}_g^\T  \mathcal{P}_{\text{diag}}\left(\M Y\M Y^\T\right) \right);
$$

\KwOut{ 
  Subspace estimate $\hat{\M U}_{g} $ and  covariance matrix estimate $\hat{\M Q}_g$.}
 \caption{RSVD-based PCA with missing data}\label{RSRS:pca}
\end{algorithm}
Note that direct computation of \(\M Q\) as the input to Algorithm~\ref{RSRS} may involve
large-scale matrix multiplications when \(d\) and \(m\) are large. In particular,
we need to compute \( \M Y\M Y^\T\), which incurs a cost of \(O(md^2p)\) floating point operations (flops) and also require
substantial memory storage.  We address these computational challenges by refining
the standard RSVD-based approach for missing-data PCA
to avoid the need for explicit computation of \(\M Y\M
Y^\T\) as follows. First, for any \(\M
\Gamma \in \mathbb{R}^{d\times \sk}\),
\bee\label{sketchingdecom:basic}
\M Q  \M \Gamma 
&= \frac{1}{mp^2} \mathcal{P}_{\text{off-diag}}\left( \M Y \M Y^\T  \right) \M \Gamma 
\\
&= \frac{1}{mp^2}\left\{ \M Y \M Y^\T - \mathcal{P}_{\text{diag}}\left( \M Y \M Y^\T \right) \right\}   \M \Gamma
= \frac{1}{mp^2}\M Y \bigl( \M Y ^\T \M \Gamma \bigr) - \frac{1}{mp^2}\mathcal{P}_{\text{diag}}\left( \M Y\M Y^\T\right) \M \Gamma,
\ee
where \(\mathcal{P}_{\mathrm{diag}}(\cdot)\) denotes the operator that preserves only the diagonal entries.
Hence, instead of explicitly forming \(\M Q\) and then computing the sketched matrix \(\M Q\M G\) in Algorithm~\ref{RSRS},
we can set \(\M \Gamma = \M G\) and compute
\begin{equation}
\label{sketchingdecom}
\frac{1}{mp^2}\M Y \left(\M Y^\T \M G \right) - \frac{1}{mp^2}\mathcal{P}_{\text{diag}}\left( \M Y\M Y^\T\right) \M G,
\end{equation}
which only requires access to \(\M Y\) and the diagonal entries of $\M Y\M Y^\T$.
In addition, the first term in the right-hand side of Eq.~\eqref{sketchingdecom} involves matrix multiplication between smaller matrices and has a computational cost of
\(O(md\sk p)\) flops while the second term requires computing the diagonal entries of \(\M Y \M Y^\T\) (which require \(O(mdp)\) flops)
followed by a matrix multiplication between $\mathcal{P}_{\text{diag}}\left( \M Y\M Y^\T \right)$ and $\M G$ (which costs \(O(d\sk)\) flops). 
Overall, Eq.~\eqref{sketchingdecom} reduces the cost for computing $\M Q \M G$ from
\(O(md^2p)\) to $O(md\sk p + d\sk)$ flops (which is significant reduction when $\sk \ll d$). The same approach can also be used 
for computing 
\(\M Q(\M Q\M G), \M Q(\M Q^2\M G), \ldots,\M Q(\M Q^{g-1}\M G)\) by simply replacing $\M \Gamma = \M G$ with \(\M \Gamma = \M Q\M G, \M Q^2\M G,\dots,\M Q^{g-1}\M G\)
in Eq.~\eqref{sketchingdecom:basic}. See Algorithm~\ref{RSRS:pca} for more details. 

The following result combines Corollary~\ref{co:l2inf_noise} with error bounds given in Corollary~4.3 of
\cite{cai2021subspace} to show that the $\hat{\M
  U}_g$ achieves the same estimation rate as that for $\hat{\M U}$ in approximating the principal subspace; note that the condition for $m$ in Eq.~\eqref{eq:sample_size} below
is identical to Eq.~(4.14) of \cite{cai2021subspace}. 
 
 \begin{theorem}\label{thm5}
   Let $\mathbf{X}$ be a $d \times m$ matrix sampled according
   to Eq.~\eqref{eq:pca_setup}. 
Let $\mu=dr^{-1}\|\M U\|_{\twoinf}^2$ denote the coherence parameter for $\M
U$, $\kappa = \lambda_1/\lambda_{k_0}$ the condition number for $\M B \M B^{\top}$, and
$s_{\ast} = \log(m+d)$. Suppose there exist constants $\tilde{c}_0 > 0$ and $\tilde{c}_1 > 0$ such that $k_0
\leq \frac{\tilde{c}_1 d}{\mu \kappa^2}$ and $m$ satisfies the sample size condition
\begin{equation}
  \label{eq:sample_size}
m \geq
\tilde{c}_0\max\Bigl\{\frac{\mu^2 \kappa^6 k_0^2s_{\ast}^6}{dp^2},\frac{\mu \kappa^5 k_0 s_{\ast}^3}{p}, \frac{\sigma^4 \kappa^2 s_{\ast}^2}
{\lambda_{k_0}^2p^2},
\frac{\sigma^2 \kappa^3 ds_{\ast}}{\lambda_{k_0} p} \Bigr\}.
\end{equation}
Now define
\begin{equation}
  \label{eq:def_E_pca}
\begin{split}
\mathscr{E} &:= \frac{\mu \kappa^2 k_0 s_{\ast}^2}{(md)^{1/2}p} +
\frac{(\mu \kappa^3 k_0)^{1/2} s_{\ast}}{(md)^{1/2}}+
\frac{\sigma^2d^{1/2}s_{\ast}}{\lambda_{k_0}^2 m^{1/2}p}+\frac{\sigma (\kappa d s_{\ast})^{1/2}}{(\lambda_{k_0} m p)^{1/2}}
+ \frac{\mu \kappa k_0}{d}.
\end{split}
\end{equation}
Let $\hat{\M U}_g$ be generated via Algorithm~\ref{RSRS:pca} (or equivalently, generated from Algorithm~\ref{RSRS} with $\hat{\M M} = \M Q$), 
$\sk \geq (1 - c_{\mathrm{gap}})^{-2} \{k_0 + \sqrt{24 k_0 \log d} + 6 \log d\}$ and
$g \geq \frac{\log d}{\log (1/\mathscr{E})}$.
Then with probability at least $1 - m^{-3}$, we have 
\begin{equation}
  \label{eq:PCA_missing}
  d_2(\hat{\M U}_g,\M U) \precsim \mathscr{E}, \quad \text{and} \quad
  d_{\twoinf}(\hat{\M U}_g, \M U) \precsim \kappa^{3/2} \log^{1/2}(m+d) \mathscr{E} \|\M U\|_{\twoinf}.
\end{equation}
Finally, if $\mathscr{E} \precsim d^{-\epsilon}$ for a fixed $\epsilon > 0$ then 
Eq.~\eqref{eq:PCA_missing} holds for all $g \geq 1 + (2 \epsilon)^{-1}$. 
\end{theorem}
\subsection{Extension: HeteroPCA with RSVD approximation}
Recent work has studied heteroskedastic PCA (HeteroPCA)~\citep{zhang2022heteroskedastic,agterberg2022entrywise,yan2024inference,zhou2025deflated},  a more elaborate PCA algorithm compared with the diagonal-deleted PCA in Section~\ref{sec:epca}. It is designed to address bias arising from heteroskedastic noise in \(\M N\) and to further enable efficient inference for the underlying principal subspace. 
Rather than directly using the leading singular vectors of the diagonal-deleted matrix \(\M Q\) to estimate \(\M U\), HeteroPCA treats the top-\(k\) eigendecomposition of \(\M Q\) as an initial spectral estimator of the underlying covariance matrix \(\M B\M B^\T = \M U\M\Lambda\M U^\T\). It then iteratively imputes the diagonal entries of \((mp^2)^{-1}\mathcal{P}_{\text{off-diag}}(\M Y\M Y^\T)\) with the diagonal entries of the current spectral estimate of $\M B\M B^\T$, followed by updating the spectral estimate of $\M B\M B^\T$ via the top-\(k\) eigendecomposition of the newly imputed matrix; see Algorithm~1 in~\citet{zhang2022heteroskedastic} for details.
As each iteration of HeteroPCA requires an SVD of a \(d \times d\) matrix, we can replace the SVD step with RSVD to improve computational efficiency. See
Algorithm~\ref{RSRS:pca:hetero} for more details. Note that Algorithm~\ref{RSRS:pca:hetero} also avoids direct computation of the matrix \(\M Y\M Y^\T\).
\begin{algorithm}[tp]
\KwIn{$\M Y\in\RR^{d\times m}$, rank $k\geq 1$, sketching
  dimension $\sk \geq k$, power iterations $g \geq 1$, number of HeteroPCA iterations  $T \geq 1$.}
Run Algorithm~\ref{RSRS:pca}\;
Obtain $\hat{\M Q}_g^{(0)}\gets \hat{\M Q}_g$ and   $\mathcal{P}_{\mathrm{diag}}(\M Y\M Y^\T)$ from Algorithm~\ref{RSRS:pca}\;
\For{$t = 1$ \KwTo $T$}{
 Generate a $d \times k$ sketching matrix $\M G^{(t)}$ whose elements are iid $\mathcal{N}(0,1)$\;
\For{$s = 1$ \KwTo $g$}{
Compute:
\bee\nonumber
\M Y^\T  \M Q_t^{s-1}\M G^{(t)}& \gets \M Y^\T \times \left( \M Q_t^{s-1}\M G^{(t)}\right);
\\ 
\M Q_t^{s}\M G^{(t)}&\gets\frac{1}{mp^2}\M Y \times \left(\M Y^\T  \M Q_t^{s-1} \M G^{(t)}\right) - \frac{1}{mp^2}\mathcal{P}_{\text{diag}}\left(\M Y\M Y^\T\right)\times  \M Q_t^{s-1} \M G^{(t)}
\\
&\quad\,\, +  \mathcal{P}_{\mathrm{diag}}\left(\hat{\M Q}_g^{(t-1)}\right)\times  \M Q_t^{s-1} \M G^{(t)};
\ee
}

Obtain the {\em exact} SVD of $\M Q_t^g\M G^{(t)}$. Let $\hat{\M U}^{(t)}_g$ be the $d\times k $ matrix  whose columns are the $k$ leading left singular vectors of $\M Q_t^g\M G^{(t)}$\;
Compute the RSVD-based rank-$k$ approximation of ${\M Q}_t$. First,  compute 
$
 \big(\hat{\M U}^{(t)}_g\big)^\T   \M Y \M Y^\T  
$, $\big(\hat{\M U}^{(t)}_g\big)^\T  \mathcal{P}_{\text{diag}}(\M Y\M Y^\T)$ and $\big(\hat{\M U}^{(t)}_g\big)^\T\mathcal{P}_{\text{diag}}(\hat{\M Q}_g^{(t-1)})$ where the matrix multiplications are proceeded from left to right.  Then compute
\bee\nonumber
\hat{\M Q}_g^{(t)} \gets & \frac{1}{mp^2}\hat{\M U}^{(t)}_g\times\left(\big(\hat{\M U}^{(t)}_g\big)^\T   \M Y \M Y^\T\right) -  \frac{1}{mp^2}\hat{\M U}^{(t)}_g\times\left(\big(\hat{\M U}^{(t)}_g\big)^\T  \mathcal{P}_{\text{diag}} \left(\M Y\M Y^\T\right) \right)
\\
& +  \hat{\M U}^{(t)}_g\times \left(\big(\hat{\M U}^{(t)}_g\big)^\T\mathcal{P}_{\text{diag}}(\hat{\M Q}_g^{(t-1)})\right);
\ee}

\KwOut{Subspace estimate $\hat{\M U}^{(T)}_{g} $ and  covariance matrix estimate $\hat{\M Q}_g^{(T)}$.}
 \caption{RSVD-based HeteroPCA}\label{RSRS:pca:hetero}
\end{algorithm}
\begin{remark}The matrix \(\M Q_t\) in Algorithm~\ref{RSRS:pca:hetero} denotes the RSVD-based approximation of the
  diagonal-imputed sample covariance matrix at the \(t\)-th iteration of HeteroPCA, i.e.,
\bee\nonumber
\M Q_t &=     \frac{1}{mp^2} \mathcal{P}_{\text{off-diag}} \bigl(\M Y\M Y^\T\bigr)
  +  \mathcal{P}_{\text{diag}}\bigl(\hat{\M Q}_g^{(t-1)}\bigr)
= \frac{1}{mp^2} \M Y \M Y^\T  -  \frac{1}{mp^2} \mathcal{P}_{\text{diag}} \bigl(\M Y\M Y^\T\bigr)  
  +  \mathcal{P}_{\text{diag}}\bigl(\hat{\M Q}_g^{(t-1)}\bigr), 
\ee
where $\hat{\M Q}_g^{(t-1)}$ is the RSVD-based rank-$k$ approximation of $\M Q_{t-1}$.
\end{remark}

{\color{black}Similar to Theorem~\ref{thm5},
we can apply Corollary~\ref{co:l2_noise} and Corollary~\ref{co:l2inf_noise} to
show that, if \(\sk\) and \(g\) are both sufficient large in each step of
Algorithm~\ref{RSRS:pca:hetero} then the RSVD-based subspace estimator
\(\hat{\M U}_g^{(T)}\) will achieves the same \(\ell_2\) and
\(\ell_{\twoinf}\) error rates as HeteroPCA (with exact SVD)
for recovering the true principal subspace $\M U$; see \citet{zhang2022heteroskedastic,yan2024inference} for theoretical results on HeteroPCA with exact SVD. 
Correspondingly, we also expect that one can use RSVD-based HeteroPCA 
to construct confidence regions for the principal subspaces and entrywise confidence intervals
for the covariance matrix (as done in \cite{yan2024inference} for HeteroPCA with exact SVD).
Finally we note that there are two iterations loops in 
Algorithm~\ref{RSRS:pca:hetero}: the outer loop 
for diagonal imputation and the inner loop for RSVD. An interesting
open problem is how to choose the number of iterations (\(g\) and \(T\)) in each loop to  
achieve the optimal trade-off between estimation accuracy and 
computational cost. Addressing this requires delicate analysis
of how \(g\) and \(T\) jointly affect the error rates of the RSVD-based HeteroPCA,
and so we leave it for future work.}

\subsection{Numerical experiments}\label{epcasimu}
\begin{figure}[t]
\includegraphics[width=0.8\linewidth]{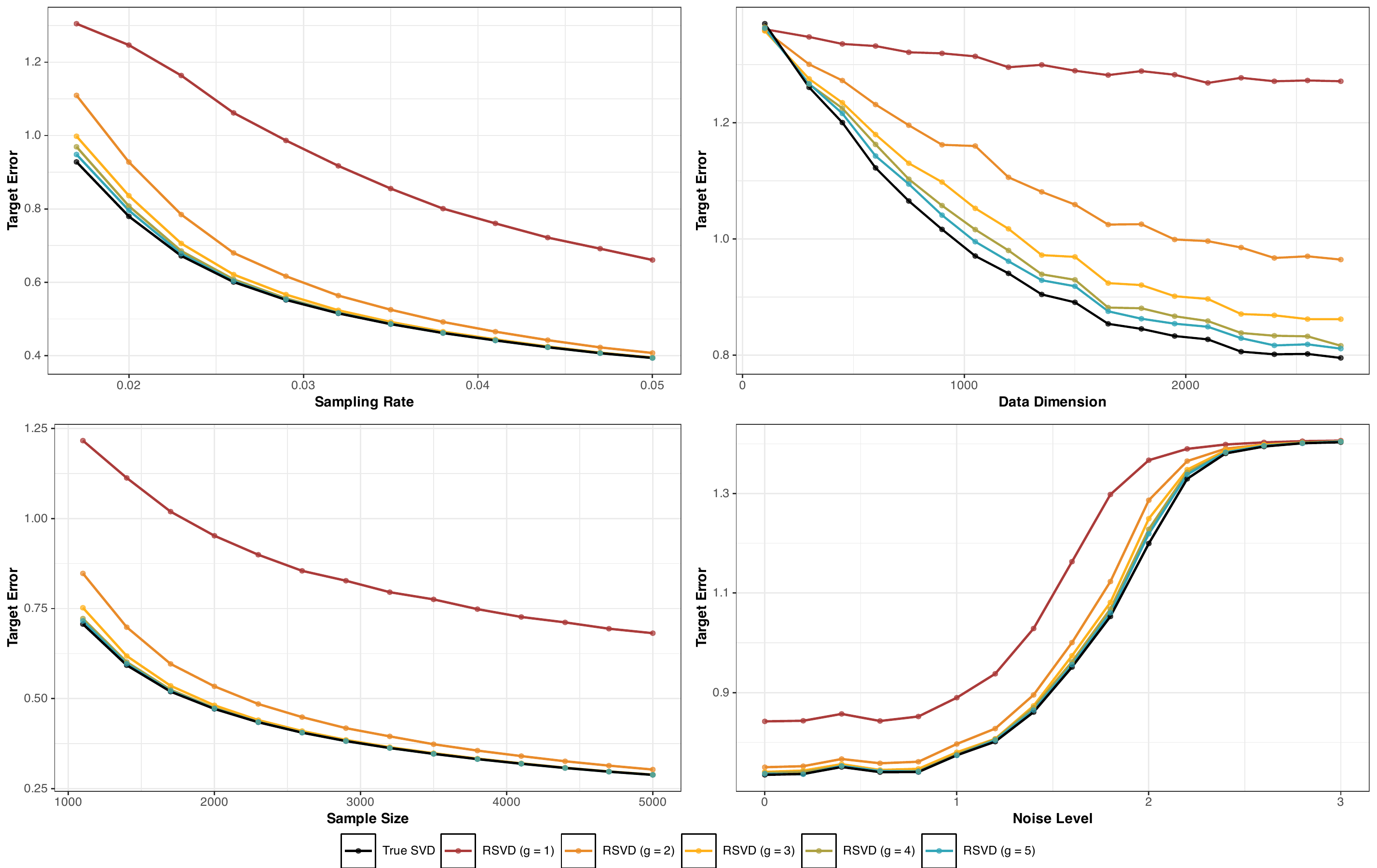} 
\centering
\caption{(Top Left) $\ell_2$ recovery errors vs. $p$ (here $d = 3000, n= 1000, \sigma =1$); (Top Right) $\ell_2$ recovery errors vs. $d$ (here $n= 1000, p =0.02, \sigma =1$); (Bottom Left) $\ell_2$ recovery errors vs. $n$ (here $ d = 3000, p =0.02, \sigma =1$); (Bottom Right) $\ell_2$ recovery errors vs. $\sigma$ (here $ d = 3000, n = 1000, p =0.02$).}
\label{fig:simu:cov}
\end{figure}
\begin{figure}[t] 
 \includegraphics[width=1\linewidth]{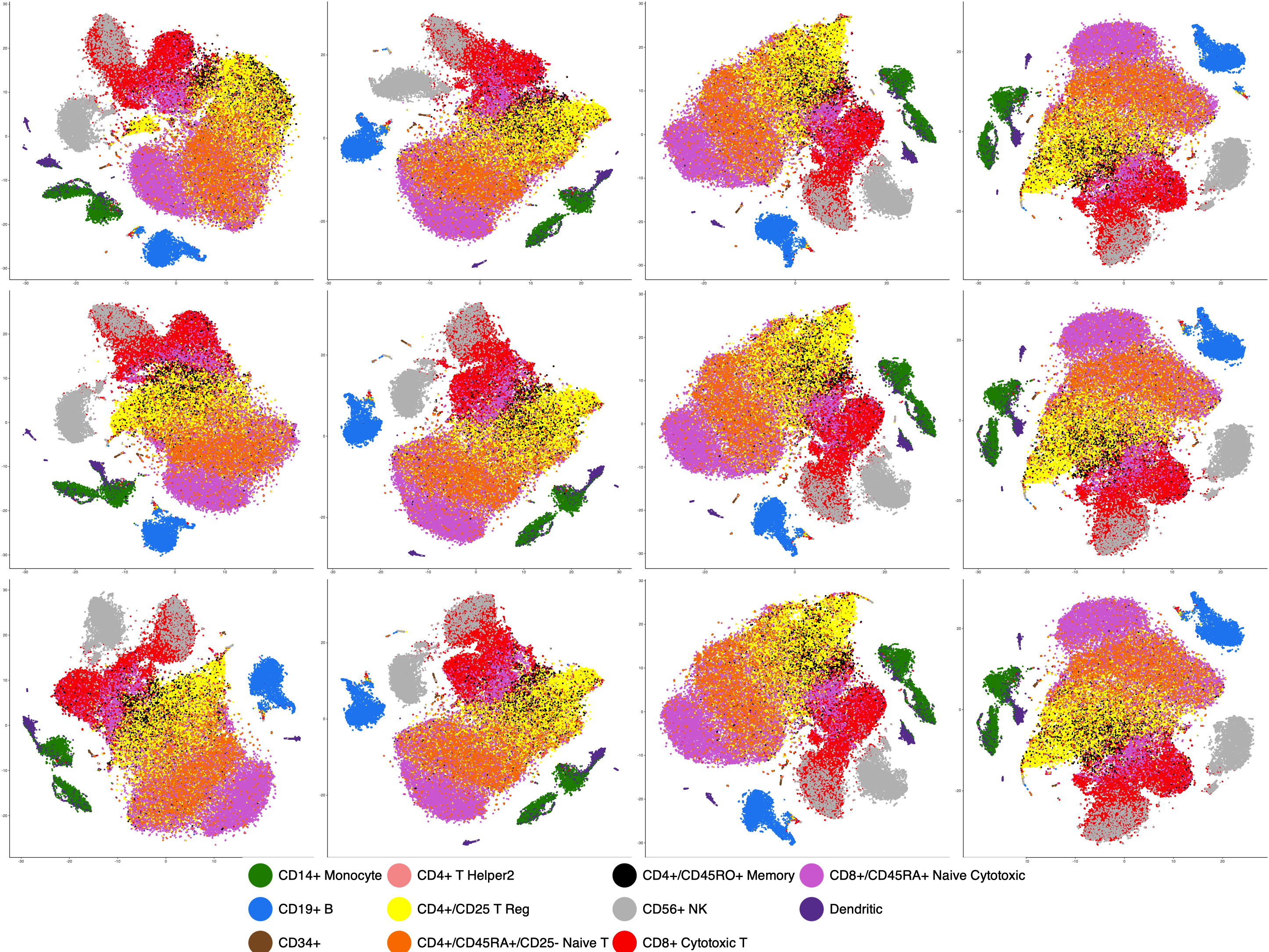}   
\caption{$t$-SNE embeddings of the 68k PBMC gene expressions,  projected  on the top 50 PCs obtained via RSVD-based PCA (Algorithm~\ref{RSRS:pca}). Top row (left to right): results for $g = 3$ with $\sk = 55$, $100$, $300$, and $1000$. Middle and bottom rows (left to right): results for $g = 4$ and $g = 5$, respectively, using the same sequence of $\sk$ values.}\label{fig:pca:real:supp}
\end{figure}
For this simulation study we used the same data generation mechanism as that
described in section $7$ of \cite{cai2021subspace} but with 
larger values of $d$. More specifically we first sample a $d \times n$
matrix $\M X^*$ whose columns are iid multivariate normal random vectors with
mean $\bm{0}$ and covariance matrix $\mathbf{U}^* \mathbf{U}^{*\top}$;
here $\mathbf{U}^*$ is a $d \times 4$ matrix whose entries are iid
standard normals. We then generate
$\mathcal{P}_{\bm{\Omega}}(\M X) = \bm{\Omega} \circ( \M X^* +
\M E)$ where the entries of $\M E$ are iid $\mathcal{N}(0,
\sigma^2)$, the entries of $\bm{\Omega}$ are iid
$\mathrm{Bernoulli}(p)$ and $\circ$ denote the Hadamard product. The
matrix $\mathcal{P}_{\bm{\Omega}}(\M X)$ represent a noisily observed
version of $\M X^*$ (with missing entries). The matrices $\M U^*, \M
  E$ and $\bm{\Omega}$ are resampled for each Monte Carlo replicate. 

Given $\mathcal{P}_{\bm{\Omega}}(\mathbf{X})$ ,we compute
$\hat{\M U}_g$ using Algorithm~\ref{RSRS:pca} with $k = 4$, $\sk =
45$ 
and $g \in \{1,2,\dots,5\}$. We then
record $d_2(\hat{\M U}_g, \M U)$; here $\M U$ denotes the $d
\times 4$ matrix whose columns are the leading eigenvectors of $\M
U^* \M U^{*\top}$. 
For comparison, we also record $d_2(\hat{\M U},\M U)$ where $\hat{\M U}$
is the $d \times 4$ matrix of leading eigenvectors of $\M Q$ and $\M
Q$ is as defined in Eq.~\eqref{eq:J_pca}. 
Figure~\ref{fig:simu:cov} reports sample means for $d_{2}(\hat{\M U}_g, \M U)$ and $d_2(\hat{\M U}, \M
U)$ as we vary the sampling probability $p$, data dimension $d$,
sample size $n$, and noise level $\sigma$; these sample means
are computed based on $500$ Monte Carlo replicates where we resampled
the matrices $\M U^*, \M E$ and $\bm{\Omega}$ in each replicate. 
Figure \ref{fig:simu:cov} shows that the RSVD estimate $\hat{\M U}_g$
yields nearly-optimal performance (as compared to
$\hat{\M U}$) when $g \geq 3$ for most settings of $p$, $d$, $\sigma$
and $n$. It is only when $d$ is large with respect to $n$
 that $\hat{\M U}_{3}$ have worse performance compared to $\hat{\M U}_{4}, \hat{\M U}_5$ and
$\hat{\M U}$. Finally, $\hat{\M U}_{1}$ is always sub-optimal. These
observations are consistent with the theoretical results presented in Section
\ref{sec:epca}.
\subsection{Additional results for Section~\ref{pca:realdata}}\label{epcasimu}

Additional results for the 68k PBMC scRNA-seq data analysis are provided. We apply Algorithm~\ref{RSRS:pca} to the full data matrix with $g = 3, 4, 5$ and  $\sk = 55,100,300,1000$. The data are projected onto the top $k = 50$ estimated principal components obtained from Algorithm~\ref{RSRS:pca}, and the resulting low-dimensional embeddings are visualized using $t$-SNE in Figure~\ref{fig:pca:real:supp}.

\section{Distributed estimation for multi-layer networks}\label{sec:dist}
Distributed estimation, also known as divide-and-conquer or aggregated
inference, is used in numerous methodological applications including PCA \cite{population_svd,
tang2021integrated,fan2019distributed,chen2021distributed}, regression
\cite{aggregated_inference_huo,ridge_regression_dobriban}, integrative data analysis
\cite{jive}, and is
also a key component underlying federated learning
\cite{federated_survey}. These types of procedures are particularly
important for analyzing large-scale datasets that are scattered across multiple
organizations or computing nodes, where both the computational
complexities and communication costs (including possible privacy constraints)
prevent the transfer of all the raw data to a single location.

We now describe how RSVD can be adapted to tackle estimation of $\M U$ in a distributed setting, thereby
reducing the communication and computation costs. 
For conciseness we will only present, as a notional example, community detection on multi-layer networks. 
More specifically, let $\M A_1, \M A_2, \dots, \M A_{L}$ be a collection of adjacency
matrices for undirected graphs, where $\M A_{\ell}$ is a stochastic blockmodel graph 
with edge probabilities $\M P_{\ell}
= \M Z \M B_{\ell} \M Z^{\top}$; here $\M Z$ denote the community
assignments and $\M B_{\ell}$ denote the block connection
probabilities. The form for $\M P _{\ell}$ indicates that the
community assignments $\M Z$ are shared between all $L$ graphs but the
connection probabilities $\M B_{\ell}$ could be different between any
pair of graphs. The main inference task is to recover, from the $\{\M
A_{\ell}\}_{\ell=1}^{L}$, the community assignments in $\M Z$; see e.g., \cite{paul2020spectral,jing2021community,chen2020global} for a few recent references. In the
event when each individual graph $\M A_{\ell}$ is very sparse, e.g., 
the average degree of $\M A_{\ell}$ is of order $O(1)$, 
consistent estimation of $\M Z$ is statistically infeasible using only
a single graph $\M A_{\ell}$ and thus it is necessary to aggregate the
$\{\M A_{\ell}\}$. One simple approach is based on first forming
$\hat{\M M} = \sum_{\ell = 1}^{L} (\M A_{\ell}^{2} - \M D_{\ell})$
where $\M D_{\ell}$ is either the $n \times n$ diagonal matrix whose diagonal
entries are the vertices degrees in $\M A_{\ell}$ or the diagonal matrix containing the diagonal entries of 
$\M A_{\ell}^2$, then extract
$\hat{\M U}$ as the leading singular vectors of $\hat{\M M}$, and
finally recover $\M Z$ by clustering the rows of $\hat{\M U}$ using
$K$-means or $K$-median clustering; see \cite{lei2020bias,cai2021subspace} for more
details. In particular the subtraction of $\{\mathbf{D}_{\ell}\}$
    corresponds to a bias-removal step and is essential 
when the $\{\M A_{\ell}\}$ are extremely sparse as then
the diagonal entries of $\sum_{i=1}^m \M A_{\ell}^2$ 
are much larger in magnitudes compared to the non-diagonal entries. 
\begin{algorithm}[tp]
\KwIn{$\{{\M A}_{\ell}\} \subset \RR^{n \times n}$, rank $k\geq 1$, sketching
  dimension $\sk \geq k$, power $g \geq 1$.}
  A central server $\mathcal{C}$ generates a seed $s$ and send to each machine $\{\mathcal{M}_{\ell}\}$\; 
     The $\{\mathcal{M}_{\ell}\}$ uses $s$ to generate the same random Gaussian matrix $\M Y_0$\;
    Each $\mathcal{M}_{\ell}$ computes $\M T_{\ell,1} \gets (\M A_{\ell}^{2} - \M D_{\ell}) \M Y_{0}$ and send it back to $\mathcal{C}$\; 
 \For{$t \gets 2 $ \KwTo $g$}{
   $\mathcal{C}$ computes $\M Y_t \gets L^{-1} \sum_{\ell} \M T_{\ell,t}$ and send it to each of the
     $\{\mathcal{M}_{\ell}\}$\;
      Each $\mathcal{M}_{\ell}$ computes $\M T_{\ell,t} \gets (\M A_{\ell}^{2} - \M D_{\ell}) \M Y_{t-1}$ and send it back to $\mathcal{C}$\;
}
   $\mathcal{C}$ computes $\hat{\M U}_g$ as the $k$ leading singular vectors of $\mathbf{Y}_g = L^{-1} \sum_{\ell} \M X_{\ell,g}$\;
   Cluster the rows of $\hat{\M U}_g$ using $K$-means clustering for some choice of $K$ \;
\KwOut{
  Estimated singular vectors $\hat{\M U}_{g}^{(k)}$ 
}
\caption{Spectral clustering for multi-layer networks using distributed RSVD}\label{distributedRSVD}
\end{algorithm}

If each $\M A_{\ell}$ is of dimensions $n \times n$, and the collection $\{\M A_{\ell}\}$ 
is stored on multiple different
machines, then calculating $\hat{\M M}$ requires sending
up to $O(Ln^2)$ bits to a central machine for aggregation, which can be prohibitive when $n$ is large and/or infeasible due to privacy constraints. These issues can be readily addressed by the RSVD estimate $\hat{\M U}_g$ as 
the use of sketching matrices $\M G$ when constructing $\hat{\M M}^{g} \M G$ alleviate the need for
transmitting the $\M A_{\ell}$ directly; see Algorithm~\ref{distributedRSVD} for more details. Let $\mathrm{nnz}(\M
A_{\ell})$ denote the number of non-zero entries in $\M A_{\ell}$ and
$N = \sum_{\ell} \mathrm{nnz}(\M A_{\ell})$ be the total number of
non-zero entries among all $\{\M A_{\ell}\}$.  In
Algorithm~\ref{distributedRSVD}, each iteration of step $5$ involves
$O(L n \sk)$ flops and transfers of $O(L n \sk)$ bits,  each iteration
of step $6$ as well as step~$3$ involves $O(N \sk)$ flops and transfers of $O(L n \sk)$
bits, Step $7$ involves $O(L n \sk + n \sk^2)$ flops, and thus
Algorithm~\ref{distributedRSVD} requires a total of $O(g(Ln + N) \sk + n \sk^2)$
flops and transfer of $O(gLn \sk)$ bits. If $g = O(\log n)$ and $\sk
\ll n$ then the computational complexity of
Algorithm~\ref{distributedRSVD} is considerably smaller than computing
the spectral clustering of $\hat{\M M}$ directly, and require transfer
of at most $O(g L n \sk)$ bits between the individual machines
$\{\mathcal{M}_{\ell}\}$ and the central server $\mathcal{C}$. Finally, we 
note that Algorithm~\ref{distributedRSVD} can also be used for distributed PCA 
\cite{chen2021distributed,fan2019distributed} by simply removing step 8 and 
changing the updates $\M T_{\ell,t} \leftarrow (\M A_{\ell}^2 - \M D_{\ell}) \M Y_t$ in step 3 to $\M T_{\ell, t} \leftarrow 
\M X_{\ell} \M X_{\ell}^{\top} \M Y_{t-1}$ where $\M X_{\ell}$ is the $n_{\ell} \times p$ data matrices whose rows represent observations stored on the $\ell$th machine and whose columns are the feature vectors for these observations. 

Let $\M M = \sum_{\ell} \M P_{\ell}^2$ and suppose that $\mathrm{rk}(\M M) = k_0$ for some finite constant $k_0$ not depending on $L$ and $n$. Suppose also that the average degrees of each $\M A_{\ell}$ is of order 
$\Theta(n \rho_n)$ for some $\rho_n \in [0,1]$ satisfying $n \rho_n = O(1)$. 
If $L^{1/2} n \rho_n \succsim \log^{1/2}(L + n)$, then from the proof of
Theorem~1 in \cite{lei2020bias}, we have
\begin{equation}
  \label{eq:lei_lin1}
    d_2(\hat{\M U}, \M U) \lesssim L^{-1/2}(n\rho_n)^{-1}\log^{1/2}(L+n)+n^{-1}
\end{equation}
with high probability, where $\M U$ contains the eigenvectors corresponding to the non-zero eigenvalues of 
$\hat{\M M}$. Furthermore, under the stronger condition that $L^{1/2} n \rho_n \succsim \log(L + n)$, 
we have by Theorem~1 in \cite{cai2021subspace} that 
\begin{equation}
  \label{eq:cai_multilayer}
  d_{\twoinf}(\hat{\M U}, \M U) \lesssim k_0^{1/2} n^{-1/2} (L^{-1/2} (n \rho_n)^{-1} \log(m+n) + n^{-1})
\end{equation}
with high probability. By Corollary~\ref{co:l2_noise} of our paper, 
$d_{2}(\hat{\M U}_g, \M U)$ and $d_{\twoinf}(\hat{\M U}_g, \M U)$ achieve the same upper bounds as that for 
$d_2(\hat{\M U}, \M U)$ and $d_{\twoinf}(\hat{\M U}, \M U)$ in 
Eq.~\eqref{eq:lei_lin1} and Eq.~\eqref{eq:cai_multilayer}, respectively, whenever 
$\hat{\M U}_g$ is computed using Algorithm~\ref{distributedRSVD} with $\sk \asymp \log n$ and 
\begin{equation}
  \label{eq:condition_g}
  g \geq \frac{2 \log(n/\sk)}{\log(1/(n^{-1} + L^{-1/2} (n \rho_n)^{-1} \log^{1/2}(L+n))}.
\end{equation}
\section{Future directions}\label{sec:discussion}
We   discuss several directions for future work in this section. 
 Firstly, 
as we shown in Section
  \ref{sec:lower_bound}, our $d_2$ and $d_{\twoinf}$ bounds and
  the corresponding phase transition are
  sharp whenever $\hat{\M M}$ satisfies the trace growth condition in Eq.~\eqref{network:setting3}. 
  This condition holds for edge-indepedent random
  graphs with homogeneous edge probabilities and it is of interest to find other 
  inference problem where this condition is also satisfied. Secondly, our 
  discussions in Section~\ref{sec:rgi} and Section~\ref{sec:mc} through Section~\ref{sec:dist}
  focus exclusively on the case where 
$\mathrm{rk}(\M M) = k_0 < n$, and this is because sharp upper bounds for 
$d_{\twoinf}(\hat{\M U}^{(k)}, \M U^{(k)})$ are available mainly when
$k \leq k_0 \ll n$. Extending the 
bounds for $d_{\twoinf}(\hat{\M U}^{(k)}, \M U^{(k)})$ 
to the case where $k_0 = n$ or $k_0 \asymp n$ is an important question 
and furthermore, when combined with the results in this paper, also leads directly to bounds for $d_{\twoinf}(\hat{\M U}_g^{(k)}, \M U^{(k)})$. 
Thirdly, while randomized subspace 
  iterations (as considered in this paper) is one of the most popular approach for RSVD, there are other
  approaches such as those based on Krylov subspaces
  \citep{musco2015randomized}; deriving $2 \to \infty$ norm bounds for RSVD using Krylov subspaces may 
  require different techniques than those presented here.
  Fourthly, many modern dataset are represented as tensors and are
  analyzed using higher-order SVD by flattening the
  tensor into matrices across different
  dimensions and then computing the truncated SVD of the resulting
  matrices \citep{de2000best,zhang2018tensor}. Perturbation analysis of RSVD for noisy tensor data 
is thus of some theoretical and
  practical interest. 

{Finally, for some random graph models, the 
  degree heterogeneity of the nodes can cause the $\ell_2$ norms of the rows of $\M U^{(k)}$ or $\hat{\M U}^{(k)}$ to vary significantly (unlike the
  delocalized setting in Section~\ref{sec:rgi} for which these row-wise $\ell_2$ norms are of order $O(n^{-1/2})$); see e.g., \citet{ke2024optimal} and \citet{cape2024robust}. 
  Under such scenario, we may be interested in row-specific entry-wise bounds as they will aid
  further theoretical analysis for RSVD-based spectral methods.  We thus devote the remaininder of this
  section to a discussion on row-specific perturbation analysis for
  $\hat{\M U}_g^{(k)}$. 

Firstly, if $g$ is sufficiently large then we can show that, 
under {\em minimal}
assumptions, the RSVD-based singular vectors $\hat{\M U}_g^{(k)}$
achieve the same row-specific entrywise error rates in recovering $\M
U^{(k)}$ as the exact singular vectors $\hat{\M U}^{(k)}$. Indeed, 
Theorem \ref{thm3ss:general} requires only a lower bound
on the sketching dimension $\tilde{k}$ and make no assumptions on the
structure of $\hat{\M M}$, while still guaranteeing that if $\hat{\zeta}_k < 1$ then
$d_{2 \to \infty}(\hat{\M U}_g^{(k)}, \hat{\M U}^{(k)})$ decays {\em
exponentially fast} as $g$ increases.  More
specifically, for any $L > 0$, we have
\bee\label{d2infboundrow}
d_{2 \to \infty}(\hat{\M U}_g^{(k)}, \hat{\M U}^{(k)}) = O(n^{-L})
\ee
with high probability, provided that $g \geq C L (\log n)/(\log
\hat{\zeta}_k^{-1})$ for a sufficiently large constant $C$ and
$\tilde{k} \succsim \log n$.  Thus we can decompose the error between
$\hat{\M U}_g^{(k)}$ and $\M U^{(k)}$ into two parts: one between
$\hat{\M U}_g^{(k)}$ and $\hat{\M U}^{(k)}$ (which holds under minimal
assumptions on $\hat{\M M}$) and the other between $\hat{\M U}^{(k)}$
and $\M U^{(k)}$ (which depends on the probabilistic model for
$\hat{\M M}$ as a noisy realization of $\M M$).  This decomposition
naturally extends our results to the case of row-specific entrywise
bound.

In particular, suppose we have a
row-specific bound between $\hat{\M U}^{(k)}$ and $\M U^{(k)}$, that is, for some orthogonal matrix $\M Q_n \in \mathbb{O}_k$ and for any $i \in [n]$ we have
$\bigl\|\M U^{(k)} - \hat{\M U}^{(k)}\M Q_n  \bigr\|_{i, \ell_2} \leq \vartheta_i$
where $\|\M M\|_{i,\ell_2}$ denotes the $\ell_2$ norm of the $i$th row of a matrix $\M M$ and $\vartheta_i$ is a quantity depending possibly on the row index $i$. 
Then, by choosing $g \geq C L (\log n)/(\log \hat{\zeta}_k^{-1})$ and $\tilde{k} \succsim \log n$, we also have for any $i \in [n]$ that
\bee\label{bounddesire}
\bigl\| \hat{\M U}_g^{(k)} - \M U^{(k)} \M Q_n' \bigr\|_{i, \ell_2}
&\leq \bigl\| \hat{\M U}_g^{(k)} - \hat{\M U}^{(k)} \M Q_n\M Q_n'\bigr\|_{i, \ell_2}
+ \bigl\| \hat{\M U}^{(k)} \M Q_n\M Q_n' - \M U^{(k)}  \M Q_n' \bigr\|_{i, \ell_2} \\
&= \bigl\|\hat{\M U}_g^{(k)} - \hat{\M U}^{(k)} \M Q_n\M Q_n' \bigr\|_{i, \ell_2}
+ \bigl\| \M U^{(k)} - \hat{\M U}^{(k)}\M Q_n\bigr\|_{i, \ell_2} \\
&= O(n^{-L}) + \vartheta_i
\ee
with high probability, where $\M Q_n'$ is a $k \times k$ orthogonal matrix
for aligning $\hat{\M U}_g^{(k)}$ and $\hat{\M U}^{(k)}$ to achieve
the rate in Eq.~\eqref{d2infboundrow}. Hence, if $g$ is {\em sufficiently large} then, with high
probability, the row-specific error between $\hat{\M U}_g^{(k)}$ and
$\M U^{(k)}$ is the same (up to an additional negligible term of order $O(n^{-L})$) 
as that between $\hat{\M U}^{(k)}$ and $\M U^{(k)}$.


Secondly, if $g \geq 1$ is arbitrary then, by a more careful inspection of the proof of Theorem~\ref{thm3ss:general}
we can also obtain the following row-specific bound between $\hat{\M U}_g^{(k)}$ and $\M U^{(k)}$; note that this new bound
holds under the same set of {\em minimal assumptions} as that for Theorem~\ref{thm3ss:general}. 


\begin{theorem}
  \label{thm:row}
Consider the setting in Theorem~\ref{thm2ss}. Define $\|\M M\|_{i,\ell_2}$ as the $\ell_2$ norm of the $i$th row of any matrix $\M M$, and for an arbitrary $\delta > 0$, define
\begin{gather}
  \nonumber
  r_{i,2} =  \frac{\sqrt{128} e (k \log \delta^{-1})^{1/2} \hat{\zeta}_k^{\tilde{g}}}{c_{\mathrm{gap}}^2 \sk^{1/2}}  +
\frac{18 n  \|\hat{\M U}^{(k)}\|_{i,\ell_2}\hat{\zeta}_k^{2\tilde{g}}
}{c_{\mathrm{gap}}^2 \sk} + \frac{36 n
  (\log \delta^{-1})^{1/2} \hat{\zeta}_k^{3\tilde{g}}}{c_{\mathrm{gap}}^3 \sk}.
\end{gather}
 Then there exists some $\M Q\in\mathbb{O}_k$ such that uniformly for all $i\in[n]$,
$$
\left\|\hat{\M U}_g^{(k)} -\hat{\M U} ^{(k)}\M Q \right\|_{i,\ell_2} \leq r_{i,2},
$$
with probability at least $1 - 4m \sk \delta - \vartheta - 2e^{-n/2}$.
\end{theorem}
Compared with the $\ell_{\twoinf}$ perturbation bound
\[
r_{\twoinf} = \frac{\sqrt{128} e (k \log \delta^{-1})^{1/2} \hat{\zeta}_k^{\tilde{g}}}{c_{\mathrm{gap}}^2 \tilde{k}^{1/2}} 
+ \frac{18 n \|\hat{\M U}^{(k)}\|_{\twoinf} \hat{\zeta}_k^{2\tilde{g}}}{c_{\mathrm{gap}}^2 \tilde{k}} 
+ \frac{36 n (\log \delta^{-1})^{1/2} \hat{\zeta}_k^{3\tilde{g}}}{c_{\mathrm{gap}}^3 \tilde{k}},
\]
we see that the only difference between $r_{i,2}$ and $r_{2, \infty}$ is the second term in $r_{i,2}$ which uses
$\|\hat{\M U}^{(k)}\|_{i,\ell_2}$ instead of $\|\hat{\M U}^{(k)}\|_{2, \infty}$.
As a result, following a similar argument as in Remark~\ref{rk:illustrate}, $r_{i,2}$ can be sharper than $r_{2,\infty}$ when
$\|\hat{\M U}^{(k)}\|_{i,\ell_2} \ll \|\hat{\M U}^{(k)}\|_{\twoinf}$, especially in the regimes
where $g$ is moderately large and $\tilde{k}$ grows slowly with $n$.  

The proof of Theorem~\ref{thm:row} is straightforward. 
In particular from the proof of Theorem~\ref{thm3ss:general} we have the decomposition
\[
\hat{\M U}_g^{(k)} - \hat{\M U}^{(k)} \M Q_{\check{\M U}_g} = \M T_1 + \M T_2 + \M T_3,
\]
for some $\M Q_{\check{\M U}_g} \in \mathbb{O}_k$, where
$\M T_3 = \check{\M U}_g^{(k)} \bigl( ( \check{\M U}_g^{(k)} )^\T \hat{\M U}_g^{(k)} - \M Q_{\check{\M U}_g} \bigr),
$
and $\check{\M U}_g^{(k)} = \hat{\M U}^{(k)} \check{\M Q}_g$ for some $\check{\M Q}_g \in \mathbb{O}_k$.
Then for Theorem~\ref{thm:row} we bound $\M T_3$ as
\[
\|\M T_3\|_{i,\ell_2} \leq \bigl\|\hat{\M U}^{(k)} \check{\M Q}_g \bigr\|_{i,\ell_2} 
\times \bigl\| ( \check{\M U}_g^{(k)})^\T \hat{\M U}_g^{(k)} - \M Q_{\check{\M U}_g} \bigr\|,
\]
while $\M T_1$ and $\M T_2$ are bounded using their $\ell_{2\to\infty}$ norms as in the proof of Theorem~\ref{thm3ss:general}.

Finally, while it is certainly possible that one can derive row-specific bounds for RSVD sharper than those in
Theorem~\ref{thm:row}, we expect that these bounds also require substantially stronger
assumptions on $\hat{\M M}$. Recall that row-specific bounds
for $\hat{\M U}^{(k)}$ (such as those for
random graphs) typically rely on probabilistic assumptions on
$\hat{\M M}$, for example that the
entries of $\hat{\M M}$ are independent or that $\hat{\M M}$ is a transformation of
such a random matrix; see, e.g., \citet{li2024minimax,ke2024optimal} and
\citet{cape2024robust}. This indicates that, to obtain sharper
row-specific bounds between $\hat{\M U}_g^{(k)}$ and $\hat{\M
U}^{(k)}$, one need to impose additional structural or distribution
assumptions on $\hat{\M M}$, and the theoretical analysis
also has to be tailored to the inference problems of interests; e.g., sharp row-specific
bounds for random graphs may require techniques that are different
from those for matrix completion.  As the primary theoretical goal of
our work is to develop a {\em general} framework for RSVD perturbation under {\it minimal} assumptions,
we leave the more refined row-specific (and problem-centric) analysis for future work.}

\section{Proofs for Section~\ref{sec:thm}}
\label{sec:proof_thm}
\subsection{Primary}
We first recall and introduce a few notations that will be used throughout this section. 
Let $\hat{\M M}$ be a $m \times n$ rectangular or asymmetric matrix and write the SVD of $\hat{\M M}$ as
\begin{equation*}
    \hat{\M M} = \hat{\M U}^{\top} \hat{\bm{{\Sigma}}} \hat{\M V}^{\top} + \hat{\M U}_{\perp}^{\top} \hat{\bm{{\Sigma}}}_{\perp} \hat{\M V}_{\perp}^{\top}
\end{equation*}
where $\hat{\M U}$ and $\hat{\M V}$ are $m \times k$ and $n \times k$ matrices containing the leading $k$ left and right singular vectors, respectively, and $\hat{\bm{\Sigma}}$ are the corresponding singular values; the remaining singular vectors and singular values are denoted by $\hat{\M U}_{\perp}, \hat{\M V}_{\perp}$ and $\hat{\bm{\Sigma}}_{\perp}$. Note that, for ease of notations, we have omitted the dependency on $k$ from all of these matrices. 
Next let $\M G$ be a $n \times \sk$ random Gaussian matrix. Then Algorithm~\ref{RSRS} computes $\hat{\M U}_g$ as the leading $k$ left singular vectors of 
$$\hat{\M M} (\hat{\M M}^{\top} \M M)^{g} \M G = \hat{\M U} \hat{\bm{\Sigma}}^{2g+1} \hat{\M V}^{\top} \M G + \hat{\M U}_{\perp} \hat{\bm{\Sigma}}_{\perp}^{2g+1} \hat{\M V}_{\perp}^{\top} \M G.$$ 
Similarly, if $\hat{\M M}$ is a $n \times n$ symmetric matrix with eigendecomposition 
$$
\hat{\M M} = \hat{\M U} \hat{\bm{\Lambda}} \hat{\M U}^{\top} + \hat{\M U}_{\perp} \hat{\bm{\Lambda}}_{\perp} \hat{\M U}_{\perp}^{\top},
$$
where $\M \Lambda$ contains the top-$k$ eigenvalues in magnitude while $\M \Lambda_\perp$ contains the rest of the eigenvalues. Then Algorithm~\ref{RSRS} computes $\hat{\M U}_g$ as the leading $k$ left singular vectors of
\begin{equation*}
\hat{\M M}^{g} \M G =  \hat{\M U} \hat{\bm{\Lambda}}^{g} \hat{\M U}^{\top} \M G + \hat{\M U}_{\perp} \hat{\bm{\Lambda}}_{\perp}^{g} \hat{\M U}_{\perp}^{\top} \M G =  \hat{\M U} \hat{\bm{\Sigma}}^{g} \hat{\M V}^{\top} \M G + \hat{\M U}_{\perp} \hat{\bm{\Sigma}}_{\perp}^{g} \hat{\M V}_{\perp}^{\top} \M G, 
\end{equation*}
where $\hat{\bm{\Sigma}} = |\hat{\bm{\Lambda}}|$ and the columns of $\hat{\M V}$ are the same as those for $\hat{\M U}$, but with their signs flipped whenever the corresponding eigenvalues in $\hat{\bm{\Lambda}}^g$ are negative.
Consolidating both of the above cases, we can define $\M Y_{g}$ as the matrix
\begin{equation}
\label{eq:Yg_def}
    \M Y_{g} = \hat{\M U} \hat{\bm{\Sigma}}^{\tilde{g}} \hat{\M V}^{\top} \M G +  \hat{\M U}_{\perp} \hat{\bm{\Sigma}}_{\perp}^{\tilde{g}} \hat{\M V}_{\perp}^{\top} \M G = \begin{cases}
\hat{\M M}(\hat{\M M}^\T\hat{\M M})^{\tilde{g}} \M G & \text{$\hat{\M M}$ is symmetric}
\\
 \hat{\M M}^{\tilde{g}} \M G & \text{otherwise}
\end{cases},
\end{equation}
where $\tilde{g} = g$ if $\hat{\M M}$ is symmetric and $\tilde{g} = 2g + 1$ otherwise. Then $\hat{\M U}_g$ is the leading left singular vectors of $\M Y_g$. More specifically, we have the following SVD
\begin{equation*}
  \M Y_g = \hat{\M U}_g \hat{\bm{\Sigma}}_g \hat{\M W}^{\top}_g + \hat{\M U}_{g, \perp} \hat{\bm{\Sigma}}_{g, \perp}
  \hat{\M W}_{g, \perp}^{\top}.
\end{equation*}
where $\hat{\M U}_g$ and $\hat{\M W}_g^{\top}$ are $m \times k$ and $n \times k$ matrices with orthonormal columns. 
We emphasize that $\hat{\bm{\Sigma}}^{g}$ and $\hat{\bm{\Sigma}}_g$ denote different quantities, i.e., $\hat{\bm{\Sigma}}^{g}$ contains the $gt$h powers of the leading $k$ singular values of $\hat{\M M}$ while $\hat{\bm{\Sigma}}_g$ contains the $k$ leading singular values of $\M Y_g$. Finally, our subsequent analysis of $\hat{\M U}_g$ is based on viewing $\M Y_g$ as an additive perturbation of $\hat{\M U} \hat{\bm{\Sigma}}^{\tilde{g}} \hat{\M V}^{\top} \M G$, and thus we also consider the SVD
\begin{equation*}
    \hat{\M U} \hat{\bm{\Sigma}}^{\tilde{g}} \hat{\M V}^{\top} \M G= \check{\M U}_g \check{\bm{\Sigma}}_g \check{\M W}_g^{\top}
\end{equation*}
where $\check{\M U}_g$ and $\check{\M W}_g$ are $m \times k$ and $n \times k$ matrices with orthonormal columns. 
\subsection{Technical lemmas}\label{app:TL}
Before commencing with the proofs of the main results,
we first state some technical lemmas~that will be
used throughout this paper. We start by listing
some basic properties of the $\ell_{\twoinf}$ norm.
\begin{lemma}\label{lm:twoinf}
For any $\M M_1 \in\RR^{d_1\times d_2}, \M M_2 \in \RR^{d_2 \times d_3}$ and $\M M_3 \in \RR^{d_4\times d_1}$, we have
\bee\nonumber
&\|\M M_1 \M M_2\|_{\twoinf} \leq \|\M M_1\|_{\twoinf}\|\M M_2\|;
\\
&\|\M M_3 \M M_1\|_{\twoinf}\leq \|\M M_3\|_\infty \|\M M_1\|_{\twoinf};
\\
&\|\M M_1\|_{\infty}\leq \sqrt{d_2}\|\M M_1\|_{\twoinf}.
\ee
For $\M U_1,\M U_2 \in \mathbb{O}_{d\times d'}$,
let $\M W_{*} := \argmin_{\M W \in \mathbb{O}_{d'}} \|\M U_1 - \M U_2\M W\|_{\F}$.
Then $$\|\M U_1^\T \M U_2 - \M W_{*}\| \leq d_2^2(\M U_1,\M U_2).$$  
\end{lemma}
For a proof of Lemma~\ref{lm:twoinf}, 
see e.g., \cite{cai2018rate} and \cite{cape2019two}. 
The following Lemma~\ref{lm:upper:hp} and Lemma~\ref{lm:basic} provide a collection
of bounds for the quantities depending on the Gaussian sketching matrix $\M G$ 
in the proofs of Theorem~\ref{thm2ss} and Theorem~\ref{thm3ss:general}. For ease of notations, 
we will fix $k$ and thus omit the index $k$ from our matrices, 
e.g., we write $\hat{\M U} $ and $\hat{\M U}_{\perp}$ in place of $\hat{\M U}^{(k)}$ and $\hat{\M U}_{\perp}^{(k)}$. 
\begin{lemma}\label{lm:upper:hp}
  Consider the setting in Theorem~\ref{thm2ss} and let $\M G$ be
  a random Gaussian sketching matrix of dimension $n \times \sk$.
  We then have $\|\M G\|\leq 3 \sqrt{n}$ with probability at least $1 - 2e^{-n/2}$. Furthermore, we also have
  \begin{equation*}
 \|\hat{\M U}_{\perp}\hat{\bm{\Sigma}}_{\perp}^{\ell}\hat{\M V}_{\perp}^\T\M G\|_{\twoinf} \leq \hat{\sigma}_{k+1}^{\ell} (2 \sk \log (1/\delta))^{1/2}
  \end{equation*}
with probability at least $1 - 2m\sk \delta$, where $\ell$ is any arbitrary and given positive integer. 
\end{lemma}
\begin{proof}[Proof of Lemma \ref{lm:upper:hp}]
We first bound $\|\M G\|$. 
Recall the non-asymptotic bound for the spectral norm of a Gaussian random matrix (see e.g. Corollary
5.35 in \cite{vershynin2010introduction}). We then have 
\bee\label{|G|2}
\|\M G\|\leq n^{1/2} + \sk^{1/2} + t,
\ee
with probability at least $1 - 2\exp(-t^2/2)$, and hence we can take $t = n^{1/2}$ in Eq.~\eqref{|G|2}
to obtain $\|\M G\|\leq
3\sqrt{n}$ with probability at least $1 - 2\exp(-n/2)$.

We now bound $\|\hat{\M U}_{\perp}\hat{\bm{\Sigma}}_{\perp}^{\ell}\hat{\M V}_{\perp}^\T\M G\|_{\twoinf}$. First, we have
\bee\label{ULUG}
\|\hat{\M U}_{\perp}\hat{\bm{\Sigma}}_{\perp}^{\ell}\hat{\M V}_{\perp}^\T\M G\|_{\twoinf} \leq \sk^{1/2}
\max_{j\in[\sk]}\|\hat{\M U}_{\perp}\hat{\bm{\Sigma}}_{\perp}^{\ell}\hat{\M V}_{\perp}^\T\bds g_j\|_{\max},
\ee
where $\bm{g}_j$ denote the $j$th column of $\M G$.
Next we note that \[\hat{\M U}_{\perp}\hat{\bm{\Sigma}}_{\perp}^{\ell}\hat{\M
V}_{\perp}^\T\bm{g}_j \sim \mathcal{N}(\bm{0}, \hat{\M
U}_{\perp}\hat{\bm{\Sigma}}_{\perp}^{2 \ell}\hat{\M U}_{\perp}^\T), \quad \text{and} \quad
\|\hat{\M U}_{\perp}\hat{\bm{\Sigma}}_{\perp}^{2\ell}\hat{\M
  U}_{\perp}^\T\|_{\max} \leq \|\hat{\bm{\Sigma}}_{\perp}\|^{2\ell} \leq \hat{\sigma}_{k+1}^{2\ell}. \]
Now for any Gaussian random variable $\xi$ with $\mathrm{Var}(\xi) \leq c$ we have
\bee\label{G:hp}
\p\big(|\xi|<t\big) \geq 1 - 2\exp(-t^2/(2c)),
\ee  for all $t > 0$; see e.g Section~2.5.1 of
\cite{vershynin2018high}.
Therefore, by taking a union over all $m$ elements of $\hat{\M U}_{\perp}\hat{\bm{\Sigma}}_{\perp}^{\ell}\hat{\M
V}_{\perp}^\T\bm{g}_j$ and then over all $j \leq \tilde{k}$, we have
\begin{equation}
  \begin{split}
  \label{eq:hatU_perp_lambda_twoinf1}
\p\Big(\max_{j \in [\tilde{k}]} \|\hat{\M U}_{\perp}\hat{\bm{\Sigma}}_{\perp}^{\ell}\hat{\M V}_{\perp}^\T \bm{g}_j \|_{\max}
<\hat{\sigma}_{k+1}^{\ell} \sqrt{2 \log  (1/ \delta)}\, \Big)  
& \geq 1 - 2m\tilde{k} \delta.
    \end{split}
\end{equation}
Combining Eq.~\eqref{eq:hatU_perp_lambda_twoinf1} and \Eq\eqref{ULUG} 
we obtain
\bee\nonumber
\|\hat{\M U}_{\perp}\hat{\bm{\Sigma}}_{\perp}^{\ell}\hat{\M V}_{\perp}^\T\M G\|_{\twoinf} <
\hat{\sigma}_{k+1}^{\ell} (2 \sk \log (1/ \delta))^{1/2}
\ee
with probability at least $1 - 2m\sk \delta$. 
\end{proof}
\begin{lemma}\label{lm:basic}
  Consider the setting in Theorem~\ref{thm2ss} where
$\M G$ is a random $n \times \sk$ Gaussian matrix with 
$$
\sk \geq (1 - c_{\mathrm{gap}})^{-2} \Big\{k + (8 k \log (1/\vartheta))^{1/2} + 2 \log(1/\vartheta) \Big\}
$$ for some arbitrary $c_{\mathrm{gap}} \in (0,1)$ and some arbitrary $\vartheta > 0$. 
Let $g$ be an arbitrary positive integer and define $\M Y_g$ as in Eq.~\eqref{eq:Yg_def}
with $\tilde{g} = g$ if $\hat{\M M}$ is symmetric and $\tilde{g} = 2g+1$ otherwise. Then
\bee\nonumber
\sigma_{k}^2(\hat{\M V}^{\top} \M G) \geq c_{\mathrm{gap}}^2 \sk \quad \text{and} \quad
\sigma_{k}^2(\M Y_{g}) \geq \sigma^2_k(\hat{\bm{\Sigma}}^{\tilde{g}} \hat{\M V}^\T\M G)  = \sigma_k^2(\check{\bm{\Sigma}}_{g}) \geq c_{\mathrm{gap}}^2 \sk \hat{\sigma}_k^{2\tilde{g}}.
 \ee
 with probability at least $1 - \vartheta$.
\end{lemma}

\begin{proof}
We first bound $\sigma_k^2(\hat{\M V}^\T\M G)$.
As $\M G$ is a $n \times \sk$ matrix whose entries are independent
standard normals, $\hat{\M V}^{\top} \M G$ is a $k \times \sk$ matrix
whose entries are also independent standard normals. 
By Theorem~II.13 in \cite{davidson2001local} we have
\bee
\label{PSI:bound}
\mathbb{P}\Bigl(\sigma_{k}(\hat{\M V}^{\top} \M G) \leq  \sk^{1/2} -
k^{1/2} - t \Bigr)
\leq e^{-t^2/2}.
\ee
Now choose an arbitrary $c_{\mathrm{gap}} \in (0, 1)$. Then
$$\sk^{1/2} - k^{1/2} - t \geq c_{\mathrm{gap}} \sk^{1/2} \Longleftrightarrow \sk \geq (1 - c_{\mathrm{gap}})^{-2} (k^{1/2} + t)^2.$$
Letting $t = (2 \log (1/\vartheta))^{1/2}$ yields
\begin{equation*}
  \mathbb{P}\left(\sigma_k(\hat{\M V}^{\top} \M G) \geq c_{\mathrm{gap}} \sk^{1/2}\right) \geq 1 - \vartheta
\end{equation*}
for $\sk \geq (1 - c_{\mathrm{gap}})^{-2} \{k + (8 k \log (1/\vartheta))^{1/2} + 2 \log (1/\vartheta)\}$. 
Next, for $\hat{\bm{\Sigma}}^{\tilde{g}} \hat{\M V}^{\top} \M G$ we have
\begin{equation}
  \label{asy:argue2}
\sigma_{k}^2(\hat{\bm{\Sigma}}^{\tilde{g}} \hat{\M V}^{\top} \M G) = \lambda_{k}(\M G^{\top} \hat{\M V} \hat{\bm{\Sigma}}^{2\tilde{g}} \hat{\M V}f^{\top} \M G) \geq \lambda_k(\M G^{\top} \hat{\M V} \hat{\M V}^{\top} \M G) \times \hat{\sigma}_k^{2 \tilde{g}} 
\end{equation}
Substituting the above bound for $\sigma_{k}(\hat{\M V}^{\top} \M G)$ into  \Eq\eqref{asy:argue2} we obtain
\bee\label{lklower:1}
\sigma_k^2\big(\hat{\bm{\Sigma}}^{\tilde{g}} \hat{\M V}^{\top} \M G) \geq c_{\mathrm{gap}}^2 \sk \hat{\sigma}_{k}^{2 \tilde{g}} 
\ee
with probability at least $1 - \vartheta$. Finally,
\begin{equation*}
  \sigma_{k}^2(\M Y_{g}) = \lambda_{k}(\M Y_{g}^{\top} \M Y_{g}) = \lambda_{k}(\M G^{\top}(\hat{\M V}^{\top} \hat{\bm{\Sigma}}^{2\ell} \hat{\M V} + \hat{\M V}_{\perp} \hat{\bm{\Sigma}}_{\perp}^{2\tilde{g}} \hat{\M V}_{\perp}^{\top}) \M G) \geq \lambda_{k}(\M G^{\top} \hat{\M V} \hat{\bm{\Sigma}}^{2\tilde{g}} \hat{\M V}^{\top} \M G) = \sigma_{k}^2(\hat{\bm{\Sigma}}^{\tilde{g}} \hat{\M V}^{\top} \M G)
\end{equation*}
and hence $\sigma_{k}^2(\M Y_{g}) \geq c_{\mathrm{gap}}^2 \sk \hat{\sigma}_{k}^{2\tilde{g}}$ with probability at least $1 -\vartheta$. 
\end{proof}
\subsection{Proof of Theorem~\ref{thm2ss}}\label{sec:pf:thm2}
For simplicity of notations, we will fix a value of $k$ and thus omit it from our notations, e.g.,
 we use $\hat{\M U}_g$ and $\hat{\M U}$ in place of $\hat{\M U}^{(k)}_g$ and $\hat{\M U}^{(k)}$, 
respectively. Furthermore let $\hat{\sigma}_i$ denote the $i$th largest singular value of $\hat{\M M}$. 
Now recall the definition of $\M Y_{g}$ from Eq.~\eqref{eq:Yg_def}. 
As $\hat{\M V}^\T\M G$ is a $k\times \sk$ matrix whose entries are iid $\mathcal{N}(0,1)$ with $\sk \geq k$, we have 
$\mathrm{rk}(\hat{\M V}^{\top} \M G) = k$ almost surely; see e.g. \cite{edelman1991distribution} for a justification of this claim. We thus have  
\[\mathrm{rk}(\hat{\M U} \hat{\bm{\Sigma}}^{\tilde{g}} \hat{\M V}^\T\M G) = \mathrm{rk}(\hat{\M U} \hat{\bm{\Sigma}}^{\tilde{g}}) = \mathrm{rk}(\hat{\M U})\]
almost surely. 
Let $\check{\M U}_g$ be the left singular vectors of 
$\hat{\M U}\hat{\bm{\Sigma}}^{\tilde{g}}\hat{\M V}^\T\M G$. 
As $\check{\M U}_g$ and $\hat{\M U}$ are both $m \times k$ matrices with orthonormal columns, $\mathrm{rk}(\check{\M U}_g) = \mathrm{rk}(\hat{\M U})$, 
and $\mathcal{C}(\check{\M U}_g) \subset \mathcal{C}(\hat{\M U})$,
we conclude that
\bee\label{UUW}
\check{\M U}_g = \hat{\M U}\check{\M Q}_g
\ee for some $\check{\M Q}_g\in\mathbb{O}_{k}$.
We can now view $\M Y_{g}$ 
as the perturbed version of $ \hat{\M U}\hat{\bm{\Sigma}}^{\tilde{g}}\hat{\M U}^\T\M G$. 
Then by the Wedin $\sin$--$\Theta$ theorem (see pages $262$ and $267$ of \cite{stewart_sun}), we have
\bee\label{e499}
                                  \vertiii{\sin \Theta(\hat{\M U}_g,\check{\M U}_g)} 
\leq \frac{ \vertiii{\hat{\M U}_{\perp}\hat{\bm{\Sigma}}_{\perp}^{\tilde{g}}\hat{\M
    V}_{\perp}^\T\M G}}{\sigma_{k}(\M Y_{g})} 
\leq \frac{ \vertiii{\hat{\bm{\Sigma}}_{\perp}^{\tilde{g}}}   
\cdot \|\M G\|}{\sigma_{k}(\M Y_{g})}
\leq \frac{3  n^{1/2} \vertiii{\hat{\bm{\Sigma}}_{\perp}^{\tilde{g}}}
}{c_{\mathrm{gap}}  \sk^{1/2} \hat{\sigma}_{k}^{\tilde{g}}} 
\ee
with probability at least $1 - \vartheta - 2e^{-n/2}$, where the last inequality follows
from Lemma~\ref{lm:upper:hp} and Lemma~\ref{lm:basic}.
\qed
\subsection{Proof of Theorem~\ref{thm3ss:general}}\label{sec:pf:ta3}
For simplicity of notations, we will fix a value of $k$ and thus omit it from our notations. 
Recall the definition of $\M Y_{g}$ in Eq.~\eqref{eq:Yg_def}. Next recall the SVD of $\M Y_{g}$ and $\hat{\M U} \hat{\bm{\Sigma}}^{\tilde{g}} \hat{\M V}^{\top} \M G$,
\bee
  \label{eq:svd_Zell}
  \M Y_{g} &= \hat{\M U}_g \hat{\bm{\Sigma}}_{g} \hat{\M W}_g^{\top} + \hat{\M U}_{g, \perp} \hat{\bm{\Sigma}}_{g, \perp} \hat{\M W}_{g, \perp}^{\top}, 
\\
 \hat{\M U} \hat{\bm{\Sigma}}^{\tilde{g}} \hat{\M V}^{\top} \M G &=
      \check{\M U}_g \check{\bm{\Sigma}}_{g} \check{\M W}_g^{\top},
\ee
where $\hat{\M U}_g$ and $\hat{\M W}_g$ are $m \times k$ and $\sk \times k$ matrices whose columns are the $k$ leading left and right singular vectors of $\M Y_{g}$, respectively. 
Now define
$$\M Q_{\check{\M U}_g} = \argmin_{\M Q\in\mathbb{O}_{k}}\|\hat{\M U}_g - \check{\M U}_g \M Q\|_\F \quad \text{and} \quad
\M Q_{\check{\M W}_g} = \argmin_{\M Q\in\mathbb{O}_{k}}\|\hat{\M W}_g - \check{\M W}_g \M Q\|_\F.$$
We derive the following decomposition:
\bee
\label{eq:Procrustean_1}
\hat{\M U}_g - \check{\M U}_g\M Q_{\check{\M U}_g}  
&= (\M I_m - \check{\M U}_g\check{\M U}_g^\T) \hat{\M U}_g
 + \check{\M U}_g (\check{\M U}_g^{\top} \hat{\M U}_g - \M Q_{\check{\M U}_g})
  \\ &= (\M I_m - \check{\M U}_g\check{\M U}_g^\T) \M Y_{g} \hat{\M W}_g \hat{\M\Sigma}_{g}^{-1} +
       \check{\M U}_g (\check{\M U}_g^{\top} \hat{\M U}_g - \M Q_{\check{\M U}_g})
\\ &=
(\M I_m - \hat{\M U}\hat{\M U}^\T) \M Y_{g} \hat{\M W}_g \hat{\M\Sigma}_{g}^{-1} + \check{\M U}_g (\check{\M U}_g^{\top} \hat{\M U}_g - \M Q_{\check{\M U}_g})
\\
&= \hat{\M U}_{\perp} \hat{\bm{\Sigma}}^{\tilde{g}}_{\perp}\hat{\M V}^\T_{\perp} \M G \hat{\M W}_g \hat{\bm{\Sigma}}_{g}^{-1} + \check{\M U}_g (\check{\M U}_g^{\top} \hat{\M U}_g - \M Q_{\check{\M U}_g}) 
\\
& = \underbrace{\hat{\M U}_{\perp}\hat{\bm{\Sigma}}_{\perp}^{\tilde{g}}\hat{\M V}_{\perp}^\T\M G\check{\M W}_g\M Q_{\check{\M W}_g} \hat{\bm{\Sigma}}_{g}^{-1}}_{\M T_1} 
\\
&\quad + \underbrace{\hat{\M U}_{\perp}\hat{\bm{\Sigma}}_{\perp}^{\tilde{g}}\hat{\M V}_{\perp}^\T\M G(\hat{\M W}_g - \check{\M W}_g \M Q_{\check{\M W}_g})\hat{\bm{\Sigma}}_{g}^{-1}}_{\M T_2}
 \\
&\quad + \underbrace{\check{\M U}_g(\check{\M U}_g^\T \hat{\M U}_g - \M Q_{\check{\M U}_g})}_{\M T_3}.
\ee 
     where $\M I_{m}$ denotes the $m \times m$ identity matrix. Recall that
     $\hat{\bm{\Sigma}}^{\tilde{g}}$ and $\hat{\bm{\Sigma}}_{g}$ represent different quantities, i.e., $\hat{\bm{\Sigma}}^{\tilde{g}}$ is the $\tilde{g}$-th power of
     $\hat{\bm{\Sigma}}$, the leading $k$ singular values of $\hat{\M M}$, while $\hat{\bm{\Sigma}}_{g}$
     are the leading
     singular values of either $\M Y_{g} = \hat{\M M}^{g} \M G$ or $\M Y_{g} = \hat{\M M} (\hat{\M M}^{\top} \hat{\M M})^{g} \M G$;
     see the discussion at the beginning of Section~\ref{sec:proof_thm}. We now bound the $2 \to \infty$ norms of $\M T_1, \M T_2,$ and $\M T_3$. 
\subsubsection{Bounding \texorpdfstring{$\|\M T_1\|_{\twoinf}$}{TEXT}}
 Let $\hat{\M L}_{\M G}\hat{\M D}_{\M G}\hat{\M R}_{\M G}^{\top}$ be the SVD of $\hat{\M G}^{\top} \hat{\M V}$ 
 where the columns of $\hat{\M L}_{\M G}\in\mathbb{O}_{\sk \times k}$ and $\hat{\M R}_{\M G}\in\mathbb{O}_{k \times k}$
     are the left and right singular vectors, respectively and the diagonal of $\hat{\M D}_{\M G}\in \RR^{k\times k}$ contains the singular values of $\M G^\T\hat{\M V}$. Recalling Eq.~\eqref{eq:svd_Zell}, we have $\check{\M W}_g = \M G^\T \hat{\M V}\hat{\bm{\Sigma}}^{\tilde{g}}\hat{\M U}^\T\check{\M U}_g\check{\M\Sigma}_{g}^{-1}$ where
     $\hat{\bm{\Sigma}}^{\tilde{g}}\hat{\M U}^\T\check{\M U}_g\check{\M\Sigma}_{g}^{-1}$ is invertible;
     see Eq.~\eqref{UUW}. Therefore
     $\check{\M W}_g$ and $\M G^\T \hat{\M V}$ share the same column space and hence
$$
\check{\M W}_g = \hat{\M L}_{\M G} \M Q_{\M G} = \M G^\T \hat{\M V} \hat{\M R}_{\M G}\hat{\M D}_{\M G}^{-1} \M Q_{\M G},
$$
for some $\M Q_{\M G}\in\mathbb{O}_k$. Using the above form for $\check{\M W}_g$, we can rewrite $\M T_1$ as  
\bee
     \M T_1 &= \hat{\M U}_{\perp}\hat{\bm{\Sigma}}_{\perp}^{\tilde{g}}\hat{\M V}_{\perp}^\T\M G\M G^\T \hat{\M V}
              \hat{\M R}_{\M G}\hat{\M D}_{\M G}^{-1} \M Q_{\M G}\M Q_{\check{\M V}_g} \hat{\M\Sigma}_{g}^{-1}\\
 &= \hat{\sigma}_{k+1}^{\tilde{g}}\cdot \bds\Upxi \cdot \hat{\M R}_{\M G}\hat{\M D}_{\M G}^{-1} \M Q_{\M G}\M Q_{\check{\M V}_g} \hat{\M\Sigma}_{g}^{-1},
\ee
where $\hat{\sigma}_{k+1}$ is the $k+1$ largest singular value of
              $\hat{\M M}$ and $$\bds\Upxi = \hat{\M U}_{\perp}(\hat{\bm{\Sigma}}_{\perp}/\hat{\sigma}_{k+1})^{\tilde{g}}
              \hat{\M V}_{\perp}^\T\M G\M G^\T \hat{\M V}.$$ Let $\hat{\bm{m}}_i$ and $\hat{\bm{v}}_j$ denote the $i$th row of $\hat{\M U}_{\perp}(\hat{\bm{\Sigma}}_{\perp}/\hat{\sigma}_{k+1})^{\tilde{g}}\hat{\M V}^\T_{\perp}$ and $j$th column of $\hat{\M V}$,
              respectively. Then the $ij$th element of $\bds\Upxi$ is of the form
\begin{equation}
\label{eq:t1_sum}
 \Upxi_{ij} =  \sum_{\ell = 1}^{\tilde{k}} \hat{\bm{m}}_i^{\top} \bm{g}_{\ell} \bm{g}_{\ell}^{\top} \hat{\bm{v}}_j
\end{equation}
where $\bm{g}_{\ell}$ is the $\ell$th column of $\M G$. 
We further denote $\Upxi_{ij\ell} = \hat{\bm{m}}_i^{\top}
\bm{g}_{\ell} \bm{g}_{\ell}^{\top}
\hat{\bm{v}}_j$.
Now note that $\hat{\M
  U}_{\perp}(\hat{\bm{\Sigma}}_{\perp}/\hat{\sigma}_{k+1})^{\tilde{g}}\hat{\M V}^\T_{\perp}
\hat{\M V} = \bds 0_{n\times k}$ and thus
$\hat{\bm{m}}_i^{\top} \hat{\bm{v}}_j = 0$ for any $(i,j)\in[n -
              k]\times [k]$. As the columns of $\M G$ are independent of $\hat{\M M}$,
              for any $\ell\in[\sk]$ and any $(i,j)\in[n -
k]\times [k]$ we have
\[
\E\big[\Upxi_{ij\ell}\big] = \mathbb{E}\big[\hat{\bm{m}}_i^{\top}
\bm{g}_{\ell} \bm{g}_{\ell}^{\top}
\hat{\bm{v}}_j\big] = \hat{\bm{m}}_i^{\top}
\mathbb{E}\big[\bm{g}_{\ell} \bm{g}_{\ell}^{\top}\big]
\hat{\bm{v}}_j = \hat{\bm{m}}_i^{\top} \hat{\bm{v}}_j = 0.
\]
Furthermore, $\|\hat{\bm{m}}_i\| \leq \|\hat{\bm{v}}_j\| = 1$ for all
$(i,j)\in[n - k]\times [k]$, which implies that $\hat{\bm{m}}_i^{\top}
\bm{g}_{\ell}$ and $\hat{\bm{v}}^\T_j \bm{g}_{\ell}$ are two
independent Gaussian random variables with variances bounded by $1$.

Denote $\|\cdot\|_{\psi_1}$ and $\|\cdot\|_{\psi_2}$ as the
sub-exponential and sub-Gaussian norms of a random variable. 
Then $\|\hat{\bm{m}}_i^{\top} \bm{g}_{\ell} \|_{\psi_2}\leq 2$,
$\|\hat{\bm{v}}^\T_j \bm{g}_{\ell}\|_{\psi_2} \leq 2$ and by
Lemma~2.7.7 in \cite{vershynin2018high}, 
\[
  \|\Upxi_{ij\ell}\|_{\psi_1} = \|\hat{\bm{m}}_i^{\top}
  \bm{g}_{\ell} \bm{g}_{\ell}^{\top}
  \hat{\bm{v}}_j\|_{\psi_1} \leq \|\hat{\bm{m}}_i^{\top}
  \bm{g}_{\ell} \|_{\psi_2}\|\hat{\bm{v}}^\T_j
  \bm{g}_{\ell}\|_{\psi_2} \leq 4\]
Therefore, given $\hat{\M M}$, $\Upxi_{ij}$ in Eq.~\eqref{eq:t1_sum}
is  a sum of {\em independent} random variables with sub-exponential
norms bounded by $4$. Then by a standard application of Bernstein's inequality (see, e.g., Theorem 2.8.1 of \cite{vershynin2018high}), we have  
\bee\nonumber
\mathbb{P}\Bigl(\max_{i,j}\big|\Upxi_{ij}\big| \geq t \Bigr)
\leq 2 m \sk \exp\Bigl[-\min\Bigl\{\frac{t^2}{128e^2 \sk}
   ,\frac{t}{16e} \Bigr\}\Bigr] 
\ee
Taking $t = (128 \sk \log (1/\delta))^{1/2} e$ 
and noting that, from our condition on $\sk$ in Theorem~\ref{thm3ss:general}, we have
$\sk \geq 2 \log (1/ \delta)$ and hence $t/(16e) \geq \log (1/\delta)$. This then implies
\bee\nonumber
\mathbb{P}\Bigl(\max_{i,j}\big|\Upxi_{ij}\big| \geq (128 \sk \log (1/\delta))^{1/2} e \Bigr)
\leq 2m\sk \delta
\ee
We therefore have
\[
\hat{\sigma}_{k+1}^{\tilde{g}} \| \bds\Upxi\|_{\twoinf} \leq
\hat{\sigma}_{k+1}^{\tilde{g}} k^{1/2} \max_{i,j}|\Upxi_{ij}| \leq
 \sqrt{128} e \hat{\sigma}_{k+1}^{\tilde{g}} (k\sk \log (1/\delta))^{1/2},
\]
with probability at least $1 - 2m\sk \delta$. 
Finally, by Lemma~\ref{lm:basic}, we have $\|\hat{\M D}_{\M G}^{-1}\| \leq
c_{\mathrm{gap}}^{-1} \sk^{-1/2}$ and $\|\hat{\M \Sigma}_{g}^{-1}\| \leq
c_{\mathrm{gap}}^{-1} \sk^{-1/2} \hat{\sigma}_{k}^{\tilde{g}}$ with probability
at least $1 - \vartheta$. Lemma~\ref{lm:twoinf} then implies
\begin{equation}
  \label{eq:mt1_final}
\|\M T_1\|_{\twoinf} \leq  \hat{\sigma}_{k+1}^{\tilde{g}} \times \|
\bds\Upxi\|_{\twoinf} \times  \|\hat{\M D}_{\M G}^{-1}\| \times \|\hat{\M
  \Sigma}_g^{-1}\|\leq \frac{ \sqrt{128} e (k \log (1/\delta))^{1/2} \hat{\sigma}_{k+1}^{\tilde{g}}}{c_{\mathrm{gap}}^2 \sk^{1/2}\hat{\sigma}_k^{\tilde{g}}}
\end{equation}
with probability at least $1 - 2m\sk\delta - \vartheta$. 
\subsubsection{Bounding \texorpdfstring{$\|\M T_2\|_{\twoinf}$}{TEXT}}
We first observe that, by Lemma~\ref{lm:twoinf}, Lemma~\ref{lm:upper:hp} and Lemma~\ref{lm:basic}, we have
\bee
\label{eq:T2_twoinf_part1}
\|\M T_2\|_{\twoinf} &\leq \|\hat{\M U}_{\perp}\hat{\bm{\Sigma}}_{\perp}^{\tilde{g}}\hat{\M V}_{\perp}^\T\M G\|_{\twoinf}  \|\hat{\M W}_g - \check{\M W}_g \M Q_{\check{\M W}_g}\|  \|\hat{\M\Sigma}_{g}^{-1}\| 
  \\ &\leq c_{\mathrm{gap}}^{-1} \hat{\zeta}_k^{\tilde{g}} (2 \log (1/\delta))^{1/2}
       \|\hat{\M W}_g - \check{\M W}_g \M Q_{\check{\M W}_g}\|
\ee
       with probability at least $1 - 2m\sk \delta - \vartheta$,
       where $\hat{\zeta}_k = {\hat{\sigma}_{k+1}}/{\hat{\sigma}_k}$. 


We first bound $\|\hat{\M W}_g
- \check{\M W}_g \M Q_{\check{\M W}_g}\|$ using the rate-optimal bound in \cite{cai2018rate} which carefully analyzes the perturbation of asymmetric matrices. Let
$$\bds X = \hat{\M U} \hat{\bds\Sigma}^{\tilde{g}} \hat{\M V}^{\top} \M G, \quad \bds Z = \hat{\M U}_{\perp} \hat{\M \Sigma}_{\perp}^{\tilde{g}} \hat{\M V}_{\perp}^{\top} \M G, \quad \hat{\bds X} = \bds X + \bds Z = \M Y_{g}$$
     Recall that the columns of $\check{\M W}_g$ and $\hat{\M W}_g$ are the leading
     right singular vectors of $\bds X$ and $\hat{\bds X}$, respectively, while the columns of $\hat{\M U}_g$ are the leading left singular vectors of $\hat{\M X}$.
Eq.~\eqref{UUW} then implies
\begin{gather*}
\bds Z_{12} := \check{\M U}_g \check{\M U}_g^{\top} \bds Z (\M I - \check{\M W}_g \check{\M W}_g^{\top}) = \bds 0, \quad  \bds Z_{21}:= (\M I - \check{\M U}_g \check{\M U}_g^{\top}) \bds Z \check{\M W}_g \check{\M W}_g^{\top} = \bds Z \check{\M W}_g \check{\M W}_g^{\top} \\
\check{\M U}_g^{\top} \hat{\bds X} \check{\M W}_g  = \check{\M U}_g^{\top} \M X \check{\M W}_g =
\check{\M \Sigma}_{g}, \quad \check{\M U}_{g, \perp}^{\top} \hat{\bds X} \check{\M W}_{g, \perp} = \check{\M U}_{g, \perp}^{\top} \bds Z \check{\M W}_{g, \perp}
\end{gather*}
where $\check{\M U}_{g, \perp}^{\top}$ and $\check{\M W}_{g, \perp}$
are matrices with orthonormal columns for the basis of $\M I - \check{\M U}_g \check{\M U}_g^{\top}$ and $\M I - \check{\M W}_g \check{\M W}_g^{\top}$, respectively. 
Let $\alpha = \sigma_{k}(\check{\M \Sigma}_{g})$ and $\beta = \|\check{\M
  U}_{g, \perp}^{\top} \hat{\bds X} \check{\M W}_{g, \perp}\| \leq
\|\bds Z\|$. We then have
\begin{equation}
  \label{eq:hat_mv_checkmv_bd1}
    \begin{split}
      \|\hat{\M W}_g - \check{\M W}_g \M Q_{\check{\M W}_g}\|& \leq \frac{2 \sqrt{2} \|\bds Z\|^2}{\alpha^2}.
\end{split}
\end{equation}
Indeed, if $\alpha^2 \geq 2 \|\M Z\|^2$ then
Eq.~\eqref{eq:hat_mv_checkmv_bd1} follows from Theorem~1 in
\cite{cai2018rate}, namely
\begin{equation*}
    \begin{split}
\|\sin \Theta(\hat{\M W}_g, \check{\M W}_g)\| &\leq \frac{\alpha \|\bds Z_{12}\| + \beta \|\bds Z_{21}\|}{\alpha^2 - \beta^2 - \min\{\|\bds Z_{12}\|^2, \|\bds Z_{21}\|^2\}} 
\\ 
&= \frac{\beta \|\bds Z_{21}\|}{\alpha^2 - \beta^2} \leq \frac{\|\bds
  Z\|^2}{\alpha^2 - \beta^2} \leq \frac{\|\bds Z\|^2}{\alpha^2 -
  \|\bds Z\|^2} \leq \frac{2 \|\bds Z\|^2}{\alpha^2}.
\end{split}
\end{equation*}
If $\alpha^2 \leq 2 \|\M Z\|^2$ then
Eq.~\eqref{eq:hat_mv_checkmv_bd1} follows from the trivial bound
$\|\sin \Theta(\hat{\M W}_g, \check{\M W}_g)\| \leq 1$. 
We therefore have, by Lemma~\ref{lm:upper:hp} and Lemma~\ref{lm:basic}, that
\begin{equation}
  \label{eq:T2_twoinf_part2}
  \|\hat{\M W}_g - \check{\M W}_g \M Q_{\check{\M W}_g}\|  \leq \frac{2 \sqrt{2} \|\hat{\bm{\Sigma}}_{\perp}\|^{2\tilde{g}}
    \cdot \|\M G\|^2}{\sigma_k^2(\check{\bm{\Sigma}}_{g})} \leq
 \frac{2 \sqrt{2} \hat{\sigma}_{k+1}^{2\tilde{g}} \times 9n}{c_{\mathrm{gap}}^2 \sk \hat{\sigma}_{k}^{2\tilde{g}}} \leq 
  \frac{18\sqrt{2} n}{c_{\mathrm{gap}}^2 \sk} \hat{\zeta}_k^{2\tilde{g}} 
\end{equation}
with probability at least $1 - \vartheta - 2e^{-n/2}$. 

Combining Eq.~\eqref{eq:T2_twoinf_part1} and Eq.~\eqref{eq:T2_twoinf_part2}, we obtain
\[\|\M T_2\|_{2 \to \infty}
\leq \frac{36 n (\log (1/\delta))^{1/2} \hat{\zeta}_k^{3\tilde{g}}}{c_{\mathrm{gap}}^3 \sk }\] 
with probability at least $1 - 2m\sk\delta - \vartheta - 2e^{-n/2}$. 
\subsubsection{Bounding \texorpdfstring{$\|\M T_3\|_{\twoinf}$}{TEXT}}
First recall Eq.~\eqref{UUW}. Then by Lemma~\ref{lm:twoinf} and Eq.~\eqref{e499} we have
\begin{equation*}\begin{split}\|\M T_3\|_{\twoinf} &\leq \|\check{\M U}_g\|_{\twoinf}  \|\check{\M U}_g^\T \hat{\M U}_g - \M Q_{\check{\M U}_g}\| 
\\ &\leq \|\check{\M U}_g\|_{\twoinf}  d^2_2(\hat{\M U}_g, \check{\M U}_g)
\leq \|\hat{\M U}\|_{\twoinf} \frac{18n \hat{\zeta}_k^{2\tilde{g}}}{c_{\mathrm{gap}}^2 \sk}
\end{split}
\end{equation*}
with probability at least $1 - \vartheta - 2e^{-n/2}$. 
\subsubsection{Putting all pieces together}
Combining the bounds for $\|\M T_1\|_{\twoinf}$ through $\M \|\M
T_3\|_{\twoinf}$ we obtain
\bee
\label{eq:t1_t3_combined}
d_{\twoinf}(\hat{\M U}_g, \hat{\M U}) &=d_{\twoinf}(\hat{\M U}_g, \check{\M U}_g)
\\
& \leq\| \hat{\M U}_g - \check{\M U}_g\M Q_{\check{\M U}_g}\|_{\twoinf} \\
&\leq \frac{ \sqrt{128} e (k \log (1/\delta))^{1/2} \hat{\zeta}_{k}^{\tilde{g}}}{c_{\mathrm{gap}}^2 \sk^{1/2}} +  \frac{36 n (\log(1/\delta))^{1/2} \hat{\zeta}_k^{3\tilde{g}}}{c_{\mathrm{gap}}^3 \sk } +
  \|\hat{\M U}\|_{2 \to \infty}
\frac{18 n \hat{\zeta}_k^{2\tilde{g}}}{{{{c_{\mathrm{gap}}^2 \sk}}}}
\ee
with probability at least $1 - 4m\sk \delta - \vartheta - 2e^{-n/2}$.  \qed 
\subsubsection{Bounding $\|(\hat{\M U}_g \hat{\M U}_g^{\top} - \hat{\M U} \hat{\M U}^{\top}) \hat{\M M}\|_{\max}$}
Let $\M R_g = \hat{\M U}_g - \check{\M U}_g \M Q_{\check{\M U}_g} = \hat{\M U}_g - \hat{\M U} \check{\M Q}_{g} 
\M Q_{\check{\M U}_g}$. As $\check{\M U}_g \M Q_{\check{\M U}_g} (\check{\M U}_g \M Q_{\check{\M U}_g})^{\top} = \hat{\M U} \hat{\M U}^{\top}$, we have
\begin{equation*}
  \begin{split}
    (\hat{\M U}_g \hat{\M U}_g^{\top}  - \hat{\M U} \hat{\M U}^{\top}) \hat{\M M}  &= \M R_g \M R_g^{\top} \hat{\M M} +  \M R_g (\check{\M Q}_{g} \M Q_{\check{\M U}_g})^{\top} \hat{\M U}^{\top} \hat{\M M} + \hat{\M U}
\check{\M Q}_{g} \M Q_{\check{\M U}_g} \M R_g^{\top} \hat{\M M}.
  \end{split}
\end{equation*}
We therefore have
\begin{equation*}
  \begin{split}
    \|(\hat{\M U}_g \hat{\M U}_g^{\top} - \hat{\M U} \hat{\M U}^{\top}) \hat{\M M}\|_{\max} &\leq \|\M R_{g}\|_{\twoinf} \bigl(\|\hat{\M M}^{\top} \M R_{g}\|_{\twoinf} +  \|\hat{\M M}^{\top} \hat{\M U}\|_{\twoinf}\bigr) + \|\hat{\M U}\|_{\twoinf} \times \|\hat{\M M} \M R_{g}\|_{\twoinf}
    \\ & \leq \|\M R_{g}\|_{\twoinf} \bigl(\|\hat{\M M}^{\top} \M R_{g}\|_{\twoinf} +  \|\hat{\M V}\|_{\twoinf}
         \times \|\hat{\bm{\Sigma}}\|\bigr) + \|\hat{\M U}\|_{\twoinf} \times \|\hat{\M M}^{\top} \M R_{g}\|_{\twoinf}. 
  \end{split}
\end{equation*}
where we have used the fact that for any matrices $\M A$ and $\M B$ for which $\M A \M B$ is well defined, $\|\M A \M B\|_{\max} \leq \|\M A\|_{\twoinf} \times \|\M B^{\top}\|_{\twoinf}$. 
By Eq.~\eqref{eq:t1_t3_combined}, we have
$\|\hat{\M R}_g\|_{\twoinf} \leq r_{\twoinf}$ with probability at least $1 - 4m\sk \delta - \vartheta - 2e^{-n/2}$. 
We now bound $\|\hat{\M M}^{\top} \M R_{g}\|_{\twoinf}$. Recalling Eq.~\eqref{eq:Procrustean_1} and 
that $ \check{\M U}_g \check{\M U}_g^{\top} = \hat{\M U} \hat{\M U}^{\top}$, we have
\begin{equation*}
  \begin{split}
    \hat{\M M}^{\top} \M R_g &= \hat{\M M}^{\top} (\M I - \hat{\M U} \hat{\M U}^{\top}) \hat{\M R}_g    +
                               \hat{\M M}^{\top} \hat{\M U} \hat{\M U}^{\top} \hat{\M R}_g \\
    &= \hat{\M M}^{\top} (\M I - \hat{\M U} \hat{\M U}^{\top}) \hat{\M U}_g +  
                               \hat{\M V} \hat{\bm{\Sigma}} \check{\M Q}_g^{\top} (\check{\M U}_g^{\top} \hat{\M U}_g - \M Q_{\check{\M U}_g}) \\
                             &= \hat{\M V}_{\perp} \hat{\bm{\Sigma}}_{\perp} \hat{\M U}_{\perp}^{\top} \M Y_{g}
                               \hat{\M W}_{g} \hat{\bm{\Sigma}}_g^{-1} + \hat{\M V} \hat{\bm{\Sigma}} \check{\M Q}_g^{\top} (\check{\M U}_g^{\top} \hat{\M U}_g - \M Q_{\check{\M U}_g}) \\
                             &= \underbrace{\hat{\M V}_{\perp} \hat{\bm{\Sigma}}_{\perp}^{\tilde{g}+1} \hat{\M V}_{\perp}^{\top} \M G \check{\M W}_{g} \M Q_{\check{\M W}_g} \hat{\bm{\Sigma}}_g^{-1}}_{\M T_1'} + \underbrace{\hat{\M V}_{\perp} \hat{\bm{\Sigma}}_{\perp}^{\tilde{g}+1} \hat{\M V}_{\perp}^{\top} \M G (\hat{\M W}_g - \check{\M W}_{g} \M Q_{\check{\M W}_g}) \hat{\bm{\Sigma}}_g^{-1}}_{\M T_2'} \\ & + \underbrace{\hat{\M V} \hat{\bm{\Sigma}} \check{\M Q}_g^{\top} (\check{\M U}_g^{\top} \hat{\M U}_g - \M Q_{\check{\M U}_g})}_{\M T_3'}  
  \end{split}
\end{equation*}
Following the same derivations as that for $\M T_1, \M T_2$ and $\M T_3$ in the proof of
Theorem~\ref{thm3ss:general}, we have
\begin{gather*}
  \|\M T_1'\|_{\twoinf} \leq \frac{\sqrt{128} e (k \log (1/\gamma))^{1/2} \hat{\sigma}_{k+1} \hat{\zeta}_k^{\tilde{g}}}
  {c_{\mathrm{gap}}^2 \sk^{1/2}}, 
 \quad \|\M T_2'\|_{\twoinf} \leq \frac{36 n (\log (1/\gamma))^{1/2} \hat{\sigma}_{k+1} \hat{\zeta}_k^{3 \tilde{g}}}{c_{\mathrm{gap}}^3 \sk}, \\
  \|\M T_3'\|_{\twoinf} \leq \frac{18 n \hat{\sigma}_1 \hat{\zeta}_k^{2g}}{c_{\mathrm{gap}}^2 \sk} \|\hat{\M V}\|_{\twoinf},
\end{gather*}
with probability at least $1 - 4 n \sk \gamma - \vartheta - 2e^{-n/2}$. 
We therefore have
\[\|\hat{\M M} \M R_g\|_{\twoinf} \leq \|\M T_1'\|_{\twoinf} + \|\M T_2'\|_{\twoinf} +
  \|\M T_3'\|_{\twoinf} \leq \hat{\sigma}_1 \tilde{r}_{\twoinf} \]
with probabilty at least $1 - 4 n \sk \gamma - \vartheta - 2e^{-n/2}$. In summary we have
\begin{equation*}
  \begin{split}
    \|(\hat{\M U}_g \hat{\M U}_g^{\top} - \hat{\M U} \hat{\M U}^{\top}) \hat{\M M}\|_{\max}
    & \leq r_{\twoinf} \times (\hat{\sigma}_1 \tilde{r}_{\twoinf} + \|\hat{\M V}\|_{\twoinf} \times
          \hat{\sigma}_1 \bigr) + \|\hat{\M U}\|_{\twoinf} \times \hat{\sigma}_1 \tilde{r}_{\twoinf} \\ &\leq                                                                                                          \hat{\sigma}_1 \bigl(r_{\twoinf} \tilde{r}_{\twoinf} + \|\hat{\M U}\|_{\twoinf} \tilde{r}_{\twoinf} +                                                                                                      \|\hat{\M V}\|_{\twoinf} r_{\twoinf}\bigr)
  \end{split}
\end{equation*}
with probability at least $1 - 4\sk(m \delta + n \gamma) - \vartheta -  2e^{-n/2}$.
\qed

\subsection{Proof of Corollary~\ref{co:l2_noise}}
For ease of notations we will omit the index $k_0$ from our matrices. 
If $\mathrm{rk}(\M M) = k_0$ then $\sigma_{k_0+1} = 0$ and, by the Davis-Kahan theorem
\citep{davis1970rotation}, $d_2(\hat{\M U}, \M U) \leq E_n/\sigma_{k_0}$. 
Next, from Corollary~\ref{co:arbitrary_noise} with $k = k_0$ and $\vartheta = n^{-3}$ we have
\begin{equation*}
  d_2(\hat{\M U}_g, \hat{\M U}) \leq d_{2}(\hat{\M U}, \M U) + 3 \sqrt{2} c_{\mathrm{gap}}^{-1}(n/\sk)^{1/2}  \zeta_{k_0}^g \leq \left\{1 +  
  3 \sqrt{2} c_{\mathrm{gap}}^{-1} (n/\sk)^{1/2} \zeta_{k_0}^{g-1}\right\}\frac{E_n}{\sigma_{k_0}} 
\end{equation*}
with probability at least $1 - 2n^{-3}$,  provided that $\sk \geq (1 - c_{\mathrm{gap}})^{-2} \{k_0 + \sqrt{24 k_0 \log n} + 6 \log n\}$. 
Let $ g_\iota =  \frac{\log\left\{3\sqrt{2}c_{\mathrm{gap}}^{-1}({n/\sk})^{1/2} \iota^{-1}\right\}}{ \log(1/\zeta_{k_0})} $. 
  Then
  \bee\nonumber
   3\sqrt{2}c_{\mathrm{gap}}^{-1} (n/\sk)^{1/2} \zeta_{k_0}^{g - 1} & =   \iota\cdot \iota^{-1}3\sqrt{2}c_{\mathrm{gap}}^{-1} (n/\sk)^{1/2} \cdot \zeta_{k_0}^{g - 1}  
=\iota\cdot \zeta_{k_0}^{-g_\iota} \cdot \zeta_{k_0}^{g - 1}  
=\iota \cdot \zeta_{k_0}^{g -1 - g_\iota},
  \ee
which yields \eqref{eq:co_main_ws}.
 
\par
Eq.~\eqref{thm2*:resopt:main}
for the strong signal regime $ \sigma_{k_0}/E_n \succsim n^{\epsilon}$ 
can be shown using the same argument as described above. First we have for some $c \in(0,1)$, $\color{black}\zeta^{-1}_{k_0} \geq \sigma_{k_0}/(2E_n)\geq cn^{\epsilon} $ when $n$ is sufficiently large. Thus we have
\bee\nonumber
g_* = \frac{\log(n )-\log(\sk)}{2 \log(1/\zeta_{k_0})}
 \leq \frac{\log(n )-\log(\sk)}{2\epsilon \log(n) +2\log c}
\leq (2\epsilon)^{-1}
\ee
for sufficiently large $n$, where the final inequality 
follows from the fact that $\log(\sk) \succsim \log\log n$. Let $\iota =  3\sqrt{2} c^{-1}_{\mathrm{gap}}$ in \eqref{eq:co_main_ws}, then we have
\bee\nonumber
g_\iota = \frac{\log(  {n/\sk}  ) }{ 2\log(1/\zeta_{k_0})} =  g_*,\ee
and when $g\in [g_*,1 + g_*]$,
\bee\nonumber
\big(1 + \iota \cdot \zeta_{k_0}^{g   - g_\iota - 1} \big)\frac{E_n}{\sigma_{k_0}} & =\frac{E_n}{\sigma_{k_0}} + 3\sqrt{2} c^{-1}_{\mathrm{gap}} \cdot \zeta_{k_0}^{g   - g_* - 1}\cdot\frac{E_n}{\sigma_{k_0}}
 \\
&\leq \frac{E_n}{\sigma_{k_0}} +6\sqrt{2} c^{-1}_{\mathrm{gap}} \cdot \zeta_{k_0}^{g   - g_* } 
\leq \frac{E_n}{\sigma_{k_0}} +12\sqrt{2} c^{-1}_{\mathrm{gap}} \cdot \left(\frac{E_n}{\sigma_{k_0}}\right)^{g   - g_* }
= O\Bigl\{ \Bigl(\frac{E_n}{\sigma_{k_0}}\Bigr)^{g   - g_* }\Bigr\}.
\ee
When $g \geq 1 + g_* +\Delta$ for any $\Delta = \omega(\log^{-1} n)$, we have
\bee\nonumber
\big(1 + \iota \cdot \zeta_{k_0}^{g   - g_\iota - 1} \big)\frac{E_n}{\sigma_{k_0}}&\leq  \frac{E_n}{\sigma_{k_0}} +6\sqrt{2} c^{-1}_{\mathrm{gap}} \cdot \zeta_{k_0}^{g   - g_* } = \frac{E_n}{\sigma_{k_0}} +6\sqrt{2} c^{-1}_{\mathrm{gap}} \cdot \zeta_{k_0}^{1 + \Delta } 
\\
&\leq  \frac{E_n}{\sigma_{k_0}} + 12 \sqrt{2} c^{-1}_{\mathrm{gap}} \cdot  c^{-\Delta}n^{-\Delta\epsilon} \frac{E_n}{\sigma_{k_0}} 
 = (1 + o(1))E_n/\sigma_{k_0} .
\ee
Combining the above results yield \eqref{thm2*:resopt:main}.

\qed 

\subsection{Proof of Corollary~\ref{co:l2inf_noise}} 
Similar to the proof of Corollary~\ref{co:l2_noise}, we will omit the index $k_0$ from our matrices. From Corollary~\ref{co:arbitrary_noise} with $k = k_0$ and $\vartheta = n^{-3}$ and $\delta = n^{-5}$ we have
\begin{equation*}
  d_{\twoinf}(\hat{\M U}_g, \M U) \leq C_* \Bigl(\frac{(k \log n)^{1/2} \zeta_{k_0}^{g}}{\sk ^{1/2}} + \frac{n u_{\twoinf} \zeta_{k_0}^{2g}}{\sk} + \frac{n (\log n)^{1/2} \zeta_{k_0}^{3g}}{\sk}\Bigr) := r_{\twoinf}
\end{equation*}
with probability at least $1 - 5n^{-3} - 2 e^{-n/2}$, where $\zeta_{k_0} = E_n/(\sigma_{k_0} - E_n)$ and 
$C_*$ is a constant depending only on $c_{\mathrm{gap}}$; 
recall that, from the statement of Corollary~\ref{co:l2_noise} we already assume $c_{\mathrm{gap}} > 0$ is a fixed but arbitrary constant. 
Now define
\begin{equation*}
  t_1 = 1 + \frac{\log(nk/\sk) + \log \log n}{2 \log(1/\zeta_{k_0})}, \quad t_2 = \frac{1}{2} + \frac{\log(n^{3/2}\|\M U\|_{\twoinf} /\sk)}{2 \log (1/\zeta_{k_0})}, \quad t_3 = \frac{1}{3} + \frac{\log(n^{3/2}/\sk) + \tfrac{1}{2} \log \log n}{3 \log (1/\zeta_{k_0})}. 
\end{equation*}
Then for $g \geq t_1$ we have
\begin{equation*}
  g \log(1/\zeta_{k_0}) \geq \log(1/\zeta_{k_0}) + \frac{1}{2}\Big\{\log(nk/\sk) + \log \log n\Big\}
\end{equation*}
which then implies
\begin{equation*}
\zeta_{k_0}^{g} \leq \zeta_{k_0} \times \frac{\sk^{1/2}}{(kn \log n)^{1/2}} \Longrightarrow
  \frac{(k \log n)^{1/2} \zeta_{k_0}^{g}}{\sk^{1/2}} \leq n^{-1/2} \zeta_{k_0} = O(n^{-1/2} E_n/\sigma_{k_0}).
\end{equation*}
Similarly, if $g \geq t_2$ then
\begin{equation*}
  \frac{n u_{\twoinf} \zeta_{k_0}^{2g}}{\sk} \leq \zeta_{k_0} \times u_{\twoinf} \times \frac{1}{n^{1/2} \|\M U\|_{\twoinf}} \leq 2 n^{-1/2} \zeta_{k_0} = O(n^{-1/2} E_n/\sigma_{k_0})
\end{equation*}
where the second inequality follows from the assumption 
$d_{\twoinf}(\hat{\M U}, \M U) \leq \|\M U\|_{\twoinf}$ in the statement of Corollary~\ref{co:l2inf_noise}. If $g \geq t_3$ then
\begin{equation*}
  \frac{n (\log n)^{1/2} \zeta_{k_0}^{3g}}{\sk} \leq n^{-1/2} \zeta_{k_0} = O(n^{-1/2} E_n/\sigma_{k_0}).
\end{equation*}
As $\sigma_{k_0} > 2 E_n$, we have $(\sigma_{k_0} - E_n)/E_n \geq \sigma_{k_0}/(2E_n)$ and thus
$$t_1 = 1 + \frac{\log(nk/\sk) + \log \log n}{2 \log (1/\zeta_{k_0})} \leq 1 +  
\frac{\log(nk/\sk) + \log \log n}{2 \log(\tfrac{1}{2} \sigma_{k_0}/E_n)},
$$
and similarly for $t_2$ and $t_3$. 
In summary, if $g \geq g_* \geq \max\{t_1, t_2, t_3\}$ then 
\[d_{\twoinf}(\hat{\M U}_g, \hat{\M U}) = O(n^{-1/2} E_n/\sigma_{k_0}) \]
with probability at least $1 - 5n^{-3} - 2e^{-n/2}$. 

Furthermore, as $\M M$ and $\hat{\M M}$ are both symmetric, we have
$u_{\twoinf} = v_{\twoinf} \leq 2\|\M U\|_{\twoinf}$ and $r_{\twoinf} = \tilde{r}_{\twoinf}$. Then by Eq.~\eqref{eq:entrywise_co_general}, for $g \geq g_*$ we have
\begin{equation*}
  \begin{split}
  \|\hat{\M U}_g \hat{\M U}_g^{\top} \hat{\M M} - \M M\|_{\max} &\leq \|\hat{\M U} \hat{\M U}^{\top} \hat{\M M} - \M M\|_{\max} + (\sigma_1 + E_n) (r_{\twoinf}^2 + 4 \|\M U\|_{\twoinf} r_{\twoinf}) \\
                                                                &\leq \|\hat{\M U} \hat{\M U}^{\top} \hat{\M M} - \M M\|_{\max} +  O((\sigma_1 + E_n)(n^{-1} (E_n/\sigma_{k_0})^2 + n^{-1/2} (E_n /\sigma_{k_0}) \|\M U\|_{\twoinf}) \\
    & \leq \|\hat{\M U} \hat{\M U}^{\top} \hat{\M M} - \M M\|_{\max} + O(\kappa n^{-1/2} E_n \|\M U\|_{\twoinf})
  \end{split}
\end{equation*}
with probability at least $1 - 9n^{-3} - 2e^{-n/2}$, where the last inequality follows from the fact that $\|\M U\|_{\twoinf} \geq k_0^{1/2} n^{-1/2}$. 
Finally, if $\sigma_{k_0} = \omega(E_n)$ and $g \geq g_* + c$ for any fixed but arbitrary $c > 0$ then
$d_{\twoinf}(\hat{\M U}_g, \hat{\M U}) = O(n^{-1/2} (E_n/\sigma_{k_0})^{1 + c}) = o(n^{-1/2} E_n/\sigma_{k_0})$ which also implies
\begin{equation*}
  \begin{split}
  \|\hat{\M U}_g \hat{\M U}_g^{\top} \hat{\M M} - \M M\|_{\max}  
                                                                &\leq \|\hat{\M U} \hat{\M U}^{\top} \hat{\M M} - \M M\|_{\max} +  o((\sigma_1 + E_n)(n^{-1} (E_n/\sigma_{k_0})^2 + n^{-1/2} (E_n /\sigma_{k_0}) \|\M U\|_{\twoinf})\bigr)\\
    & \leq \|\hat{\M U} \hat{\M U}^{\top} \hat{\M M} - \M M\|_{\max} + o(\kappa n^{-1/2} E_n \|\M U\|_{\twoinf})
  \end{split}
\end{equation*}
\section{Proofs for Section~\ref{sec:rgi}, Section~\ref{sec:mc} and Section~\ref{sec:epca}}  \label{sec:pfa}
\subsection{Proof of Theorem~\ref{thm:lower}}
\label{sec:proof_lower}

As $d_{\twoinf}(\hat{\M U}_g, \M U) \geq n^{-1/2} d_{2}(\hat{\M U}_g, \M U)$ always holds,
we will only derive the lower bound for $d_2(\hat{\M U}_g, \M U)$. It is sufficient to lower bound $\|(\M I - \hat{\M U} \hat{\M U}^{\top}) \hat{\M U}_g \hat{\M U}_g^{\top}\|$, which serves as a lower bound for $d_2^2(\hat{})$.
Note that
\begin{equation*}
    \begin{split}
      (\M I - \hat{\M U} \hat{\M U}^{\top}) \hat{\M U}_g \hat{\M U}_g^{\top} &= \hat{\M U}_{\perp} \hat{\M U}_{\perp}^{\top} \hat{\M M}^{g} \M G \M T (\M T^{\top} \M G^{\top} \hat{\M M}^{2g} \M G \M T)^{-1} \M T^{\top} \M G^{\top} \hat{\M M}^{g}
    \end{split}
\end{equation*}
where $\M T$ is a $\tilde{k} \times k$ matrix with orthonormal columns such that the column space of $\hat{\M M}^{g} \M G \M T$ is the same as that for $\hat{\M U}_g$. Such $\M T$ always exists as the column space of $\hat{\M M}^g\M G\M T$ must include the column space of $\hat{\M U}_g$ with $\sk \geq k$. We then have
\begin{equation*}
    \begin{split}
      \|(\M I - \hat{\M U} \hat{\M U}^{\top}) \hat{\M U}_g \hat{\M U}_g^{\top}\| &=
      \lambda_{\max}\bigl(\hat{\M U}_{\perp}^{\top} \hat{\M M}^{g} \M G \M T (\M T^{\top} \M G^{\top} \hat{\M M}^{2g} \M G \M T)^{-1} \M T^{\top} \M G^{\top} \hat{\M M}^{g} \hat{\M U}_{\perp}\bigr)
      \\ &= \lambda_{\max}\bigl( \hat{\bm{\Lambda}}_{\perp}^{g} \hat{\M U}_{\perp}^{\top} \M G \M T (\M T^{\top} \M G^{\top} \hat{\M M}^{2g}
           \M G \M T)^{-1} \M T^{\top} \M G^{\top}  \hat{\M U}_{\perp} \hat{\bm{\Lambda}}_{\perp}^{g} \bigr)
      \\ & \geq \frac{\lambda_{\max}\bigl( \hat{\bm{\Lambda}}_{\perp}^{g} \hat{\M U}_{\perp}^{\top} \M G \M T \M T^{\top} \M G^{\top}  \hat{\M U}_{\perp} \hat{\bm{\Lambda}}_{\perp}^{g}\bigr)}
           {\|\M T^{\top} \M G^{\top} \hat{\M M}^{2g} \M G \M T\|} \\
    &\geq \frac{\lambda_{\max}\bigl(\M T^{\top} \M G^{\top} \hat{\M U}_{\perp} \hat{\bm{\Lambda}}_{\perp}^{2g} \hat{\M U}_{\perp}^{\top} \M G \M T\bigr)}{\|\M G^{\top} \hat{\M M}^{2g} \M G\|}. 
    \end{split}
\end{equation*}
Now $\hat{\M U}_{\perp} \hat{\bm{\Lambda}}_{\perp}^{2g} \hat{\M U}_{\perp}^{\top} = \hat{\M M}^{2g} - \hat{\M U} \hat{\bm{\Lambda}}^{2g} \hat{\M U}^{\top}$ and hence
\begin{equation*}
    \begin{split}
    \|(\M I - \hat{\M U} \hat{\M U}^{\top}) \hat{\M U}_g \hat{\M U}_g^{\top}\|  
    & \geq \frac{\lambda_{\max}(\M T^{\top} \M G^{\top} \hat{\M M}^{2g} \M G \M T) - \|\hat{\bm{\Lambda}}\|^{2g} \times \|\M T^{\top} \M G^{\top} \hat{\M U}\|^2}{\|\M G^{\top} \hat{\M M}^{2g} \M G\|}
   \\ & \geq \frac{\lambda_{\max}(\M G^{\top} \hat{\M M}^{2g} \M G) - \|\hat{\bm{\Lambda}}\|^{2g} \times \|\M G^{\top} \hat{\M U}\|^2}{\|\M G^{\top} \hat{\M M}^{2g} \M G\|}
    \end{split}
\end{equation*}
where the last inequality follows from the fact that $\M T$ is a partial isometry mapping the column space of $\hat{\M M}^{g} \M G$ to
$\hat{\M U}_g$ and hence 
$\lambda_{\max}(\M T \M G^{\top} \hat{\M M}^{2g} \M G \M T) = \lambda_{\max}(\M G^{\top} \hat{\M M}^{2g} \M G).$
Furthermore, as the diagonal entries of a matrix is majorized by its eigenvalues, we have
\begin{equation*}
    \|(\M I - \hat{\M U} \hat{\M U}^{\top}) \hat{\M U}_g \hat{\M U}_g^{\top}\| \geq \frac{\bigl(\max \bm{g}_i^{\top} \hat{\M M}^{2g} \bm{g}_i\bigr) - \|\hat{\bm{\Lambda}}\|^{2g} \times \|\M G^{\top} \hat{\M U}\|^2}{\|\M G^{\top} \hat{\M M}^{2g} \M G\|}
\end{equation*}
where $\bm{g}_i$ is the $i$th column of $\M G$ and the maximum is taken over all $i \leq \tilde{k}$.
We next recall our assumption on $\hat{\M M}$, namely that 
$\mathbb{E}[\mathrm{tr} \, \hat{\M M}^{2g}] \geq c_g(n^{g+1} \rho_n^{g} + (n \rho_n)^{2g})$
for some constant $c > 0$ and $n \rho_n \succsim n^{\beta}$ for some $\beta > 0$.  
Fix a $g < \beta^{-1}$. By Markov's inequality, there exists a constant $C > 0$ such that
\begin{equation}
\label{eq:lower_bound_trace}
    \mathrm{tr} \, \hat{\M M}^{2g} \geq C c_g n^{1 + \beta g}
\end{equation}
with probability at least $1 - p_0$. Furthermore, recall that we had assumed $\|\M M\| \asymp n \rho_n$ 
and $\lambda_{k_0}/E_n \asymp (n \rho_n)^{1/2}$,
which together implies $\|\hat{\M M}\| \asymp n \rho_n $ with probability at least $1 - n^{-3}$. 
Let $\mathcal{E}$  be the event that Eq.~\eqref{eq:lower_bound_trace} holds together with $\|\hat{\M M}\| \asymp n \rho_n$; note that $\mathbb{P}(\mathcal{E}) \geq 1 - p_0 - n^{-3}$. 
Then by the Hanson-Wright inequality \citep{hanson_wright}, we have
\begin{equation*}
  \begin{split}
\p\bigl(\bm{g}_{i}^{\top} \hat{\M M}^{2g} 
    \bm{g}_{i}  \geq \mathrm{tr} \, \hat{\M M}^{2g} -t \mid \mathcal{E} \bigr) &\geq 1 -
                                                                \exp\Big\{-C_2\min\Big(\frac{t^2}
                                                                {\|\hat{\M M}^{2g}\|^2_\F},\frac{t}{\|\hat{\M M}^{2g}\|}\Big)\Big\}
\\
&\geq 1 - \exp\Big\{-C_3\min\Big(\frac{t^2}{n^{2g + 1}\rho_n^{2g}},\frac{t}{(n\rho_n)^{2g}}\Big)\Big\},
    \end{split}
  \end{equation*}
for any $t \geq 0$; here $C_2 \geq 0$ and $C_3 \geq 0$ are constants
not depending on $n$. We thus have, by a union bound over all $i  \leq \sk$ that
\begin{equation*}
\p\bigl(\max_{i} \bm{g}_{i}^{\T} \hat{\M M}^{2g}
\bm{g}_{i}  \geq \mathrm{tr}\, \hat{\M M} -t \mid \mathcal{E} \bigr) \geq 1 - \sk \exp\Big\{-C_3\min\Big(\frac{t^2}{n^{2g + 1}\rho_n^{2g}},\frac{t}{(n\rho_n)^{2g}}\Big)\Big\}.
\end{equation*}
We now choose $t = C_4n^{g + 1/2}\rho_n^g\log^{1/2}{n}$ for some sufficiently large $C_4
\geq 0$. Then, conditional on $\mathcal{E}$, we have for sufficiently large $n$ that
\bee\label{lower:final2}
\max_{i \leq \sk} \bm{g}_{i}^{\top} \hat{\M M}^{2g}
\bm{g}_{i}  \geq \mathrm{tr} \hat{\M M}^{2g} - C_4n^{g + 1/2}\rho_n^g\log^{1/2}n 
\geq \frac{1}{2} \mathrm{tr} \hat{\M M}^{2g}
\ee
with high probability. 
Furthermore, from Lemma~\ref{lm:basic} we have
\begin{equation*}
  \begin{split}
\|\M G^{\top} \hat{\M M}^{2g} \M G\| &\leq \|\M G^{\top} \hat{\M U}\|^{2} \times \|\hat{\bm{\Lambda}}\|^{2g} + \|\M G\|^2 \times \|\hat{\bm{\Lambda}}_{\perp}\|^{2g}
\\ & \precsim \tilde{k} \|\bm{\Lambda}\|^{2g} + n \hat{\lambda}_{k+1}^{2g} 
\\
&\precsim 
    \tilde{k} (n \rho_n)^{2g} + n^{g+1} \rho_n^{g} 
\\
&\precsim (\log n)(n\rho_n)^{2g} + n^{g+1} \rho_n^{g}
  \end{split}
  \end{equation*}
with probability at least $1 - n^{-3}$, where the final inequality is because $\sk \asymp \log n = o(n^{\epsilon})$ for any $\epsilon > 0$.  
Combining the above bounds we have (conditional on $\mathcal{E}$),
\begin{equation}
    \label{eq:lower_bound_main}
    \|(\M I - \hat{\M U} \hat{\M U}^{\top}) \hat{\M U}_g \hat{\M U}_g^{\top}\| \geq
    \frac{\tfrac{1}{2} \mathrm{tr} \, \hat{\M M}^{2g} - C_0 \tilde{k} \|\hat{\bm{\Lambda}}\|^{2g}}{C_1 n^{g+1} \rho_n^{g} + C_1 (\log n)(n\rho_n)^{2g}}
\end{equation}
with probability at least $1 - 2n^{-3}$. Finally, we also have for sufficiently large $n$ that
\begin{equation}
  \label{eq:lower_tr_Mhat}
  \begin{split}
    \frac{1}{2} \mathrm{tr} \hat{\M M}^{2g} - C_0 \tilde{k} \|\hat{\bm{\Lambda}}\|^{2g} &                                                                                    
                                                                                                                                                                                   \geq \frac{1}{4} \mathrm{tr} \hat{\M M}^{2g} \geq \frac{1}{4} c_g n^{g+1} \rho_n^{g} 
  \end{split}
\end{equation}
Substituting Eq.~\eqref{eq:lower_tr_Mhat} into Eq.~\eqref{eq:lower_bound_main}, and 
unconditioning with respect to $\mathcal{E}$, we have for all $g < \beta^{-1}$, $(\log n)(n\rho_n)^{2g}\prec n^{g + 1}\rho_n^g$ and thus
\begin{equation*}
  \begin{split}
    \|(\M I - \hat{\M U} \hat{\M U}^{\top}) \hat{\M U}_g \hat{\M U}_g^{\top}\| &
   \geq \frac{\mathrm{tr} \, \hat{\M M}^{2g} - C_0 \tilde{k} \|\hat{\bm{\Lambda}}\|^{2g}}{C_5 n^{g+1} \rho_n^{g}}
    \geq \frac{ c_g n^{g+1} \rho_n^{g}}{4 C_5 n^{g+1} \rho_n^{g}}
\geq C_{LB} c_g 
  \end{split}
\end{equation*}
with probability at least $1 - p_0 - 2n^{-3}$, where $C_{LB}$ is a constant depending only on $p_0$.  Similarly, if $\beta^{-1}\leq g < 1 +\beta^{-1}$, we have $(\log n)(n\rho_n)^{2g}\succ n^{g + 1}\rho_n^g$ and thus
$$
 \|(\M I - \hat{\M U} \hat{\M U}^{\top}) \hat{\M U}_g \hat{\M U}_g^{\top}\| 
   \geq \frac{\mathrm{tr} \, \hat{\M M}^{2g} - C_0 \tilde{k} \|\hat{\bm{\Lambda}}\|^{2g}}{C_6 (\log n)(n\rho_n)^{2g}}
    \geq \frac{ c_g n^{g+1} \rho_n^{g}}{(\log n)(n\rho_n)^{2g}} \geq C_{LB}c_g(\log n)^{-1}(n)^{  1- g\beta}.
$$
Note that   $C_{LB}$ and $c_g$ written in the theorem are the squared roots of the corresponding constants written here.
\qed
\subsection{Proof of Theorem~\ref{thm:mc} and Corollary~\ref{co:entrywise}}
Let $\hat{\M M} = p^{-1}
  \hat{\M T}$ and $\mathbf{M} = \E\big[p^{-1} \hat{\M T}\big] 
= \E[p^{-1} \{{\mathcal{P}}_{\M \Omega}(\M T+ \M N)\}] = \M T$. Now
define $\M U$, $\hat{\M U}$ and $\hat{\M U}_g$ accordingly, where, for simplicity of notations we have dropped the index $k_0$ from these matrices.  
Finally let $\hat{\M T}_g = p^{-1} \hat{\M U}_g \hat{\M U}_g^{\top} \hat{\M T}$. 
\subsubsection{Bounding $\|\mathbf E\|$}
We have
\bee\label{E:decom:mc}
\M E = \hat{\M M} - \M M = \underbrace{\frac{1}{p}\bigl\{{\mathcal{P}}_{\M\Omega}(\M T)  - p\M T\bigr\}}_{\M E_1} + \underbrace{\frac{1}{p}\mathcal{P}_{\M \Omega}(\M N)}_{\M E_2}.
\ee
Now $\M E_1$ is a random symmetric matrix whose upper
triangular entries are independent mean $0$ random variables with
$$\frac{1}{\|\M T\|_{\max}}\max_{(i,j)\in[n]^2}\big|p[\M E_1]_{ij}\big|= \frac{1}{p\|\M T\|_{\max}}\max_{(i,j)\in[n]^2}\Big\{(1 - 2p)\big|[\M T]_{ij}
\big|,p\big|[\M T]_{ij}\big|\Big\} \leq 1.$$
Furthermore we also have
\bee\label{variance:derive}
\max_{i}\sum_{j=1}^n \frac{\E\big|p[\M E_1]_{ij}\big|^2}{\|\M T\|_{\max}^2} &=
\max_{i}\sum_{j=1}^n\frac{p \bigl\{[\M T]_{ij} - p[\M
  T]_{ij}\bigr\}^2 + (1 -p)\bigl\{ -p[\M T]_{ij}\bigr\}^2}{\|\M T\|_{\max}^2}
\\
&\leq \max_{i}\sum_{j=1}^n\bigl[\{p(1 - p)^2 + (1 - p)p^2\}\bigr]
\leq np.
\ee
By Remark 3.13 in \cite{bandeira2016sharp}, there exists a universal constant $c>0$ such that for  $t = (np)^{1/2}$,
\bee\nonumber
\p\Big[\frac{p}{\|\M T\|_{\max}}\|{\M E_1}\|\geq 4(np)^{1/2} + t\Big]\leq ne^{-t^2/c}
\ee
which immediately implies, 
$
\|\M E_1\|\precsim (n/p)^{1/2} \|\M T\|_{\max}
$ with probability at least $1 - \tfrac{1}{2}n^{-3}$.
On the other hand,
$
p^{-1}\|\M E_2\| \leq \|\mathcal{P}_{\M\Omega}(\M N)\| \precsim \sigma (n/p)^{1/2},
$
with probability at least $1 - \tfrac{1}{2}n^{-3}$,  see e.g., Lemma
13 in \cite{abbe2020entrywise}. 
\Eq\eqref{E:decom:mc} therefore implies
\bee\label{E:bound:mc}
\|\M E\| \precsim (\sigma+\|\M T\|_{\max}) (n/p)^{1/2}
\ee 
with probability at least $1 - n^{-3}$. We can thus choose
$E_n = C_{\text{MC}}(\sigma+\|\M T\|_{\max})(n/p)^{1/2}$
for some finite constant $ C_{\text{MC}} > 2$. In summary we have
\begin{equation}
  \label{eq:SNR_proof_matrix_completion}
\frac{E_n}{|\lambda_{k_0}(\M T)|} \precsim \frac{n^{1/2}(\sigma + \|\M
  T\|_{\max})}{p^{1/2}|\lambda_{k_0}(\M T)|} \precsim \kappa^{-1} (\log n)^{-1/2} 
\end{equation}
where the last inequality follows from the assumption in Eq.~\eqref{eq:cond_matrix_comp1}.
\subsubsection{Bounding $d_2(\hat{\mathbf{U}}_{g},  \mathbf{U})$, 
   $d_{2\to\infty}(\hat{\mathbf{U}}_g, \mathbf{U})$ and $\|\hat{\mathbf{T}}_g - \mathbf{T}\|_{\text{F}}$}
Eq.~\eqref{eq:SNR_proof_matrix_completion} implies $|\lambda_{k_0}(\M T)|/E_n \rightarrow \infty$ as $n \rightarrow \infty$. Let $g \geq \frac{\log n}{\log(|\lambda_{k_0}(\M T)|/E_n)}$. Then by Corollary~\ref{co:l2_noise} and Corollary~\ref{co:l2inf_noise} we have
\bee\label{two:twoinf:mc}
d_{2}(\hat{\M U}_g,\M U) &\leq d_2(\hat{\M U}, \M U) + o\Bigl(\frac{E_n}{|\lambda_{k_0}(\M T)|}\Bigr) \\
d_{\twoinf}(\hat{\M U}_g,\M U) &\leq d_{\twoinf}(\hat{\M U}, \M U) + o\Bigl(\frac{n^{-1/2} E_n}{|\lambda_{k_0}(\M T)|}\Bigr)
\ee
with probability at least $1 - n^{-3}$.
Similarly, by Eq.~\eqref{eq:entrywise_guarantee}, we have
\begin{equation}
  \label{eq:entry_mc1}
  \|\hat{\M T}_g - \M T\|_{\max} \leq \|\hat{\M T}_S - \M T\|_{\max} + o\bigl(n^{-1/2} E_n \kappa \|\M U\|_{\twoinf}\bigr)
\end{equation}
with probability at least $1 - n^{-3}$, where $\hat{\M T}_S = p^{-1} \hat{\M T}^{(k_0)}$ is the truncated rank-$k_0$ SVD of $\hat{\M T}$.
Eq.~\eqref{eq:matrix_completion} follows directly from Eq.~\eqref{two:twoinf:mc} together with the bounds for
$d_{\twoinf}(\hat{\M U}_g, \M U)$ and $\|\hat{\M T}_S - \M T\|_{\max}$ given in Theorem~3.4 of \cite{abbe2020entrywise}. 
\subsubsection{Entrywise limiting distribution}
Fix $g \geq g_*:= \frac{\log n}{\log(|\lambda_{k_0}(\M T)/E_n)}$. 
For ease of exposition we say that an event
$\mathcal{E}$ happens with high probability (whp) if $\mathcal{E}$ happens
with probability at least $1 - Cn^{-3}$. Here $C > 0$ is an
arbitrary constant that can change from line to line.
First recall the definition of $v_{ij}^*$ in Eq.~\eqref{eq:def_vstar}. As the entries of $\M T$ are homogeneous,
the $\{[\M T^2]_{ii}\}$ also homogeneous. Now $[\M T^2]_{ii} = \|[\M U
\bm{\Lambda}]_{i}\|^2$ and hence (as $\M T$ has bounded condition number),
$\min_i\|[\M U]_i\|^2\asymp \max_i\|[\M U]_i\|^2$
where $\|[\M U]_i\|^2$ is the squared $\ell_2$ norm of the $i$th row of $\M U$. 
We therefore have
\bee\label{homog:con}
v_{ij}^* &\geq p^{-1} \{\min_{k \ell} (1-p) T_{k \ell}^2 +  \sigma^2\}(\|[\M U]_{i}\|^2 + \|[\M U]_j\|^2)
\\ & \succsim p^{-1} \{\min_{k \ell} (1-p) T_{k \ell}^2 +  \sigma^2\} \|\M U\|_{\twoinf}^2 \asymp
p^{-1}(\|\M T\|_{\max}^2 + \sigma^2) \|\M U\|_{\twoinf}^2
\ee
where the final ``equality'' is due to the assumption that $p$ is bounded away from $1$. 
Recalling the expression for $E_n$ given after Eq.~\eqref{E:bound:mc}, we have
\begin{equation*}
  \begin{split}
 (v_{ij}^*)^{-1/2} \times o(n^{-1/2} E_n \kappa \|\M U\|_{\twoinf}) &= o(n^{-1/2} p^{1/2} (\|\M T\|_{\max}^2 + \sigma^2)^{-1/2} E_n)
\\ &= o(\|\M T\|_{\max}^2 + \sigma^2)^{-1/2} (\|\M T\|_{\max} + \sigma)) 
= o(1)
  \end{split}
\end{equation*}
Therefore, by Eq.~\eqref{eq:entry_mc1}, we obtain
\begin{equation}
  \label{eq:clt_mc_vstar}
(v_{ij}^*)^{-1/2}([\hat{\M T}_g - \M T]_{ij}) = (v_{ij}^*)^{-1/2}([\hat{\M T}_S - \M T]_{ij}) + o(1) \rightsquigarrow \mathcal{N}(0,1)
\end{equation}
where the convergence in disitribution of $[\hat{\M T}_S - \M T]_{ij}$ is precisely Theorem~4.12 in \cite{chen2021spectral}. 
\subsubsection{Confidence interval}
We now derive Eq.~\eqref{eq:mc_clt2}. This is equivalent to showing that
$\hat{v}_{ij}/v_{ij}^* \rightarrow 1$ in probability where
$\hat{v}_{ij}$ is defined in Eq.~\eqref{def:empv}. Our derivations will proceed in three steps.
\\

\noindent{\bf Step 1.} We first consider a truncated version of $\hat{\M T}$. Recall that the (upper triangular) entries of $\M N$, denoted $\eta_{ij}$, are iid $\mathcal{N}(0, \sigma^2)$. Let $\tilde{\eta}_{ij} := \eta_{ij}\mathbb{I}\{|\eta_{ij}|
\leq 5 \sigma (\log n)^{1/2}\}$ and let $\tilde{\M N}$ be the matrix
whose entries are the $\tilde{\eta}_{ij}$. Now define
\bee\label{truncation:ET}
\tilde{\M T} := \mathcal{P}_{\M \Omega}(\M T + \tilde{\M N}) \text{ and }\tilde{\M E} = \tilde{\M T}/p - \M T.
\ee
By standard tail bounds for Gaussian distribution, we have 
\bee\nonumber
\p\big(\max_{ij} |\eta_{ij}| \leq 5\sigma (\log n)^{1/2} \big)  \geq 1 - n^2\p\big(|\eta_{ij}| > 5\sigma (\log n)^{1/2} \big) \geq 1 - n^{-10}.
\ee We thus have $\tilde{\M E} = \M E$ and $\tilde{\M T} = \hat{\M T}$ whp.
See \cite[Section~3.2.3]{chen2021spectral} for more details. 
Next let $v_{ij}$ be the variance of $[\M U\M U^\T \tilde{\M E}]_{ij} +
[\tilde{\M E}\M U\M U^\T ]_{ij}$. We can then use a similar argument to the proof of Theorem 4.12 in \cite{chen2021spectral}, to bound $v_{ij}$ from below. More specifically, for all $k, \ell$ we have  $$\zeta_{k \ell}^2 = ([\M
U]_k^{\top} [\M U]_{\ell})^2 \leq \|[\M U]_{k}\|^2 \times \|[\M
U]_{\ell}\|^2$$ and hence, 
following the proof of Lemma 4.19 in \cite{chen2021spectral}, we have
\bee\label{vij:def}
v_{ij} &=
\sum_{ \ell\neq j}\E\big([\tilde{\M E}]_{i \ell}^2\big)\zeta_{k \ell}^2 + \sum_{ \ell \neq i}\E\big([\tilde{\M E}]_{\ell
 j}^2\big)\zeta_{i \ell}^2  + \E\big([\tilde{\M
 E}]_{ij}^2\big)\big\{\zeta_{ii} +
\zeta_{jj}\big\}^2
\\
&\geq \Big\{\min_{(k, \ell)\in[n^2]}\E\big([\tilde{\M E}]_{k \ell}^2\big)\Big\}\Big\{\|[\M U]_i\|^2 + \|[\M U]_j\|^2\Big\}.
\ee
As the distribution of $\eta_{ij}$ is symmetric around $0$, the
distribution of $\tilde{\eta}_{ij}$ is also symmetric around $0$ and 
$\mathbb{E}[\tilde{\eta}_{ij}] = 0$. We therefore have
\bee\label{tE:res}
\E([\tilde{\M E}]_{ij}^2) &= p^{-1}\E\big\{([\M T]_{ij} - p[\M T]_{ij} + [\tilde{\M N}]_{ij})^2\big\} + (1-p)[\M T]_{ij}^2
\\
&=\frac{1-p}{p}[\M T]_{ij}^2 + \frac{1}{p}\E[\tilde{\eta}_{ij}^2]
\geq \frac{1-p}{p}\min_{k\ell}|T_{k\ell}|^2 + \frac{1}{2p}\sigma^2,
\ee
provided that $n$ is sufficiently large. The inequality in the above
display is derived as follows. By the Cauchy--Schwarz inequality we have
\bee
\label{tE:res:2}
\frac{1}{\sigma^2}\E\big[\eta_{ij}^2\mathbb{I}(|\eta_{ij}|> 5\sigma (\log n)^{1/2})\big]
\leq \Bigl(\E\bigl[(\eta_{ij}/\sigma)^4\bigr] \p\bigl[|\eta_{ij}/\sigma|> 5 (\log n)^{1/2} \bigr]\Bigr)^{1/2} = o(1),
\ee
and hence
\bee\nonumber
\E\big([\tilde{\M N}]_{ij}^2/\sigma^2\big) &= \sigma^{-2}\E\Big[\eta_{ij}^2\big\{1 - \mathbb{I}(|\eta_{ij}|> 5\sigma
(\log n)^{1/2} \big\}\Big]
= 1 + o(1).
\ee
Under the conditions of Theorem \ref{thm:mc} we have that $p$ is
bounded away from $1$ and $\M T$ is homogeneous, and thus
\bee\label{vij_bound}
v_{ij} \succsim p^{-1}(\|\M T\|_{\max}^2 + \sigma^2)\{\|[\M U]_i\|^2 + \|[\M U]_j\|^2\}
\ee
We next show that $v_{ij}/v_{ij}^* \rightarrow 1$ in probability
and hence $v_{ij}$ can be replaced by $v_{ij}^*$ in
without changing the limit result in Eq.~\eqref{eq:mc_clt1}; 
here $v_{ij}^*$ is defined in \eqref{eq:def_vstar}.
By \Eq\eqref{tE:res} and \Eq\eqref{tE:res:2}, we have
\bee\nonumber
\big|\E [{\M E}]_{k \ell}^2 - \E [\tilde{\M E}]_{k \ell}^2\big|  &=
\Big|\frac{1}{p}\E[\eta_{kl}^2] - \frac{1}{p}\E[\tilde{\eta}_{k \ell}^2\big)\Big|
\leq \frac{\sigma^2}{p} \Bigl(\E\big[(\eta_{kl}/\sigma)^4\big] \p\big(|\eta_{kl}| > 5\sigma (\log n)^{1/2}
\big)\Bigr)^{1/2}
\ee
and thus
\bee\label{eemax}
\max_{(k, \ell)\in[n]^2}\big|\E [{\M E}]_{k \ell}^2 - \E [\tilde{\M
  E}]_{k \ell}^2\big| = o\big(p^{-1} \sigma^2).
\ee
Combining Eq.~\eqref{tE:res}, Eq.~\eqref{vij_bound} and Eq.~\eqref{eemax}, we obtain
\bee\nonumber
\Big|\frac{v_{ij}^* - v_{ij}}{v_{ij}}\Big| &\leq\Big\{ \max_{(k,
  \ell)\in[n]^2}\big|\E [{\M E}]_{k \ell}^2 - \E [\tilde{\M
  E}]_{k \ell}^2\big|\Big\}\cdot\frac{\Big[\sum_{ \ell \not =
    j}\zeta_{\ell j}^2 + \sum_{\ell \not = i}\zeta_{i \ell}^2  +
  \big\{\zeta_{ii} + \zeta_{jj}\big\}^2\Big]}{v_{ij}}
\\
&\leq \frac{2 \max_{(k, \ell)\in[n]^2}\big|\E [{\M E}]_{kl}^2 - \E
  [\tilde{\M E}]_{k \ell}^2\big|\cdot\Big\{\|[\M U]_i\|^2 + \|[\M U]_j\|^2
  \Big\}}{\min_{(k, \ell)\in[n^2]}\E\big([\tilde{\M E}]_{kl}^2\big)\cdot\Big\{\|[\M U]_i\|^2 + \|[\M U]_j\|^2\Big\}}
= o(1),
\ee
as desired. 
\noindent{\bf Step 2.} We now consider an unbiased estimator for $v_{ij}$, namely
\bee\nonumber
\tilde{v}_{ij} &:= \sum_{\ell \neq j}[\tilde{\M E}]_{i
  \ell}^2\zeta_{\ell j}^2 + \sum_{\ell \not = i}[\tilde{\M E}]_{lj}^2
\zeta_{i \ell}^2  + [\tilde{\M E}]_{ij}^2(\zeta_{ii} + \zeta_{jj})^2.
\ee
Then by a standard application of Bernstein's inequality, we have $|\tilde{v}_{ij} - v_{ij}| = o(v_{ij}^*)$ with high probability. See the derivation of Eq.~(4.171) in \cite{chen2021spectral} for more details. 
Next let $
\breve{\M E} := \hat{\M U}_g\hat{\M U}_g^\T\tilde{\M T}/p - p^{-1} \hat{\M T}
$
where $\tilde{\M T}$ is defined in Eq.~\eqref{truncation:ET} and let
\bee\nonumber
\breve{v}_{ij} = \sum_{\ell \not = j}[\breve{\M E}]_{il}^2\zeta_{\ell
  j}^2 + \sum_{\ell \not = i}[\breve{\M E}]_{\ell j}^2\zeta_{i \ell}^2  +
[\breve{\M E}]_{ij}^2(\zeta_{ii} + \zeta_{jj})^2.
\ee
Then whp $\breve{\M E} = \hat{\M E} = \hat{\M T}_g -
p^{-1} \hat{\M T}$ and $\breve{v}_{ij} =
\hat{v}_{ij}$.
As $\tilde{\M E} = \M E$ whp, we have by Eq.~\eqref{eq:matrix_completion} that
\bee\nonumber
\|\breve{\M E} - \tilde{\M E}\|_{\max} = \|\hat{\M E} - \M E\|_{\max} = \|\hat{\M T}_g - \M T\| \precsim p^{-1/2} (\|\M T\|_{\max} + \sigma) (n \log n)^{1/2} \|\M U\|_{\twoinf}^2
\ee
whp. Therefore, by standard tail bounds for Gaussian random variables, we have
\bee\nonumber
\|\breve{\M E}\|_{\max} &\precsim \|\tilde{\M E}\|_{\max} + \|\breve{\M E} - \tilde{\M E}\|_{\max}
\\
&\precsim p^{-1/2}  {(\|\M T\|_{\max} + \sigma)} (n \log n)^{1/2}\|\M U\|_{\twoinf} +
p^{-1}\Big\{\|\M T\|_{\max} + \sigma\sqrt{\log n}\Big\}
\\
&\precsim \frac{1}{p}\Big\{\|\M T\|_{\max} + \sigma\sqrt{\log n}\Big\},
\ee
whp. Furthermore, we have by \Eq\eqref{two:twoinf:mc} that
\bee\nonumber
\|\hat{\M U}_g\|_{\twoinf} \leq \|\M U\|_{\twoinf} +
d_{\twoinf}(\hat{\M U}_g,\M U) = (1 + o(1)) \|\M U\|_{\twoinf})
\ee
whp and hence
\bee\nonumber
\|\hat{\M U}_g\hat{\M U}_g^\T - \M U\M U^\T\|_{\max} \leq \|\hat{\M U}_g\hat{\M U}_g^\T - \hat{\M U}\hat{\M U}^\T\|_{\max} + \|\hat{\M U} \hat{\M U}^{\top} - \M U \M U^{\top}\|_{\max} = (1 + o(1)) \|\hat{\M U} \hat{\M U}^{\top} - \M U \M U^{\top}\|_{\max}
\ee
whp. 
In summary we have whp that
\bee\nonumber
\|\hat{\M U}_g\hat{\M U}_g^\T\|_{\max} \leq \|\hat{\M U}_g\hat{\M U}_g^\T - \M U\M U^\T\|_{\max} + \|\M U\M U^\T\|_{\max} \precsim \|\M U\|_{\twoinf}^2.
\ee
\noindent{\bf Step 3.}
Finally we bound $|\hat{v}_{ij} - v_{ij}^*|$
using the same arguments as that presented in \Eq(4.176)--\Eq(4.177) of
\cite{chen2021spectral}, but with terms depending on $\hat{\M U}$
replaced by terms depending on $\hat{\M
  U}_g$.
More specifically let $\hat{\zeta}_{k \ell} = [\hat{\M U}_g \hat{\M
    U}_g]^{\top}_{k \ell}$. We then have, after some tedious algebra, that
\bee\label{vij2}
|\hat{v}_{ij} - \tilde{v}_{ij}| &= |\breve{v}_{ij} - \tilde{v}_{ij}|
\\
&\leq\Big| \sum_{\ell \not = j}\bigl([\breve{\M E}]_{i \ell}^2 \hat{\zeta}_{\ell j}^2
 - [\tilde{\M E}]_{i \ell}^2 \zeta_{\ell j}^2\bigr)
\Big| + \Big|\sum_{\ell \not = i}\bigl([\breve{\M
    E}]_{\ell j}^2 \hat{\zeta}_{i \ell}^2 - [\tilde{\M E}]_{\ell j}^2\zeta_{i \ell}^2  \bigr)\Big|
+ 2\Big| [\breve{\M E}]_{ij}^2\hat{\zeta}_{ii}\hat{\zeta}_{jj} - [\tilde{\M E}]_{ij}^2{\zeta}_{ii}{\zeta}_{jj}\Big|
\\
&\precsim(\|\tilde{\M E}\|_{\max} + \|\breve{\M
  E}\|_{\max})\|\breve{\M E} - \tilde{\M E}\|_{\max}\|\hat{\M U}_g\hat{\M U}_g^\T\|_{\max} 
\\
&+ (\|\hat{\M U}_g\hat{\M U}_g^\T\|_{\max} + \|{\M U}{\M U}^\T\|_{\max})\|\hat{\M U}_g\hat{\M U}_g^\T - \M U\M U^\T\|_{\max}\|\tilde{\M E}\|^2
\\
&+(\|\tilde{\M E}\|_{\max} + \|\breve{\M E}\|_{\max})\|\breve{\M E} - \tilde{\M E}\|_{\max}\|\hat{\M U}_g\hat{\M U}_g^\T\|_{\max}^2
\\
&+\|\tilde{\M E}\|_{\max}^2(\|\hat{\M U}_g\hat{\M U}_g^\T\|_{\max} + \|{\M U}{\M U}^\T\|_{\max})\|\hat{\M U}_g\hat{\M U}_g^\T - \M U\M U^\T\|_{\max}
\\ & \precsim \|\tilde{\M E}\|_{\max}
\|\hat{\M T}_S - \M T\|_{\max}
\|\M U \M U^{\top}\|_{\max} + \|\M U \M U^{\top}\|_{\max} \|\hat{\M U} \hat{\M U}^{\top} - \M U \M U^{\top}\|_{\max} \|\tilde{\M E}\|^2 \\ &+ \|\tilde{\M E}\|_{\max} \|\M T_S - \M T\|_{\max} \|\M U \M U^{\top}\|_{\max}^2 
+ \|\tilde{\M E}\|_{\max}^2 \|\M U \M U^{\top}\|_{\max} \|\hat{\M U} \hat{\M U}^{\top} - \M U \M U^{\top}\|_{\max}
\ee
Now define
\begin{gather*}
  B = \|\tilde{\M E}\|_{\max} \precsim p^{-1}\{\|\M T\|_{\max} + \sigma \log^{1/2}{n}\}, \quad
  \sigma^2_{E} = \max_{k \ell} \mathbb{E}\bigl([\tilde{\M E}]_{k \ell}^2\bigr) \precsim p^{-1}(\|\M T\|_{\max}^2 + \sigma^2). 
\end{gather*}
Then, following the same steps as that for bounding
Eq.~(4.177) in \cite{chen2021spectral} (where we had assumed that $k_0 \precsim 1$ and $\kappa \precsim 1$,
which then implies $\|\M U\|_{\twoinf} \asymp n^{-1/2}$) we have
\begin{equation}
  \label{eq:vhat-vstar_ij}
  \begin{split}
  |\hat{v}_{ij} - v_{ij}| &\precsim \frac{B \sigma_E  (\log n)^{1/2}}{n^{3/2}}
  +  \frac{\sigma_E^3 (\log n)^{1/2}}{|\lambda_{k_0}(\M T)| n^{1/2}}
    \\ &  \precsim  \frac{(\log n) \{\|\M T\|_{\max} + \sigma\}^2 }{(np)^{3/2}}
         + \frac{(\log n)^{1/2} \{\|\M T\|_{\max} + \sigma\}^3}{n^{1/2}p^{3/2} |\lambda_{k_0}(\M T)|}
  \end{split}
\end{equation}
Combining Eq.~\eqref{eq:vhat-vstar_ij} and Eq.~\eqref{homog:con} we obtain
\begin{equation*}
  \frac{|\hat{v}_{ij} - v_{ij}|}{v_{ij}} \precsim \frac{\log n}{(n p)^{1/2}} +
  \frac{(n \log n)^{1/2} \{\|\M T\|_{\max} + \sigma\}}{p^{1/2} |\lambda_{k_0}(\M T)|} = o(1)
\end{equation*}
where the final equality follows from the assumptions C3 and C4 in Theorem~\ref{thm:mc}. 
Therefore, by Slutsky's Theorem and Eq.~\eqref{eq:clt_mc_vstar}, we have
\bee\nonumber
\frac{[\hat{\M T}_g - \M T]_{ij}}{\hat{v}_{ij}^{1/2}}  = \Bigl(\frac{v_{ij}}{\hat{v}_{ij}}\Bigr)^{1/2} \frac{[\hat{\M T}_g - \M T]_{ij}}{v_{ij}^{1/2}}\rightsquigarrow \mathcal{N}(0,1),
\ee
as desired.
\subsection{Proof of Theorem~\ref{thm5}}
Let $\mGs = \M B\M F(\M B\M F)^{\top}$. Then with probability at least $1 - d^{-10}$, we have
\bee\label{basic:pca}
\quad |\lambda_i(\mGs)|\asymp m\lambda_i \,\,
\text{for $i \in [k_0]$}.
\ee
and conditional on a given $\M Q^*$, with probability at least $1 - d^{-10}$,
\bee
\label{eq:En_pca}
\|\bds{\M Q} - \mGs\|\precsim m\lambda_{k_0} \mathscr{E}.
\ee
See Eq.~(33) and Eq.~(37) in the supplementary material of \cite{cai2021subspace} for derivations of the above bounds; note that the quantities $n, r$ and $\mathcal{E}_{\mathrm{ce}}$ in \cite{cai2021subspace} corresponds to the quantities
$m, k_0$ and $\mathcal{E}$ in the current paper.
Let $E_n = m \lambda_k \mathscr{E}$. Then, recalling the condition for $m$ in Eq.~\eqref{eq:sample_size}, we have
$\lambda_{k_0}(\M Q^*)/E_n \asymp \mathscr{E}^{-1} \succsim \log^{2}(m+d)$ with probability at least $1 - d^{-6}$.   
Therefore, by applying Corollary~\ref{co:l2_noise} and Corollary~\ref{co:l2inf_noise} with
$\M M = \M Q^*$, $\hat{\M M} = \M Q$, $\sk \geq (1 + c_{\mathrm{gap}}) \max\{k_0, \log n\}$ and
$g \geq g_* = \frac{\log d}{\log(1/\mathscr{E})}$, we have with probability at least $1 - d^{-5}$ that
\begin{equation}
\label{eq:proof_pca_1}  
d_{2}(\hat{\M U}_g, \M U) \leq d_{2}(\hat{\M U}, \M U) + \mathscr{E}, \quad \text{and} \quad
d_{\twoinf}(\hat{\M U}_g, \M U) \leq d_{\twoinf}(\hat{\M U}, \M U) + o(d^{-1/2} \mathscr{E})
\end{equation}
Eq.~\eqref{eq:PCA_missing} follows directly from Eq.~\eqref{eq:proof_pca_1} together with bounds for
$d_{2}(\hat{\M U}, \M U)$ and $d_{\twoinf}(\hat{\M U}, \M U)$ given in
Corollary~4.3 of \cite{cai2021subspace}. 
\qed
\end{document}